\documentclass[12pt,reqno]{amsart}

\catcode`\@=11

\thinlines
\newskip\Einheit \Einheit=0.6cm
\newcount\xcoord \newcount\ycoord
\newdimen\xdim \newdimen\ydim \newdimen\PfadD@cke \newdimen\Pfadd@cke
\def\PfadDicke#1{\PfadD@cke#1 \divide\PfadD@cke by2 \Pfadd@cke\PfadD@cke \multiply\PfadD@cke by2}
\long\def\LOOP#1\REPEAT{\def\BODY{#1}\ITERATE}
\def\ITERATE{\BODY \let\next\ITERATE \else\let\next\relax\fi \next}
\let\REPEAT=\fi
\def\Punkt{\hbox{\raise-2pt\hbox to0pt{\hss\scriptsize$\bullet$\hss}}}
\def\DuennPunkt(#1,#2){\unskip
  \raise#2 \Einheit\hbox to0pt{\hskip#1 \Einheit
          \raise-2.5pt\hbox to0pt{\hss\normalsize$\bullet$\hss}\hss}}
\def\NormalPunkt(#1,#2){\unskip
  \raise#2 \Einheit\hbox to0pt{\hskip#1 \Einheit
          \raise-3pt\hbox to0pt{\hss\large$\bullet$\hss}\hss}}
\def\DickPunkt(#1,#2){\unskip
  \raise#2 \Einheit\hbox to0pt{\hskip#1 \Einheit
          \raise-4pt\hbox to0pt{\hss\Large$\bullet$\hss}\hss}}
\def\Kreis(#1,#2){\unskip
  \raise#2 \Einheit\hbox to0pt{\hskip#1 \Einheit
          \raise-4pt\hbox to0pt{\hss\Large$\circ$\hss}\hss}}
\def\Diagonale(#1,#2)#3{\unskip\leavevmode
  \xcoord#1\relax \ycoord#2\relax
      \raise\ycoord \Einheit\hbox to0pt{\hskip\xcoord \Einheit
         \unitlength\Einheit
         \line(1,1){#3}\hss}}
\def\AntiDiagonale(#1,#2)#3{\unskip\leavevmode
  \xcoord#1\relax \ycoord#2\relax \advance\xcoord by -0.05\relax
      \raise\ycoord \Einheit\hbox to0pt{\hskip\xcoord \Einheit
         \unitlength\Einheit
         \line(1,-1){#3}\hss}}
\def\Pfad(#1,#2),#3\endPfad{\unskip\leavevmode
  \xcoord#1 \ycoord#2 \thicklines\ZeichnePfad#3\endPfad\thinlines}
\def\ZeichnePfad#1{\ifx#1\endPfad\let\next\relax
  \else\let\next\ZeichnePfad
    \ifnum#1=1
      \raise\ycoord \Einheit\hbox to0pt{\hskip\xcoord \Einheit
         \vrule height\Pfadd@cke width1 \Einheit depth\Pfadd@cke\hss}%
      \advance\xcoord by 1
    \else\ifnum#1=2
      \raise\ycoord \Einheit\hbox to0pt{\hskip\xcoord \Einheit
        \hbox{\hskip-\PfadD@cke\vrule height1 \Einheit width\PfadD@cke depth0pt}\hss}%
      \advance\ycoord by 1
    \else\ifnum#1=3
      \raise\ycoord \Einheit\hbox to0pt{\hskip\xcoord \Einheit
         \unitlength\Einheit
         \line(1,1){1}\hss}
      \advance\xcoord by 1
      \advance\ycoord by 1
    \else\ifnum#1=4
      \raise\ycoord \Einheit\hbox to0pt{\hskip\xcoord \Einheit
         \unitlength\Einheit
         \line(1,-1){1}\hss}
      \advance\xcoord by 1
      \advance\ycoord by -1
    \fi\fi\fi\fi
  \fi\next}
\def\hSSchritt{\leavevmode\raise-.4pt\hbox to0pt{\hss.\hss}\hskip.2\Einheit
  \raise-.4pt\hbox to0pt{\hss.\hss}\hskip.2\Einheit
  \raise-.4pt\hbox to0pt{\hss.\hss}\hskip.2\Einheit
  \raise-.4pt\hbox to0pt{\hss.\hss}\hskip.2\Einheit
  \raise-.4pt\hbox to0pt{\hss.\hss}\hskip.2\Einheit}
\def\vSSchritt{\vbox{\baselineskip.2\Einheit\lineskiplimit0pt
\hbox{.}\hbox{.}\hbox{.}\hbox{.}\hbox{.}}}
\def\DSSchritt{\leavevmode\raise-.4pt\hbox to0pt{%
  \hbox to0pt{\hss.\hss}\hskip.2\Einheit
  \raise.2\Einheit\hbox to0pt{\hss.\hss}\hskip.2\Einheit
  \raise.4\Einheit\hbox to0pt{\hss.\hss}\hskip.2\Einheit
  \raise.6\Einheit\hbox to0pt{\hss.\hss}\hskip.2\Einheit
  \raise.8\Einheit\hbox to0pt{\hss.\hss}\hss}}
\def\dSSchritt{\leavevmode\raise-.4pt\hbox to0pt{%
  \hbox to0pt{\hss.\hss}\hskip.2\Einheit
  \raise-.2\Einheit\hbox to0pt{\hss.\hss}\hskip.2\Einheit
  \raise-.4\Einheit\hbox to0pt{\hss.\hss}\hskip.2\Einheit
  \raise-.6\Einheit\hbox to0pt{\hss.\hss}\hskip.2\Einheit
  \raise-.8\Einheit\hbox to0pt{\hss.\hss}\hss}}
\def\SPfad(#1,#2),#3\endSPfad{\unskip\leavevmode
  \xcoord#1 \ycoord#2 \ZeichneSPfad#3\endSPfad}
\def\ZeichneSPfad#1{\ifx#1\endSPfad\let\next\relax
  \else\let\next\ZeichneSPfad
    \ifnum#1=1
      \raise\ycoord \Einheit\hbox to0pt{\hskip\xcoord \Einheit
         \hSSchritt\hss}%
      \advance\xcoord by 1
    \else\ifnum#1=2
      \raise\ycoord \Einheit\hbox to0pt{\hskip\xcoord \Einheit
        \hbox{\hskip-2pt \vSSchritt}\hss}%
      \advance\ycoord by 1
    \else\ifnum#1=3
      \raise\ycoord \Einheit\hbox to0pt{\hskip\xcoord \Einheit
         \DSSchritt\hss}
      \advance\xcoord by 1
      \advance\ycoord by 1
    \else\ifnum#1=4
      \raise\ycoord \Einheit\hbox to0pt{\hskip\xcoord \Einheit
         \dSSchritt\hss}
      \advance\xcoord by 1
      \advance\ycoord by -1
    \fi\fi\fi\fi
  \fi\next}
\def\Koordinatenachsen(#1,#2){\unskip
 \hbox to0pt{\hskip-.5pt\vrule height#2 \Einheit width.5pt depth1 \Einheit}%
 \hbox to0pt{\hskip-1 \Einheit \xcoord#1 \advance\xcoord by1
    \vrule height0.25pt width\xcoord \Einheit depth0.25pt\hss}}
\def\Koordinatenachsen(#1,#2)(#3,#4){\unskip
 \hbox to0pt{\hskip-.5pt \ycoord-#4 \advance\ycoord by1
    \vrule height#2 \Einheit width.5pt depth\ycoord \Einheit}%
 \hbox to0pt{\hskip-1 \Einheit \hskip#3\Einheit 
    \xcoord#1 \advance\xcoord by1 \advance\xcoord by-#3 
    \vrule height0.25pt width\xcoord \Einheit depth0.25pt\hss}}
\def\Gitter(#1,#2){\unskip \xcoord0 \ycoord0 \leavevmode
  \LOOP\ifnum\ycoord<#2
    \loop\ifnum\xcoord<#1
      \raise\ycoord \Einheit\hbox to0pt{\hskip\xcoord \Einheit\Punkt\hss}%
      \advance\xcoord by1
    \repeat
    \xcoord0
    \advance\ycoord by1
  \REPEAT}
\def\Gitter(#1,#2)(#3,#4){\unskip \xcoord#3 \ycoord#4 \leavevmode
  \LOOP\ifnum\ycoord<#2
    \loop\ifnum\xcoord<#1
      \raise\ycoord \Einheit\hbox to0pt{\hskip\xcoord \Einheit\Punkt\hss}%
      \advance\xcoord by1
    \repeat
    \xcoord#3
    \advance\ycoord by1
  \REPEAT}
\def\Label#1#2(#3,#4){\unskip \xdim#3 \Einheit \ydim#4 \Einheit
  \def\lo{\advance\xdim by-.5 \Einheit \advance\ydim by.5 \Einheit}%
  \def\llo{\advance\xdim by-.25cm \advance\ydim by.5 \Einheit}%
  \def\loo{\advance\xdim by-.5 \Einheit \advance\ydim by.25cm}%
  \def\o{\advance\ydim by.25cm}%
  \def\ro{\advance\xdim by.5 \Einheit \advance\ydim by.5 \Einheit}%
  \def\rro{\advance\xdim by.25cm \advance\ydim by.5 \Einheit}%
  \def\roo{\advance\xdim by.5 \Einheit \advance\ydim by.25cm}%
  \def\l{\advance\xdim by-.30cm}%
  \def\r{\advance\xdim by.30cm}%
  \def\lu{\advance\xdim by-.5 \Einheit \advance\ydim by-.6 \Einheit}%
  \def\llu{\advance\xdim by-.25cm \advance\ydim by-.6 \Einheit}%
  \def\luu{\advance\xdim by-.5 \Einheit \advance\ydim by-.30cm}%
  \def\u{\advance\ydim by-.30cm}%
  \def\ru{\advance\xdim by.5 \Einheit \advance\ydim by-.6 \Einheit}%
  \def\rru{\advance\xdim by.25cm \advance\ydim by-.6 \Einheit}%
  \def\ruu{\advance\xdim by.5 \Einheit \advance\ydim by-.30cm}%
  #1\raise\ydim\hbox to0pt{\hskip\xdim
     \vbox to0pt{\vss\hbox to0pt{\hss$#2$\hss}\vss}\hss}%
}
\catcode`\@=12

\numberwithin{equation}{section}

\newtheorem{Theorem}{Theorem}

\newtheorem{Proposition}[Theorem]{Proposition}
\newtheorem{Lemma}[Theorem]{Lemma}

\newtheorem{Conjecture}[Theorem]{Conjecture}

\theoremstyle{definition}

\theoremstyle{remark}

\catcode`\@=11
\newwrite\Seiten
\immediate\openout\Seiten detsurv.sei\relax

\def\machSeite#1{\@ifundefined{@p#1}%
{\write\Seiten{\string\edef\csname p#1\endcsname{\thepage}}%
\@namedef{@p#1}{}}%
{\write\Seiten{\string\edef\csname p#1\endcsname{\csname p#1\endcsname , 
\thepage}}}%
}

\def\@secnumfont{\bfseries}
\def\section{\@startsection{section}{1}%
  \z@{.7\linespacing\@plus\linespacing}{.5\linespacing}%
  {\normalfont\bfseries}}
\def\subsection{\@startsection{subsection}{2}%
  \z@{.5\linespacing\@plus.7\linespacing}{-.5em}%
  {\normalfont\scshape}}
\catcode`\@=12

{\obeylines
\gdef\MATH{\begingroup\parindent0pt\parskip0pt plus 0.25pt\obeylines%
        \def^^M{\vskip5pt}%
        \obeyspaces\tt\small}%
}
\def\goodbreakpoint{\par\penalty-5000%
         \vrule height10pt depth2pt width0pt\leavevmode}
\def\endMATH{\endgroup}
\def\MATHphi{\leavevmode
        \hbox to 0pt{\hbox to 5.24995pt{\hss$\phi$\hss}\hss}}
\def\MATHGamma{\leavevmode
        \hbox to 0pt{\hbox to 5.24995pt{\hss$\Gamma$\hss}\hss}}
\def\MATHpi{\leavevmode
        \hbox to 0pt{\hbox to 5.24995pt{\hss$\pi$\hss}\hss}}
\def\MATHinfty{\leavevmode
        \hbox to 0pt{\hbox to 5.24995pt{\hss$\infty$\hss}\hss}}
\def\MATHhStrich{\leavevmode
        \hbox to 0pt{\hbox to 5.24995pt{\vrule height4.5pt depth-3.5pt width5.24995pt}\hss}}
\def\MATHluEck{\leavevmode
        \hbox to 0pt{\hbox to 5.24995pt{\hskip2.12497pt
         \vrule height4.5pt depth1pt width1pt
         \vrule height4.5pt depth-3.5pt width2.12498pt}\hss}}
\def\MATHruEck{\leavevmode
        \hbox to 0pt{\hbox to 5.24995pt{%
         \vrule height4.5pt depth-3.5pt width2.12497pt
         \vrule height4.5pt depth1pt width1pt
         \hskip2.12498pt}\hss}}
\def\MATHloEck{\leavevmode
        \hbox to 0pt{\hbox to 5.24995pt{\hskip2.12497pt
         \vrule height9pt depth-3.5pt width1pt
         \vrule height4.5pt depth-3.5pt width2.12498pt}\hss}}
\def\MATHroEck{\leavevmode
        \hbox to 0pt{\hbox to 5.24995pt{%
         \vrule height4.5pt depth-3.5pt width2.12497pt
         \vrule height9pt depth-3.5pt width1pt
         \hskip2.12498pt}\hss}}
\def\MATHvStrich{\leavevmode
        \hbox to 0pt{\hbox to 5.24995pt{\hskip2.12497pt
         \vtop to 0pt{\hsize1pt\vss%
                \vrule height17pt depth6pt width1pt\vskip8pt\vss\par}%
         \hskip2.12498pt}\hss}}
\def\MATHtStueck{\leavevmode
        \hbox to 0pt{\hbox to 5.24995pt{%
         \vrule height4.5pt depth-3.5pt width2.12497pt
         \vrule height4.5pt depth2pt width1pt
         \vrule height4.5pt depth-3.5pt width2.12498pt}\hss}}
\def\MATHbackslash{\leavevmode
        \hbox to 0pt{\hbox to 5.24995pt{\hss$\backslash$\hss}\hss}}
\def\MATHlbrace{\leavevmode
        \hbox to 0pt{\hbox to 5.24995pt{\hss$\{$\hss}\hss}}
\def\MATHrbrace{\leavevmode
        \hbox to 0pt{\hbox to 5.24995pt{\hss$\}$\hss}\hss}}
\def\MATHkleiner{\leavevmode
        \hbox to 0pt{\hbox to 5.24995pt{\hss$\langle$\hss}\hss}}
\def\MATHgroesser{\leavevmode
        \hbox to 0pt{\hbox to 5.24995pt{\hss$\rangle$\hss}\hss}}
\def\MATHhoch{\leavevmode
        \hbox to 0pt{\hbox to 5.24995pt{\hss$^\land$\hss}\hss}}
\def\MATHtief{\leavevmode
        \hbox to 0pt{\hbox to 5.24995pt{\hss\vrule height0pt depth.8pt width3pt\hss}\hss}}

\def\al{\alpha}
\def\be{\beta}
\def\ga{\gamma}
\def\de{\delta}

\def\la{\lambda}

\def\si{\sigma}

\def\om{\omega}
\def\Ga{\Gamma}
\def\De{\Delta}

\def\La{\Lambda}

\def\R{{\mathbb R}}

\def\Z{{\mathbb Z}}

\def\today{\ifcase\month\or
 January\or February\or March\or April\or May\or June\or
 July\or August\or September\or October\or November\or December\fi
 \space\number\day, \number\year}

\def\({\left(}
\def\){\right)}
\def\[{\left[}
\def\]{\right]}

\def\Per{\operatorname{Per}}
\def\Tr{\operatorname{Tr}}

\def\inv{\operatorname{inv}}
\def\maj{\operatorname{maj}}
\def\des{\operatorname{des}}

\def\prodl{\prod\limits}

\def\3{\ss}
\let\vv\v
\def\v#1{{\vert #1\vert}}
\def\Tr{\operatorname{Tr}}
\def\po#1#2{(#1)_{#2}}
\def\fl#1{\left\lfloor#1\right\rfloor}
\def\cl#1{\left\lceil#1\right\rceil}

\def\complexi{{\sqrt{-1}}}
\def\complexikl{{(\sqrt{-1})}}
\def\bk{\operatorname{bk}}
\def\CT{\operatorname{CT}}
\def\NC{\operatorname{NC}}
\def\rank{\operatorname{rank}}
\def\match{\operatorname{NCmatch}}
\def\Pf{\operatorname{Pf}}

\begin{document}

\newbox\Adr
\setbox\Adr\vbox{
\centerline{Institut f\"ur Mathematik der Universit\"at Wien,}
\centerline{Strudlhofgasse 4, A-1090 Wien, Austria.}
\centerline{E-mail: {\tt\footnotesize kratt@pap.univie.ac.at}}
\centerline{WWW: \footnotesize\tt http://radon.mat.univie.ac.at/People/kratt}
}

\title{Advanced Determinant Calculus}

\author[C.~Krattenthaler]{C.~Krattenthaler$^\dagger$\\[18pt]\box\Adr}

\address{Institut f\"ur Mathematik der Universit\"at Wien,
Strudlhofgasse 4, A-1090 Wien, Austria.\newline
e-mail: KRATT@Pap.Univie.Ac.At\\
WWW: \tt http://radon.mat.univie.ac.at/People/kratt}

\dedicatory{Dedicated to the pioneer of determinant evaluations 
(among many other things),\\
George Andrews}

\thanks{$^\dagger$ Research partially supported by the Austrian
Science Foundation FWF, grants P12094-MAT and P13190-MAT}


\begin{abstract}
The purpose of this article is threefold. First, it provides the
reader with a few useful and efficient tools which should enable her/him to
evaluate nontrivial determinants for the case such a determinant
should appear in her/his research. Second, it lists a number of such
determinants that have been already evaluated, together with
explanations which tell in which
contexts they have appeared. Third, it points out references 
where further such determinant evaluations can be found.
\end{abstract}

\maketitle

\section{Introduction} 
Imagine, you are working on a problem. As things develop it 
turns out that, in order to solve your problem, you need to evaluate 
a certain determinant. Maybe your determinant is
\begin{equation} \label{eq:Cauchyspezial}
\det_{1\le i,j,\le n}\(\frac {1} {i+j}\),
\end{equation}
or
\begin{equation} \label{eq:MacMahon}
\det_{1\le i,j\le n}\(\binom {a+b}{a-i+j} \),
\end{equation}
or it is possibly
\begin{equation} \label{eq:MRR1}
\det_{0\le i,j\le n-1}\(\binom {\mu+i+j}{2i-j} \),
\end{equation}
or maybe
\begin{equation} \label{eq:pentagon1}
\det_{1\le i,j\le n}\(\binom {x+y+j}{x-i+2j}-\binom
{x+y+j}{x+i+2j}\).
\end{equation}
Honestly, which ideas would you have? (Just to tell you that I do not
ask for something impossible: Each of these four determinants can be
evaluated in ``closed form". If you want to see the solutions
immediately, plus information where these determinants come from,
then go to \eqref{eq:Cauchy}, \eqref{eq:MacMahon-conj}/\eqref{eq:PP2}, 
\eqref{eq:MRR2}/\eqref{eq:MRR},
respectively \eqref{eq:pentagon}.)

Okay, let us try some row and column manipulations. Indeed, although
it is not completely trivial (actually, it is quite a challenge), 
that would work for the first two determinants, \eqref{eq:Cauchyspezial} 
and \eqref{eq:MacMahon}, although I do not recommend that. 
However, I do not recommend {\em at all\/} that you try
this with the latter two determinants, \eqref{eq:MRR1} and
\eqref{eq:pentagon1}. I promise that you will fail. (The
determinant \eqref{eq:MRR1} does not look much more complicated than
\eqref{eq:MacMahon}. Yet, it is.)

So, what should we do instead?

Of course, let us look in the literature! Excellent idea. We may have
the problem of not knowing where to start looking. Good starting
points are certainly classics like 
\machSeite{MuirAB}\cite{MuirAB},
\machSeite{MuirAC}\cite{MuirAC},
\machSeite{MuirAD}\cite{MuirAD},
\machSeite{PascAA}\cite{PascAA}
and 
\machSeite{TurnAA}\cite{TurnAA}\footnote{Turnbull's book \cite{TurnAA}
does in fact contain rather lots of very general identities satisfied by
determinants, than determinant ``evaluations" in the
strict sense of the word. However, suitable specializations of these
general identities do also yield ``genuine" evaluations, see for
example Appendix~B. Since the value of this book may not be easy to
appreciate because of heavy notation, we refer the reader to
\cite{LeclAA} for a clarification of the notation and 
a clear presentation of many such identities.}.
This will lead to the first success, as \eqref{eq:Cauchyspezial} does
indeed turn up there (see 
\machSeite{MuirAB}\cite[vol.~III, p.~311]{MuirAB}). Yes, you will also
find evaluations for 
\eqref{eq:MacMahon} (see e.g\@. 
\machSeite{OstrAA}\cite{OstrAA}) and \eqref{eq:MRR1} (see
\machSeite{MiRRAD}\cite[Theorem~7]{MiRRAD}) in the existing literature. But at the time of the
writing you will not, to the best of my knowledge, find an evaluation
of \eqref{eq:pentagon1} in the literature.

The purpose of this article is threefold. First, I want to describe
a few useful and efficient tools which should enable you to
evaluate nontrivial determinants (see Section~\ref{sec2}). 
Second, I provide a list containing 
a number of such
determinants that have been already evaluated, together with
explanations which tell in which
contexts they have appeared (see Section~\ref{sec3}). 
Third, even if you should not find your
determinant in this list, I point out references 
where further such determinant evaluations can be found, maybe your
determinant is there.

{\em Most important} of all is that I want to convince you that, today,
$$\text{{\em Evaluating determinants is not\/} (okay: may not be) {\em
difficult!\/}}$$ 
When George Andrews, who must be rightly called 
{\em the pioneer of determinant evaluations}, in the seventies astounded the
combinatorial
community by his highly nontrivial determinant evaluations (solving
difficult enumeration problems on plane partitions), it {\em was} really
difficult. His method (see Section~\ref{sec:LU} for a description)
required a good ``guesser" and an excellent
``hypergeometer" (both of which he was and is). While at that time
especially to be the latter was quite a task, in the meantime 
{\em both\/} guessing and evaluating binomial and hypergeometric
sums has been largely {\em trivialized}, as both can be done (most of the time)
{\em completely automatically}. 
For guessing (see Appendix~A) 
this is due to tools
like  {\tt Superseeker}\footnote{\label{foot:Integer}the 
electronic version of the ``Encyclopedia of Integer Sequences" 
\machSeite{SlPlAA}%
\machSeite{SloaAA}%
\cite{SlPlAA,SloaAA}, 
written and developed by Neil Sloane and
Simon Plouffe; see {\tt
http://www.research.att.com/\~{}njas/sequences/ol.html}},
{\tt gfun} and {\tt Mgfun}\footnote{\label{foot:gfun}written by Bruno Salvy
and Paul Zimmermann, respectively Frederic Chyzak; available from {\tt
http://pauillac.inria.fr/algo/libraries/libraries.html}} 
\machSeite{SaZiAA}%
\machSeite{ChyzAA}%
\cite{SaZiAA,ChyzAA}, and {\tt
Rate}\footnote{\label{foot:Rate}written in {\sl Mathematica} 
by the author; available from {\tt
http://radon.mat.univie.ac.at/People/kratt}; the {\sl Maple} equivalent
{\tt GUESS} by Fran\c cois B\'eraud and Bruno Gauthier is available from 
{\tt http://www-igm.univ-mlv.fr/\~{}gauthier}} (which
is by far the most primitive of the three, but it is the most effective
in this context). For ``hypergeometrics" this is due to the
``WZ-machinery"\footnote{\label{foot:WZ}{\em Maple} 
implementations written by Doron
Zeilberger are available from {\tt
http://www.math.temple.edu/\~{}zeilberg}, {\sl Mathematica} implementations 
written by Peter Paule, Axel Riese, Markus Schorn, Kurt Wegschaider
are available 
from {\tt
http://www.risc.uni-linz.ac.at/research/combinat/risc/software}} 
(see 
\machSeite{PeWZAA}%
\machSeite{WiZeAC}%
\machSeite{ZeilAM}%
\machSeite{ZeilAN}%
\machSeite{ZeilAV}%
\cite{PeWZAA,WiZeAC,ZeilAM,ZeilAN,ZeilAV}). And even if you
should meet a case where the WZ-machinery should exhaust your
computer's capacity, then there are still computer algebra
packages like HYP and HYPQ\footnote{written in {\sl Mathematica} 
by the author; available from
{\tt http://radon.mat.univie.ac.at/People/kratt}}, or 
HYPERG\footnote{written in {\sl Maple} by Bruno Ghauthier; available from
{\tt http://www-igm.univ-mlv.fr/\~{}gauthier}}, which make {\em you} an
expert hypergeometer, as these packages comprise
large parts of the present hypergeometric knowledge, and, thus,
enable you to conveniently manipulate binomial and
hypergeometric series (which George
Andrews did largely by hand) on the computer. 
Moreover, as of today, there are a few new (perhaps just overlooked) insights
which make life easier in many cases. It is these which form large
parts of Section~\ref{sec2}.

So, if you see a determinant, don't be frightened, evaluate it
yourself!

\section{Methods for the evaluation of determinants} 
\label{sec2}

In this section I describe a few useful methods and theorems which
(may) help you to evaluate a determinant. As was mentioned already in
the Introduction, it is always possible that simple-minded things
like doing some {\em row} and/or {\em column operations}, or applying 
{\em Laplace expansion}
may produce an (usually inductive) evaluation of a determinant. Therefore, you
are of course advised to try such things first. What I am mainly addressing
here, though, is the case where that first, ``simple-minded" attempt
failed. (Clearly, there is no point in addressing row and column
operations, or Laplace expansion.)

Yet, we must of course start (in Section~\ref{sec:standard}) 
with some standard determinants, such as
the {\em Vandermonde determinant\/} or {\em Cauchy's double alternant}. 
These are of course well-known. 

In Section~\ref{sec:general} we
continue with some {\em general determinant evaluations} that
generalize the evaluation of the Vandermonde determinant,
which are however apparently not equally well-known, although they should
be. {\em In
fact, I claim that about 80~\% of the determinants that you meet in
``real life," and which can apparently be evaluated, are a special
case of just the very first of these} (Lemma~\ref{lem:Krat1}; see
in particular Theorem~\ref{thm:PP} and the subsequent remarks).
Moreover, as is demonstrated in Section~\ref{sec:general}, it is pure
routine to check whether a determinant is a special case of one of
these general determinants. Thus, it can be really considered as a
``method" to see if a determinant can be evaluated by one of the
theorems in Section~\ref{sec:general}.

The next method which I describe is the so-called {\em ``condensation
method"} (see Section~\ref{sec:cond}), a method which allows to
evaluate a determinant inductively (if the method works).

In Section~\ref{sec:ident}, a method, which I call the {\em ``identification
of factors" method}, is described. This method has been extremely
successful recently. It is based on a very simple idea, which comes
from one of the standard proofs of the Vandermonde determinant
evaluation (which is therefore described in
Section~\ref{sec:standard}).

The subject of Section~\ref{sec:diff} is a method
which is based on finding one
or more {\em differential or difference equations} for the matrix of which
the determinant is to be evaluated.

Section~\ref{sec:LU} contains a short description of George Andrews'
favourite method, which basically consists of {\em explicitly} doing the
{\em LU-factorization} of the matrix of which the determinant is to be
evaluated.

The remaining subsections in this section are conceived as a complement
to the preceding. In Section~\ref{sec:Hankel} a special type of
determinants is addressed, {\em Hankel determinants}. (These are
determinants of the form $\det_{1\le i,j\le n}(a_{i+j})$, and are
sometimes also called {\em persymmetric} or {\em Tur\'anian
determinants}.) As is explained
there, you should expect that a Hankel determinant evaluation is to be
found in the domain of {\em orthogonal polynomials} and {\em continued
fractions}. Eventually, in Section~\ref{sec:misc} a few further,
possibly useful results are exhibited.

Before we finally move into the subject, it must be pointed out that the
methods of determinant evaluation as presented here are ordered
according to the conditions a determinant must satisfy so that the
method can be applied to it, from ``stringent" to ``less stringent". I.~e.,
first come the methods which require that the matrix of which the
determinant is to be taken satisfies a lot of conditions (usually: it
contains a lot of parameters, at least, implicitly), and in the end
comes the method (LU-factorization) which requires nothing. In fact,
this order (of methods) is also the order in which I recommend that
you try them on your determinant.
That is, what I suggest is (and this is the rule I follow):

\begin{enumerate}
\item[(0)] First try some simple-minded things ({\em row} and {\em column
operations},
{\em Laplace expansion}). Do not waste too much time. If you
encounter a {\em Hankel-determinant\/} then see Section~\ref{sec:Hankel}.
\item[(1)] If that fails, check whether your determinant is a special case of one of
the {\em general determinants} in Sections~\ref{sec:general} (and
\ref{sec:standard}).
\item[(2)] If that fails, see if the {\em condensation method\/} (see
Section~\ref{sec:cond}) works. (If necessary, try
to introduce more parameters into your determinant.)
\item[(3)] If that fails, try the {\em ``identification of factors"} method (see
Section~\ref{sec:ident}). Alternatively, and in particular if your
matrix of which you want to find the determinant is the matrix
defining a system of differential or difference equations, try the 
{\em differential/difference equation method\/} of Section~\ref{sec:diff}.
(If necessary, try
to introduce a parameter into your determinant.)
\item[(4)] If that fails, try to work out the {\em LU-factorization} of
your determinant (see Section~\ref{sec:LU}).
\item[(5)] If all that fails, then we are really in trouble. Perhaps
you have to put more efforts into determinant manipulations (see
suggestion (0))?
Sometimes it is worthwile to interpret the matrix whose 
determinant you want to know as a
linear map and try to find a basis on which this map acts
triangularly, or even diagonally
(this requires that the eigenvalues of the matrix are ``nice";
see 
\machSeite{DiFrAA}%
\machSeite{DiFrAB}%
\machSeite{KraHAA}%
\machSeite{KrSlAA}%
\machSeite{WilRAA}%
\cite{DiFrAA,DiFrAB,KraHAA,KrSlAA,WilRAA} for examples where that worked).
Otherwise, maybe something from Sections~\ref{sec:misc} or \ref{sec3} helps? 
\end{enumerate}

A final remark: It was indicated that some of the methods
require that your determinant contains (more or less) parameters.
Therefore it is always a good idea to:
$$\text {\em Introduce more parameters into your determinant!}$$ 
(We address this in more detail in the last paragraph of
Section~\ref{sec:standard}.)
The more parameters you can play with, the more likely you will be
able to carry out the determinant evaluation. (Just to mention a few
examples: The condensation method needs, at least, two parameters. The
``identification of factors" method needs, at least, one parameter,
as well as the differential/difference equation method in
Section~\ref{sec:diff}.)

\subsection{A few standard determinants} \label{sec:standard}

Let us begin with a short proof of the {\em Vandermonde determinant
evaluation}
\begin{equation} \label{eq:Vandermonde}
\det_{1\le i,j\le n}\(X_i^{j-1}\)=\prod _{1\le i<j\le n}
^{}(X_j-X_i).
\end{equation}

Although the following proof is well-known, it makes still sense to quickly
go through it because,
by extracting the essence of it, we will be able to build a very
powerful method out of it (see Section~\ref{sec:ident}).

If $X_{i_1}=X_{i_2}$ with $i_1\ne i_2$, then the Vandermonde
determinant \eqref{eq:Vandermonde} certainly vanishes because in that
case two rows of the determinant are identical. Hence,
$(X_{i_1}-X_{i_2})$ divides the determinant as a polynomial 
in the $X_i$'s. But that means that the complete product $\prod
_{1\le i<j\le n} (X_j-X_i)$ (which is exactly the right-hand side
of \eqref{eq:Vandermonde}) must divide the determinant.

On the other hand, the determinant is a polynomial in the $X_i$'s of
degree at most $\binom n2$. Combined with the previous observation,
this implies that the determinant equals the right-hand side product
times, possibly, some constant. To compute the constant, compare
coefficients of $X_1^0X_2^1\cdots X_n^{n-1}$ on both sides of 
\eqref{eq:Vandermonde}. This completes the proof of
\eqref{eq:Vandermonde}.

\medskip
At this point, let us extract the essence of this proof as we will
come back to it in Section~\ref{sec:ident}. The basic steps are:
{\em 
\begin{enumerate}
\item Identification of factors
\item Determination of degree bound
\item Computation of the multiplicative constant.
\end{enumerate}
}

\medskip

An immediate generalization of the Vandermonde determinant evaluation 
is given by the proposition below. It can be proved in just the same
way as the above proof of the Vandermonde determinant evaluation
itself.

\begin{Proposition} \label{prop:Vandermonde-allg}
Let $X_1,X_2,\dots,X_n$ be
indeterminates. If $p_1,p_2,\dots,p_{n}$ are polynomials of the form
$p_j(x)=a_{j}x^{j-1}+{}$lower terms, then
\begin{equation} \label{eq:Vandermonde-allg}
\det_{1\le i,j\le n}\(p_j(X_i)\)=a_1a_2\cdots a_n\prod _{1\le i<j\le n}
^{}(X_j-X_i).
\end{equation}
\quad \quad \qed
\end{Proposition}

The following
variations of the Vandermonde determinant evaluation are equally easy
to prove.
\begin{Lemma}The following identities hold true:
\begin{equation} \label{eq:Cn}
\det_{1\le i,j\le n}(X_i^j-X_i^{-j})=(X_1\cdots X_n)^{-n}\prod
_{1\le i<j\le n} ^{}\big((X_i-X_j)(1-X_iX_j)\big)
\ \prod _{i=1} ^{n}(X_i^2-1),
\end{equation}
\vskip-5pt
\begin{multline} \label{eq:Bn}
\det_{1\le i,j\le n}(X_i^{j-1/2}-X_i^{-(j-1/2)})
\\=(X_1\cdots
X_n)^{-n+1/2}\prod
_{1\le i<j\le n} ^{}\big((X_i-X_j)(1-X_iX_j)\big)
\ \prod _{i=1} ^{n}(X_i-1),
\end{multline}
\vskip-5pt
\begin{equation} \label{eq:Dn}
\det_{1\le i,j\le n}(X_i^{j-1}+X_i^{-(j-1)})=2\cdot(X_1\cdots X_n)^{-n+1}\prod
_{1\le i<j\le n} ^{}\big((X_i-X_j)(1-X_iX_j)\big)
,
\end{equation}
\vskip-5pt
\begin{multline}
\det_{1\le i,j\le n}(X_i^{j-1/2}+X_i^{-(j-1/2)})\\=(X_1\cdots
X_n)^{-n+1/2}\prod
_{1\le i<j\le n} ^{}\big((X_i-X_j)(1-X_iX_j)\big)
\ \prod _{i=1} ^{n}(X_i+1).
\end{multline}
\quad \quad \qed
\end{Lemma}
We remark that the evaluations \eqref{eq:Cn}, \eqref{eq:Bn}, 
\eqref{eq:Dn} are basically the
Weyl denominator factorizations of types $C$, $B$, $D$, respectively
(cf\@. 
\machSeite{FuHaAA}%
\cite[Lemma~24.3, Ex.~A.52, Ex.~A.62, Ex.~A.66]{FuHaAA}).
For that reason they may be called the {\em
``symplectic"}, the {\em ``odd orthogonal"}, and the {\em ``even orthogonal"
Vandermonde determinant evaluation}, respectively.

If you encounter generalizations of such determinants of the form
$\det_{1\le i,j\le n}(x_i^{\la_j})$ or 
$\det_{1\le i,j\le n}(x_i^{\la_j}-x_i^{-\la_j})$, etc., then you
should be aware that what you encounter is basically {\em Schur
functions}, {\em characters for the symplectic groups}, 
or {\em characters for the orthogonal groups}
(consult 
\machSeite{FuHaAA}%
\machSeite{MacdAC}%
\machSeite{ProcAK}%
\cite{FuHaAA,MacdAC,ProcAK} for more information on these
matters; see in particular
\machSeite{MacdAC}%
\cite[Ch.~I, (3.1)]{MacdAC}, 
\machSeite{FuHaAA}%
\cite[p.~403, (A.4)]{FuHaAA},
\machSeite{FuHaAA}%
\cite[(24.18)]{FuHaAA},
\machSeite{FuHaAA}%
\cite[(24.40) + first paragraph on
p.~411]{FuHaAA}, 
\machSeite{ProcAK}%
\cite[Appendix~A2]{ProcAK}, 
\machSeite{FuHaAA}%
\cite[(24.28)]{FuHaAA}).
In this context, one has to also mention Okada's general results on
evaluations of determinants and Pfaffians (see Section~\ref{sec:misc} for
definition) in 
\machSeite{OkadAI}%
\cite[Sec.~4]{OkadAI} and 
\machSeite{OkKrAA}%
\cite[Sec.~5]{OkKrAA}.

Another standard determinant evaluation is the evaluation of
{\em Cauchy's double alternant\/} (see
\machSeite{MuirAB}%
\cite[vol.~III, p.~311]{MuirAB}),
\begin{equation} \label{eq:Cauchy}
\det_{1\leq i,j\leq n}\left(
		\frac{1}{X_i+Y_j} 
\right)  =
	\frac{
		\prod_{1\leq i<j\leq n}(X_i-X_j)(Y_i-Y_j)
	}{
		\prod_{1\leq i, j\leq n}(X_i+Y_j)
	}.
\end{equation}
Once you have seen the above proof of the Vandermonde determinant
evaluation, you will immediately know how to prove this determinant
evaluation. 

On setting $X_i=i$ and $Y_i=i$, $i=1,2,\dots,n$, in
\eqref{eq:Cauchy}, we
obtain the evaluation of our first determinant in the Introduction,
\eqref{eq:Cauchyspezial}. 
For the evaluation of
a mixture of Cauchy's double alternant and Vandermonde's determinant
see 
\machSeite{BaFoAA}%
\cite[Lemma~2]{BaFoAA}.

Whether or not you tried to evaluate 
\eqref{eq:Cauchyspezial} directly, here is an {\em important lesson to be
learned} (it was already mentioned earlier): 
To evaluate \eqref{eq:Cauchyspezial} directly is quite
difficult, whereas proving its generalization \eqref{eq:Cauchy} is
almost completely trivial. Therefore, it is always
a good idea to try to {\em introduce more parameters into your
determinant}.
(That is, in a way such that the more
general determinant still evaluates nicely.) More
parameters mean that you have more objects at your disposal to play
with. 

The most stupid way to introduce parameters is to just write $X_i$
instead of the row index $i$, or write $Y_j$ instead of the column
index $j$.\footnote{Other common examples of introducing more
parameters are: Given that the $(i,j)$-entry of your determinant is a
binomial such as $\binom {i+j}{2i-j}$, try $\binom {x+i+j}{2i-j}$
(that works; see \eqref{eq:MRR}), 
or even $\binom {x+y+i+j}{y+2i-j}$ (that does not work;
but see \eqref{eq:MacMahon}), or $\binom {x+i+j}{2i-j}+\binom
{y+i+j}{2i-j}$ (that works; see \eqref{eq:AB}, and consult
Lemma~\ref{lem:GoJa} and the remarks thereafter). However, sometimes
parameters have to be introduced in an unexpected way, see
\eqref{det:Bombieri}.
(The parameter $x$ was introduced into a
determinant of Bombieri, Hunt and van der Poorten, which is obtained
by setting $x=0$ in \eqref{det:Bombieri}.)} 
For the determinant \eqref{eq:Cauchyspezial} even both
simultaneously was possible. For the determinant \eqref{eq:MacMahon}
either of the two (but not both) would work. On the contrary, there
seems to be no nontrivial way to introduce more parameters in the
determinant \eqref{eq:pentagon1}. This is an
indication that the evaluation of this determinant is
in a different category of difficulty of evaluation.
(Also \eqref{eq:MRR1} belongs to this ``different category". It is
possible to introduce one more parameter, see \eqref{eq:AB}, but it
does not seem to be possible to introduce more.)

\subsection{A general determinant lemma, plus variations and
generalizations} \label{sec:general}
In this section I present an apparently not so well-known 
determinant evaluation that generalizes 
Vandermonde's determinant, and some companions. As Lascoux pointed
out to me, most of these determinant evaluations can be derived from the
evaluation of a certain determinant of minors of a given matrix 
due to Turnbull 
\machSeite{TurnAB}%
\cite[p.~505]{TurnAB}, see Appendix~B.
However, this (these) determinant
evaluation(s) {\em deserve(s) to be better known}. Apart from the
fact that there are numerous applications of it (them)
which I am aware of, my proof is that I
meet very often people who stumble across a special case of this (these)
determinant evaluation(s), and then have a hard time to actually do the
evaluation because, usually, their special case does not show the
hidden general structure which is lurking behind. On the
other hand, as I will demonstrate in a moment, 
if you know this (these) determinant evaluation(s) then it is a
matter completely mechanical in nature to see whether 
it (they) is (are) applicable to your
determinant or not. If one of them is applicable, you are immediately done.

The determinant evaluation of which I am talking is the 
determinant lemma from 
\machSeite{KratAM}%
\cite[Lemma~2.2]{KratAM} given below. Here, 
and in the following, empty products (like
$(X_i+A_n)(X_i+A_{n-1})\dotsb(X_i+A_{j+1})$ for $j=n$) equal 1 by convention.

\begin{Lemma} \label{lem:Krat1}
Let $X_1,\dots,X_n$, $A_2,\dots,A_n$, and $B_2,
\dots ,B_n$ be indeterminates. Then there holds
\begin{multline} \label{eq:Krat1}
\det_{1\le i,j\le n}\Big((X_i+A_n)(X_i+A_{n-1})\cdots(X_i+A_{j+1})
(X_i+B_j)(X_i+B_{j-1})\cdots (X_i+B_2)\Big)\\
\hskip2cm =\prod _{1\le i<j\le n} ^{}(X_i-X_j)\prod _{2\le i\le j\le n}
^{}(B_i-A_j). 
\end{multline}
\quad \quad \qed
\end{Lemma}

Once you have guessed such a formula, it is easily proved. In the
proof in 
\machSeite{KratAM}%
\cite{KratAM} the determinant is reduced to a determinant
of the form \eqref{eq:Vandermonde-allg} by suitable column
operations. Another proof, discovered by Amdeberhan (private
communication), is by condensation, see Section~\ref{sec:cond}.
For a derivation from the above mentioned evaluation of a determinant of
minors of a given matrix, due to Turnbull, see Appendix~B.

Now let us see what the value of this formula is, by checking if it
is of any use in the case of the second determinant in the
Introduction, \eqref{eq:MacMahon}. The recipe that you should follow
is:

{\em

\begin{enumerate}
\item Take as many factors out of rows and/or columns of your
determinant, so that all denominators are cleared.
\item Compare your result with the determinant in \eqref{eq:Krat1}.
If it matches, you have found the evaluation of your determinant.
\end{enumerate}
}
\noindent
Okay, let us do so:
\begin{align*} 
\det_{1\le i,j\le n}&\(\binom {a+b}{a-i+j} \)=
\prod _{i=1} ^{n}\frac {(a+b)!} {(a-i+n)!\,(b+i-1)!}\\
&\hskip2cm
\times
\det_{1\le i,j\le n}\big((a-i+n)(a-i+n-1)\cdots(a-i+j+1)\\
&\hskip4cm
\cdot(b+i-j+1)(b+i-j+2)\cdots(b+i-1)\big)\\
&\quad =(-1)^{\binom n2}\prod _{i=1} ^{n}\frac {(a+b)!}
{(a-i+n)!\,(b+i-1)!}\\
&\hskip2cm
\times
\det_{1\le i,j\le n}\big((i-a-n)(i-a-n+1)\cdots(i-a-j-1)\\
&\hskip4cm
\cdot(i+b-j+1)(i+b-j+2)\cdots(i+b-1)\big).
\end{align*}
Now compare with the determinant in \eqref{eq:Krat1}. Indeed, the
determinant in the last line is just the special case $X_i=i$,
$A_j=-a-j$, $B_j=b-j+1$. Thus, by \eqref{eq:Krat1}, we have a result
immediately. A particularly attractive way to write it is
displayed in \eqref{eq:MacMahon-conj}.

Applications of Lemma~\ref{lem:Krat1} are abundant, see
Theorem~\ref{thm:PP} and the remarks accompanying it.

\medskip
In 
\machSeite{KratAO}%
\cite[Lemma~7]{KratAO}, a determinant evaluation is given
which is closely related to Lemma~\ref{lem:Krat1}. It was used there
to establish enumeration results about {\em shifted plane partitions
of trapezoidal shape}. It is the first result in the lemma below. It
is ``tailored" for the use in the context of $q$-enumeration. 
For plain enumeration, one would use the second result. This is a
limit case of the first (replace $X_i$ by $q^{X_i}$, $A_j$ by
$-q^{-A_j}$ and $C$ by $q^C$ in \eqref{eq:Krat2}, divide both sides by
$(1-q)^{n(n-1)}$, and then let $q\to1$). 
\begin{Lemma}\label{lem:Krat2}
Let $X_1,X_2,\dots,X_n,A_2,\dots,A_n$ be
indeterminates. Then there hold
\begin{multline} \label{eq:Krat2}
\det_{1\le i,j\le n} \big(( {C}/ {X_i}+A_n)( {C}/ {X_i}+A_{n-1})\dotsb( {C}/
{X_i}+A_{j+1})\\
\cdot (X_i+A_n)(X_i+A_{n-1})\dotsb(X_i+A_{j+1})\big)
=\prod _{i=2} ^{n}A_i^{i-1}\prod _{1\le i<j\le n}
(X_i-X_j)(1-C/X_iX_j),
\end{multline}
and
\begin{multline} \label{eq:Krat2a}
\det_{1\le i,j\le n} \big(( {X_i}-A_n-C)( {X_i}-A_{n-1}-C)\dotsb(
{X_i}-A_{j+1}-C)\\
\cdot (X_i+A_n)(X_i+A_{n-1})\dotsb(X_i+A_{j+1})\big)
=\prod _{1\le i<j\le n} (X_j-X_i)(C-X_i-X_j).
\end{multline}
\quad \quad \qed
\end{Lemma}
(Both evaluations are in fact
special cases in disguise of \eqref{eq:Vandermonde-allg}. Indeed, the
$(i,j)$-entry of the determinant in \eqref{eq:Krat2} is a polynomial
in $X_i+C/X_i$, while the 
$(i,j)$-entry of the determinant in \eqref{eq:Krat2a} is a polynomial
in $X_i-C/2$, both of degree $n-j$.)

The standard application of Lemma~\ref{lem:Krat2} is given in
Theorem~\ref{thm:shifted}.

\medskip
In 
\machSeite{KratAP}%
\cite[Lemma~34]{KratAP}, a common generalization of Lemmas~\ref{lem:Krat1}
and \ref{lem:Krat2} was given. In order to have a convenient statement of 
this determinant evaluation, we define the {\em degree} of 
a {\em Laurent polynomial} $p(X)=\sum _{i=M}
^{N}a_ix^i$, $M,N\in \Z$, $a_i\in \R$ and $a_N\ne0$, to be 
$\deg p:=N$.
\begin{Lemma} \label{lem:Krat3}
Let $X_1,X_2,\dots,X_n,A_2,A_3,\dots,A_n,C$ be
indeterminates. If $p_0,p_1,\dots,p_{n-1}$ are Laurent polynomials with
$\deg p_j\le j$ and $p_j(C/X)=p_j(X)$ for $j=0,1,\dots,n-1$, then
\begin{multline} \label{eq:Krat3}
\det_{1\le i,j\le n}\big((X_i+A_n)(X_i+A_{n-1})\dotsb(X_i+A_{j+1})\\
\cdot
(C/X_i+A_n)(C/X_i+A_{n-1})\dotsb(C/X_i+A_{j+1})\cdot p_{j-1}(X_i)\big)\\
=\prod _{1\le i<j\le n} ^{}(X_i-X_j)(1-C/X_iX_j)\prod _{i=1}
^{n}A_i^{i-1}\prod _{i=1} ^{n}p_{i-1}(-A_i)\ .
\end{multline}
\quad \quad \qed
\end{Lemma}

Section~3 contains several determinant evaluations which are
implied by the above determinant lemma, see Theorems~\ref{thm:PP5},
\ref{thm:PP4} and \ref{thm:PP6}.

Lemma~\ref{lem:Krat1} does indeed come out of the above
Lemma~\ref{lem:Krat3} by setting $C=0$ and
$$p_{j}(X)=\prod _{k=1} ^{j}(B_{k+1}+X).$$
Obviously, Lemma~\ref{lem:Krat2} is the special case
$p_j\equiv 1$, $j=0,1,\dots,n-1$.
It is in fact worth stating the $C=0$ case of Lemma~\ref{lem:Krat3}
separately.

\begin{Lemma} \label{lem:Krat3a}
Let $X_1,X_2,\dots,X_n,A_2,A_3,\dots,A_n$ be
indeterminates. If $p_0,p_1,\dots,p_{n-1}$ are polynomials with
$\deg p_j\le j$ for $j=0,1,\dots,n-1$, then
\begin{multline} \label{eq:Krat3a}
\det_{1\le i,j\le n}\big((X_i+A_n)(X_i+A_{n-1})\dotsb(X_i+A_{j+1})
\cdot p_{j-1}(X_i)\big)\\
=\prod _{1\le i<j\le n} ^{}(X_i-X_j)\prod _{i=1} ^{n}p_{i-1}(-A_i)\ .
\end{multline}
\quad \quad \qed
\end{Lemma}

Again, Lemma~\ref{lem:Krat3} is tailored for applications in
$q$-enumeration. So, also here, it may be convenient to state the
according limit case that is suitable for plain enumeration
(and perhaps other applications).
\begin{Lemma} \label{lem:Krat5}
Let $X_1,X_2,\dots,X_n,A_2,A_3,\dots,A_n,C$ be
indeterminates. If $p_0,p_1,\dots,\break p_{n-1}$ are polynomials with
$\deg p_j\le 2j$ and $p_j(C-X)=p_j(X)$ for $j=0,1,\dots,n-1$, then
\begin{multline} \label{eq:Krat5}
\det_{1\le i,j\le n}\big((X_i+A_n)(X_i+A_{n-1})\dotsb(X_i+A_{j+1})\\
\cdot
(X_i-A_n-C)(X_i-A_{n-1}-C)\dotsb(X_i-A_{j+1}-C)\cdot p_{j-1}(X_i)\big)\\
=\prod _{1\le i<j\le n} ^{}(X_j-X_i)(C-X_i-X_j)
\prod _{i=1} ^{n}p_{i-1}(-A_i)\ .
\end{multline}
\quad \quad \qed
\end{Lemma}

In concluding, I want to mention that, now since more than ten years,
I have a different common generalization of Lemmas~\ref{lem:Krat1} and
\ref{lem:Krat2} (with some
overlap with Lemma~\ref{lem:Krat3}) in my
drawer, without ever having found use for it. Let us nevertheless
state it here; maybe it is exactly the key to the solution of
a problem of yours.

\begin{Lemma} \label{lem:Krat6}
Let $X_1,\dots,X_n$, $A_2,\dots,A_n$, $B_2,
\dots B_n$, $a_2,\dots,a_n$, $b_2,
\dots b_n$, and $C$ be indeterminates. Then there holds
\begin{multline} \label{eq:Krat6}
\det_{1\le i,j\le n}\(
\begin{cases} (X_i+A_n)\cdots(X_i+A_{j+1})
(C/X_i+A_n)\cdots(C/X_i+A_{j+1})\\
(X_i+B_j)\cdots (X_i+B_2)
(C/X_i+B_j)\cdots (C/X_i+B_2)&j<m\\
(X_i+a_n)\cdots(X_i+a_{j+1})
(C/X_i+a_n)\cdots(C/X_i+a_{j+1})\\
(X_i+b_j)\cdots (X_i+b_2)
(C/X_i+b_j)\cdots (C/X_i+b_2)&j\ge m
\end{cases}
\)\\
=\prod _{1\le i<j\le n} ^{}(X_i-X_j)(1-C/X_iX_j)
\prod _{2\le i\le j\le m-1} ^{}(B_i-A_j)(1-C/B_iA_j)\kern2cm\\
\times
\prod _{i=2} ^{m}\prod _{j=m} ^{n}(b_i-A_j)(1-C/b_iA_j)
\prod _{m+1\le i\le j\le n} ^{}(b_i-a_j)(1-C/b_ia_j)\\
\times
\prod _{i=2} ^{m}(A_i\cdots A_n)\prod _{i=m+1} ^{n}(a_i\cdots a_n) 
\prod _{i=2} ^{m-1}(B_2\cdots B_i)\prod _{i=m} ^{n}(b_2\cdots b_i). 
\end{multline}
\quad \quad \qed
\end{Lemma}

The limit case which goes with this determinant lemma is the
following. (There is some overlap with Lemma~\ref{lem:Krat5}.)
\begin{Lemma} \label{lem:Krat7}
Let $X_1,\dots,X_n$, $A_2,\dots,A_n$, $B_2,
\dots ,B_n$, $a_2,\dots,a_n$, $b_2,
\dots ,b_n$, and $C$ be indeterminates. Then there holds
\begin{multline} \label{eq:Krat7}
\det_{1\le i,j\le n}\(
\begin{cases} (X_i+A_n)\cdots(X_i+A_{j+1})
(X_i-A_n-C)\cdots(X_i-A_{j+1}-C)\\
(X_i+B_j)\cdots (X_i+B_2)
(X_i-B_j-C)\cdots (X_i-B_2-C)&j<m\\
(X_i+a_n)\cdots(X_i+a_{j+1})
(X_i-a_n-C)\cdots(X_i-a_{j+1}-C)\\
(X_i+b_j)\cdots (X_i+b_2)
(X_i-b_j-C)\cdots (X_i-b_2-C)&j\ge m
\end{cases}
\)\\
=\prod _{1\le i<j\le n} ^{}(X_i-X_j)(C-X_i-X_j)
\prod _{2\le i\le j\le m-1} ^{}(B_i-A_j)(B_i+A_j+C)\kern3cm\\
\times
\prod _{i=2} ^{m}\prod _{j=m} ^{n}(b_i-A_j)(b_i+A_j+C)
\prod _{m+1\le i\le j\le n} ^{}(b_i-a_j)(b_i+a_j+C). 
\end{multline}
\quad \quad \qed
\end{Lemma}

If you are looking for more determinant evaluations of such a general type,
then you may want to look at 
\machSeite{SchlAB}%
\cite[Lemmas~A.1 and A.11]{SchlAB} and
\machSeite{SchlAE}%
\cite[Lemma~A.1]{SchlAE}.

\subsection{The condensation method} \label{sec:cond}

This is Doron Zeilberger's favourite method. It allows (sometimes) to
establish an elegant, {\em effortless} inductive proof of a determinant
evaluation, in which the only task is to guess the result correctly.

The method is often attributed
to Charles Ludwig Dodgson 
\machSeite{DodgAA}%
\cite{DodgAA}, better known as Lewis Carroll.
However, the identity on which it is based seems to be actually due
to P.~Desnanot (see 
\machSeite{MuirAB}%
\cite[vol.~I, pp.~140--142]{MuirAB}; with the first rigorous proof
being probably due to Jacobi, see 
\machSeite{BresAO}%
\cite[Ch.~4]{BresAO} and
\machSeite{KnutAF}%
\cite[Sec.~3]{KnutAF}). This
identity is the following.

\begin{Proposition} \label{prop:cond}
Let $A$ be an $n\times n$ matrix. Denote the submatrix of $A$ in which
rows $i_1,i_2,\dots,i_k$ and columns $j_1,j_2,\dots,j_k$ are 
omitted by $A_{i_1,i_2,\dots,i_k}^{j_1,j_2,\dots,j_k}$. Then there holds
\begin{equation} \label{eq:cond}
\hfill \det A\cdot \det A_{1,n}^{1,n}=\det A_{1}^{1}\cdot \det A_n^n-
\det A_1^n\cdot \det A_n^1.
\end{equation}
\quad \quad  \qed
\end{Proposition}

So, what is the point of this identity? Suppose you are given a
family $(\det M_n)_{n\ge0}$ of determinants, $M_n$ being an
$n\times n$ matrix, $n=0,1,\dots$. Maybe
$M_n=M_n(a,b)$ is the matrix underlying the
determinant in \eqref{eq:MacMahon}. Suppose
further that you have already worked out a conjecture for the
evaluation of $\det M_n(a,b)$ (we did in fact already evaluate this
determinant in Section~\ref{sec:general}, but let us ignore that for
the moment),
\begin{equation} \label{eq:MacMahon-conj}
\det M_n(a,b):=\det_{1\le i,j\le n}\(\binom {a+b}{a-i+j}\)\overset{\text
{?}}{=}
\prod _{i=1} ^{n}\prod _{j=1} ^{a}\prod _{k=1} ^{b}\frac {i+j+k-1}
{i+j+k-2}.
\end{equation}
Then you have already proved your conjecture, once you observe that 
\begin{align} \notag
\big(M_n(a,b)\big)_n^n&=M_{n-1}(a,b), \\
\notag
\big(M_n(a,b)\big)_1^1&=M_{n-1}(a,b), \\
\notag
\big(M_n(a,b)\big)_n^1&=M_{n-1}(a+1,b-1), \\
\notag
\big(M_n(a,b)\big)_1^n&=M_{n-1}(a-1,b+1), \\
\label{eq:minors}
\big(M_n(a,b)\big)_{1,n}^{1,n}&=M_{n-2}(a,b).
\end{align}
For, because of \eqref{eq:minors}, Desnanot's
identity \eqref{eq:cond}, with $A=M_n(a,b)$, 
gives a recurrence which expresses $\det M_n(a,b)$ in terms of quantities
of the form $\det M_{n-1}(\,.\,)$ and $\det M_{n-2}(\,.\,)$.
So, it just remains to check the conjecture \eqref{eq:MacMahon-conj} 
for $n=0$ and $n=1$, and that the right-hand side of
\eqref{eq:MacMahon-conj} satisfies the same recurrence, because that completes
a perfect induction with respect to $n$.
(What we have described here is basically the contents of
\machSeite{ZeilBL}%
\cite{ZeilBL}. For a bijective proof of Proposition~\ref{prop:cond}
see 
\machSeite{ZeilBP}%
\cite{ZeilBP}.)

Amdeberhan (private communication) discovered that in fact the
determinant evaluation \eqref{eq:Krat1} itself (which we used to
evaluate the determinant \eqref{eq:MacMahon} for the first time) can
be proved by condensation. The reader will easily figure out the
details. Furthermore, the condensation method also proves the
determinant evaluations \eqref{eq:Krat} and \eqref{eq:qKrat}. (Also
this observation is due to Amdeberhan 
\machSeite{AmdeAD}%
\cite{AmdeAD}.) At another
place, condensation was used by Eisenk\"olbl 
\machSeite{EisTAA}%
\cite{EisTAA} in order
to establish a conjecture by Propp 
\machSeite{PropAA}%
\cite[Problem~3]{PropAA} 
about the enumeration of rhombus
tilings of a hexagon where some triangles along the border of the
hexagon are missing.

\medskip
The reader should observe that crucial for a successful application
of the method is the existence of (at least) {\em two parameters} (in
our example these are $a$ and $b$),
which help to still stay within the same family of matrices
when we take minors of our original
matrix (compare \eqref{eq:minors}). (See the last paragraph of
Section~\ref{sec:standard} for a few hints of how to introduce more
parameters into your determinant, in the case that you are short of
parameters.)
Obviously, aside from the fact that we
need at least two parameters, we can hope for a success of
condensation only if our determinant is of a special kind.

\subsection{The ``identification of factors" method}
\label{sec:ident}

This is the method that I find most convenient to work with, once you
encounter a determinant that is not amenable to an evaluation 
using the previous recipes. It is best to explain this method along with 
an example. So, let us
consider the determinant in \eqref{eq:MRR1}. Here it is, together
with its, at this point, unproven evaluation,
\begin{multline} \label{eq:MRR2}
\det_{0\le i,j\le n-1}\(\binom {\mu+i+j}{2i-j} \)\\
=(-1)^{\chi(n\equiv 3\mod 4)}2^{\binom {n-1}2}\prod _{i=1} ^{n-1}
\frac {\(\mu+i+1\)_{\fl{(i+1)/2}}\,
\(-\mu-3n+i+\frac {3} {2}\)_{\fl{i/2}}} {(i)_i},
\end{multline}
where $\chi(\mathcal A)=1$ if $\mathcal A$ is
true and $\chi(\mathcal A)=0$ otherwise, and
where the {\em shifted factorial\/} $(a)_k$ is defined by 
$(a)_k:=a(a+1)\cdots(a+k-1)$, $k\ge1$, and $(a)_0:=1$.

As was already said in the
Introduction, this determinant belongs to a different category of
difficulty of evaluation, so that nothing what was presented so far will
immediately work on that determinant.

Nevertheless, I claim that the procedure which we chose to evaluate the
Vandermonde determinant works also with the above determinant.
To wit:
{\em 
\begin{enumerate}
\item Identification of factors
\item Determination of degree bound
\item Computation of the multiplicative constant.
\end{enumerate}
}
You will say: `A moment please! The reason that this procedure worked
so smoothly for the Vandermonde determinant is that there are so many
(to be precise: $n$) variables at our disposal. On the contrary, the
determinant in \eqref{eq:MRR2} has exactly {\em one} (!) variable.'
Yet --- and this is the point that I want to make here --- it works,
{\em in spite of having just one variable at our disposal\/}!.

What we want to prove in the first step is that the right-hand side
of \eqref{eq:MRR2} divides the determinant. For example, we would
like to prove that $(\mu+n)$ divides the determinant (actually,
$(\mu+n)^{\fl{(n+1)/3}}$, we will come to that in a moment). Equivalently,
if we set $\mu=-n$ in the determinant, then it should vanish. How
could we prove that? Well, if it vanishes then there must be a linear
combination of the columns, or of the rows, that vanishes. So, let us
find such a linear combination of columns or rows. Equivalently, for
$\mu=-n$ we find a vector in the kernel of the matrix in \eqref{eq:MRR2}, 
respectively its transpose. More generally (and this addresses that
we actually want to prove that $(\mu+n)^{\fl{(n+1)/3}}$ divides the
determinant):

\medskip
{\em \leftskip1cm
\rightskip1cm
\noindent
For proving that $(\mu+n)^E$ divides the
determinant, we find $E$ linear independent vectors in the
kernel.
\par}
\medskip
\noindent
(For a formal justification that this does indeed suffice,
see Section~2 of 
\machSeite{KratBI}%
\cite{KratBI}, and in particular the Lemma in that
section.)

Okay, how is this done in practice? You go to your computer, crank
out these vectors in the kernel, for
$n=1,2,3,\dots$, and try to make a guess what they are in general.
To see how this works, let us do it in our example.
What the computer gives is the following (we are using {\sl
Mathematica} here):

\MATH
\goodbreakpoint%
In[1]:= V[2]
\goodbreakpoint%
Out[1]= %
\MATHlbrace 0, c[1]%
\MATHrbrace 
\goodbreakpoint%
In[2]:= V[3]
\goodbreakpoint%
Out[2]= %
\MATHlbrace 0, c[2], c[2]%
\MATHrbrace 
\goodbreakpoint%
In[3]:= V[4]
\goodbreakpoint%
Out[3]= %
\MATHlbrace 0, c[1], 2 c[1], c[1]%
\MATHrbrace 
\goodbreakpoint%
In[4]:= V[5]
\goodbreakpoint%
Out[4]= %
\MATHlbrace 0, c[1], 3 c[1], c[3], c[1]%
\MATHrbrace 
\goodbreakpoint%
In[5]:= V[6]
\goodbreakpoint%
Out[5]= %
\MATHlbrace 0, c[1], 4 c[1], 2 c[1] + c[4], c[4], c[1]%
\MATHrbrace 
\goodbreakpoint%
In[6]:= V[7]
\goodbreakpoint%
Out[6]= %
\MATHlbrace 0, c[1], 5 c[1], c[3], -10 c[1] + 2 c[3], -5 c[1] + c[3], c[1]%
\MATHrbrace 
\goodbreakpoint%
In[7]:= V[8]
\goodbreakpoint%
Out[7]= %
\MATHlbrace 0, c[1], 6 c[1], c[3], -25 c[1] + 3 c[3], c[5], -9 c[1] + c[3], c[1]%
\MATHrbrace 
\goodbreakpoint%
In[8]:= V[9]
\goodbreakpoint%
Out[8]= %
\MATHlbrace 0, c[1], 7 c[1], c[3], -49 c[1] + 4 c[3], 
 
\MATHgroesser      -28 c[1] + 2 c[3] + c[6], c[6], -14 c[1] + c[3], c[1]%
\MATHrbrace 
\goodbreakpoint%
In[9]:= V[10]
\goodbreakpoint%
Out[9]= %
\MATHlbrace 0, c[1], 8 c[1], c[3], -84 c[1] + 5 c[3], c[5], 
 
\MATHgroesser      196 c[1] - 10 c[3] + 2 c[5], 98 c[1] - 5 c[3] + c[5], -20 c[1] + c[3], 
 
\MATHgroesser      c[1]%
\MATHrbrace 
\goodbreakpoint%
In[10]:= V[11]
\goodbreakpoint%
Out[10]= %
\MATHlbrace 0, c[1], 9 c[1], c[3], -132 c[1] + 6 c[3], c[5], 
 
\MATHgroesser      648 c[1] - 25 c[3] + 3 c[5], c[7], 234 c[1] - 9 c[3] + c[5], 
 
\MATHgroesser      -27 c[1] + c[3], c[1]%
\MATHrbrace 

\endMATH
\noindent
Here, $V[n]$ is the generic vector (depending on the indeterminates $c[i]$) 
in the kernel of the matrix in \eqref{eq:MRR2} with $\mu=-n$. 
For convenience, let us denote this
matrix by $M_n$.

You do not have to stare at these data for long to see that, in
particular, 

\begin{enumerate}
\item [] the vector $(0,1)$ is in the kernel of $M_2$, 
\item [] the vector $(0,1,1)$ is in the kernel of $M_3$, 
\item [] the vector $(0,1,2,1)$ is in the kernel of $M_4$,
\item [] the vector $(0,1,3,3,1)$ is in the kernel of $M_5$ (set $c[1]=1$ and $c[3]=3$),
\item [] the vector $(0,1,4,6,4,1)$ is in the kernel of $M_6$ (set $c[1]=1$ and
$c[4]=4$), etc. 
\end{enumerate}
\noindent
Apparently, 
\begin{equation} \label{eq:vector}
\textstyle \(0,\binom {n-2}0,\binom {n-2}1,\binom {n-2}2, \dots ,
\binom {n-2}{n-2}\)
\end{equation}
is in the kernel of $M_n$. That was easy! But we need more linear
combinations. Take a closer look, and you will see that the pattern
persists (set $c[1]=0$ everywhere, etc.). It will take you no time to work
out a full-fledged conjecture for $\fl{(n+1)/3}$ linear independent
vectors in the kernel of $M_n$.

Of course, there remains something to be proved. We need to actually
prove that our guessed vectors are indeed in the kernel. E.g., in
order to prove that the vector \eqref{eq:vector} is in the kernel, we
need to verify that
$$\sum _{j=1} ^{n-1}\binom {n-2}{j-1}\binom {-n+i+j}{2i-j}=0$$
for $i=0,1,\dots,n-1$.
However, verifying binomial identities is pure routine today, by means of
Zeilberger's algorithm 
\machSeite{ZeilAM}%
\machSeite{ZeilAV}%
\cite{ZeilAM,ZeilAV} (see Footnote~\ref{foot:WZ} in the
Introduction). 

Next you perform the same game with the other factors of the
right-hand side product of \eqref{eq:MRR2}. This is not much more
difficult. (See Section~3 of 
\machSeite{KratBI}%
\cite{KratBI} for details. There,
slightly different vectors are used.)

Thus, we would have finished the first step, ``identification of
factors," of our plan: We have proved that the right-hand side of
\eqref{eq:MRR2} divides the determinant as a polynomial in $\mu$.

The second step, ``determination of degree bound," 
consists of determining the (maximal) degree in $\mu$ of determinant and
conjectured result. As is easily seen, this is $\binom n2$ in each
case.

The arguments thus far show that the determinant in \eqref{eq:MRR2}
must equal the right-hand side times, possibly, some constant. To
determine this constant in the third step,
``computation of the multiplicative constant," 
one compares coefficients of $x^{\binom n2}$
on both sides of \eqref{eq:MRR2}. This is an enjoyable
exercise. (Consult 
\machSeite{KratBI}%
\cite{KratBI} if you do not want to do it yourself.)
Further successful applications of this procedure can be found in 
\machSeite{CiEKAA}%
\machSeite{CiKrAC}%
\machSeite{EisTAB}%
\machSeite{KratBG}%
\machSeite{KratBD}%
\machSeite{KratBH}%
\machSeite{KrZeAA}%
\machSeite{KupeAD}%
\machSeite{PoorAB}%
\cite{CiEKAA,CiKrAC,EisTAB,KratBG,KratBD,KratBH,KrZeAA,KupeAD,PoorAB}.

\medskip
Having done that, let me point out that most of the individual
steps in this sort of
calculation can be done (almost) {\em automatically}. 
In detail, what did we do? We had to

\medskip
\begin{enumerate}
\item  {\em Guess} the result. (Indeed, without the result we could not
have got started.)
\item  {\em Guess} the vectors in the kernel.
\item  Establish a {\em binomial\/} ({\em hypergeometric}) {\em identity}.
\item  Determine a degree bound.
\item  Compute a particular value or coefficient in order to determine the
multiplicative constant.
\end{enumerate}
\medskip

As I explain in Appendix~A, {\em guessing can be largely
automatized}. It
was already mentioned in the Introduction that {\em proving binomial\/} ({\em
hypergeometric})
{\em identities can be done by the computer}, thanks to the {\em
``WZ-machinery"}
\machSeite{PeWZAA}%
\machSeite{WiZeAC}%
\machSeite{ZeilAM}%
\machSeite{ZeilAN}%
\machSeite{ZeilAV}%
\cite{PeWZAA,WiZeAC,ZeilAM,ZeilAN,ZeilAV} (see Footnote~\ref{foot:WZ}). 
Computing the degree bound
is (in most cases) so easy that no computer is needed. (You may use it
if you want.) It is only the determination of the multiplicative
constant (item 5 above) 
by means of a special evaluation of the determinant or the
evaluation of a special coefficient (in our example we determined the
coefficient of $\mu^{\binom n2}$) for which I am not able to offer a
recipe so that things could be carried out on a computer.

\medskip
The reader should notice that crucial for a successful application
of the method is the existence of (at least) {\em one parameter} (in
our example this is $\mu$) to be able to apply the polynomiality
arguments that are the ``engine" of the method. If there is no
parameter (such as in the determinant in
Conjecture~\ref{conj:Bombieri},
or in the determinant \eqref{eq:Okada} which would
solve the problem of $q$-enumerating totally symmetric plane
partitions), then we even cannot get
started. (See the last paragraph of
Section~\ref{sec:standard} for a few hints of how to introduce a
parameter into your determinant, in the case that you are short of a
parameter.)

\medskip
On the other hand, a significant advantage of the ``identification of
factors method" is that not only is it capable of proving evaluations of 
the form
$$\det(M)=\text {CLOSED FORM},$$
(where CLOSED FORM means a product/quotient of ``nice" factors, such
as \eqref{eq:MRR2} or \eqref{eq:MacMahon-conj}), but
also of proving evaluations of the form
\begin{equation} \label{eq:UGLY}
\det(M)=\text {(CLOSED FORM)}\times \text {(UGLY POLYNOMIAL)},
\end{equation}
where, of course, $M$ is a matrix containing (at least) one
parameter, $\mu$ say. Examples of such determinant evaluations
are \eqref{eq:TSSCPP1}, \eqref{eq:qTSSCPP1}, \eqref{eq:tsscpp5} 
or \eqref{eq:FuKr2}.
(The UGLY POLYNOMIAL in \eqref{eq:TSSCPP1}, \eqref{eq:qTSSCPP1}
and \eqref{eq:FuKr2} is
the respective sum on the right-hand side, which in neither case
can be simplified). 

How would one approach the proof of such an evaluation?
For one part, we already know. ``Identification of factors" 
enables us to show that (CLOSED FORM) divides $\det(M)$ as a
polynomial in $\mu$. Then, comparison of degrees in $\mu$ on both
sides of \eqref{eq:UGLY} yields that (UGLY POLYNOMIAL) is a (at this
point unknown)
polynomial in $\mu$ of some maximal degree, $m$ say. How can we
determine this polynomial? Nothing ``simpler" than that: We find $m+1$
values $e$ such that we are able to evaluate $\det(M)$ at $\mu=e$. If
we then set $\mu=e$ in \eqref{eq:UGLY} and solve for (UGLY
POLYNOMIAL), then we obtain evaluations of (UGLY POLYNOMIAL) at
$m+1$ different values of $\mu$. Clearly, this suffices to find (UGLY
POLYNOMIAL), e.g., by Lagrange interpolation. 

I put ``simpler" in
quotes, because it is here where the crux is: We may not be able to find
{\em enough} such special evaluations of $\det(M)$. 
In fact, you may object: `Why all these complications?
If we should be able to find $m+1$ special values of $\mu$ for
which we are able to evaluate $\det(M)$, then what prevents us from
evaluating $\det(M)$ as a whole, for generic $\mu$?' When I am talking
of evaluating $\det(M)$ for $\mu=e$, then what I have in mind is that
the evaluation of $\det(M)$ at $\mu=e$ is ``nice" (i.e., 
gives a ``closed form," with no ``ugly"
expression involved, such as in \eqref{eq:UGLY}), which is easier to
identify (that is, to guess; see Appendix~A) 
and in most cases easier to prove.
By experience, such evaluations are rare. Therefore, the above
described procedure will only work if the degree of (UGLY POLYNOMIAL)
is not too large. (If you are just a bit short of evaluations, then
finding other informations about (UGLY POLYNOMIAL), like the leading
coefficient, may help to overcome the problem.)

To demonstrate this procedure by going through a concrete example is 
beyond the scope of this article. We refer the reader to
\machSeite{CiKrAA}%
\machSeite{FiscAA}%
\machSeite{FuKrAC}%
\machSeite{FuKrAD}%
\machSeite{KratBG}%
\machSeite{KratBD}%
\cite{CiKrAA,FiscAA,FuKrAC,FuKrAD,KratBG,KratBD} 
for places where this procedure
was successfully used to solve difficult enumeration problems on
rhombus tilings, respectively prove a conjectured constant term
identity.

\subsection{A differential/difference equation method}
\label{sec:diff}

In this section I outline a method for the
evaluation of determinants, often used by
Vitaly Tarasov and Alexander Varchenko, which, as the preceding method, 
also requires (at least) one parameter. 

Suppose we are given a matrix $M=M(z)$, depending on the parameter
$z$, of which we want to compute the determinant. 
Furthermore, suppose we know that $M$ satisfies a differential
equation of the form 
\begin{equation} \label{eq:diff-eq}
\frac {d} {dz}M(z)=T(z)M(z),
\end{equation} 
where $T(z)$ is some other known matrix. Then, by elementary linear
algebra, we obtain a differential equation for the determinant,
\begin{equation} \label{eq:diff-det}
\frac {d} {dz}\det M(z) = \Tr(T(z))\cdot \det M(z),
\end{equation}
which is usually easy to solve. (In fact, the differential operator
in \eqref{eq:diff-eq} and \eqref{eq:diff-det} could be replaced by
{\em any} operator. In particular, we could replace $d/dz$ by the
{\em difference} operator with respect to $z$, in which case
\eqref{eq:diff-det} is usually easy to solve as well.)

Any method is best illustrated by an example. Let us try this method
on the determinant \eqref{eq:MacMahon}. Right, we did already
evaluate this determinant twice (see Sections~\ref{sec:general} and
\ref{sec:cond}), but let us pretend that we have forgotten all this.

Of course, application of the method to \eqref{eq:MacMahon} itself
does not seem to be extremely promising, because that would involve
the differentiation of binomial coefficients. So, let us first take
some factors out of the determinant (as we also did in
Section~\ref{sec:general}),
\begin{align*} 
\det_{1\le i,j\le n}&\(\binom {a+b}{a-i+j} \)=
\prod _{i=1} ^{n}\frac {(a+b)!} {(a-i+n)!\,(b+i-1)!}\\
&\hskip2cm
\times\det_{1\le i,j\le n}\big((a-i+n)(a-i+n-1)\cdots(a-i+j+1)\\
&\hskip4cm
\cdot(b+i-j+1)(b+i-j+2)\cdots(b+i-1)\big).
\end{align*}
Let us denote the matrix underlying the 
determinant on the right-hand side of this equation
by $M_n(a)$. In order to apply the above method, we have need for a
matrix $T_n(a)$ such that
\begin{equation} \label{eq:diff-Tn}
\frac {d} {da}M_n(a)=T_n(a)M_n(a).
\end{equation} 
Similar to the procedure of Section~\ref{sec:LU}, the best idea is to
go to the computer, crank out $T_n(a)$ for $n=1,2,3,4,\dots$, and,
out of the data, make a guess for $T_n(a)$. Indeed, it suffices that
I display $T_5(a)$,
\begin{multline*} 
\left(\matrix 
  \frac {1} {1 + a + b} + \frac {1} {2 + a + b} + \frac {1} {3 + a + b} + 
    \frac {1} {4 + a + b}&\frac {4} {4 + a + b}&
   -\frac {{6}} {3 + a + b} + \frac {6} {4 + a + b}
\\ 
0&   \frac {1} {1 + a + b} + \frac {1} {2 + a + b} + \frac {1} {3 + a + b}&
   \frac {3} {3 + a + b}\\
    0&0&\frac {1} {1 + a + b} + \frac {1} {2 + a + b}\\
    0&0&0\\ 
  0&0&0
\endmatrix\right.
\\
\left.\begin{matrix} 
   \frac {4} {2 + a + b} - \frac {8} {3 + a + b} + \frac {4} {4 + a + b}&
   -\frac {1} {1 + a + b} + \frac {3} {2 + a + b} - \frac {3} {3 + a + b} + 
    \frac {1} {4 + a + b}
\\
-\frac {{3}} {2 + a + b} + \frac {3} {3 + a + b}&
   \frac {1} {1 + a + b} - \frac {2} {2 + a + b} + \frac {1} {3 + a + b}
\\
\frac {2} {2 + a + b}&
   -\frac {1} {1 + a + b} + \frac {1} {2 + a + b}
\\
\frac {1} {1 + a + b}&\frac {1} {1 + a + b}
\\
0&0
\end{matrix}\right)
\end{multline*}
(in this display, the first line contains columns $1,2,3$ of $T_5(a)$,
while the second line contains the remaining columns),
so that you are forced to conclude that, apparently, it must be true
that
$$T_n(a)=\(\binom {n-i}{j-i}\sum _{k=0} ^{n-i-1}\binom {j-i-1}k
\frac {(-1)^k} {a+b+n-i-k}\)_{1\le i,j,\le n}.$$
That \eqref{eq:diff-Tn} holds with this choice of $T_n(a)$ is then
easy to verify. Consequently, by means of \eqref{eq:diff-det}, we
have 
$$\frac {d} {da}\det M_n(a)=\bigg(\sum _{\ell=1} ^{n-1}
\frac {n-\ell} {a+b+\ell}\bigg)\det M_n(a),$$
so that
\begin{equation} \label{eq:const}
M_n(a)=\text{constant}\cdot\prod _{\ell=1} ^{n-1}(a+b+\ell)^{n-\ell}.
\end{equation}
The constant is found to be $(-1)^{\binom n2}\prod _{\ell=0}
^{n-1}\ell!$, e.g., by dividing both sides of
\eqref{eq:const} by $a^{\binom n2}$, letting $a$ tend to infinity, and
applying \eqref{eq:Vandermonde-allg} to the remaining determinant.

More sophisticated applications of this method (actually, of a
version for {\em systems} of {\em difference} operators) can be found in
\machSeite{TaVaAA}%
\cite[Proof of Theorem~5.14]{TaVaAA} and 
\machSeite{TaVaAB}%
\cite[Proofs of Theorems~5.9, 5.10, 5.11]{TaVaAB}, in the context
of the {\em Knizhnik--Zamolodchikov equations}.

\subsection{LU-factorization} \label{sec:LU}

This is George Andrews' favourite method. Starting point is the
well-known fact (see 
\machSeite{GantAA}%
\cite[p.~33ff]{GantAA}) that, given a square matrix $M$, 
there exists, under suitable, not very stringent conditions (in
particular, these are satisfied if all top-left principal minors of
$M$ are nonzero),
a unique lower triangular matrix $L$ and a unique upper diagonal matrix
$U$, the latter with all entries along the diagonal equal to 1, such that
\begin{equation} \label{eq:LU-echt}
M=L\cdot U.
\end{equation}
This unique factorization of the matrix $M$ is known as the
{\em L}(ower triangular){\em U}(pper triangular)-{\em factorization of
$M$}, or as well as the {\em Gau{\ss} decomposition} of $M$.

Equivalently, for a square matrix $M$ (satisfying these conditions) 
there exists a 
unique lower triangular matrix $L$ and a unique upper triangular matrix
$U$, the latter with all entries along the diagonal equal to 1, such that
\begin{equation} \label{eq:LU1}
M\cdot U=L.
\end{equation}
Clearly, once you know $L$ and $U$, the determinant of $M$ is easily
computed, as it equals the product of the diagonal entries of $L$.

Now, let us suppose that we are 
given a family $(M_n)_{n\ge0}$ of matrices, where $M_n$ is an
$n\times n$ matrix, $n=0,1,\dots$, of which we want to
compute the determinant. Maybe $M_n$ is the determinant in
\eqref{eq:MRR1}. By the above, we know that (normally) there exist 
uniquely determined matrices $L_n$ and $U_n$, $n=0,1,\dots$, $L_n$
being lower triangular, $U_n$ being upper triangular with all
diagonal entries equal to 1, such that
\begin{equation} \label{eq:LU}
M_n\cdot U_n=L_n.
\end{equation}
However, we do not know what the matrices $L_n$ and $U_n$ are. What
George Andrews does is that he goes to his computer, cranks out $L_n$
and $U_n$ for $n=1,2,3,4,\dots$ (this just amounts to solving a
system of linear equations), and, out of the data, tries to guess
what the coefficients of the matrices $L_n$ and $U_n$ are. Once he
has worked out a guess, he somehow proves that his guessed matrices
$L_n$ and $U_n$ do indeed satisfy \eqref{eq:LU}. 

This program is carried out in 
\machSeite{AnStAA}%
\cite{AnStAA}
for the family of determinants in
\eqref{eq:MRR1}. As it turns out, guessing is really
easy, while the underlying hypergeometric identities which are needed
for the proof of \eqref{eq:LU} are (from a hypergeometric viewpoint)
quite interesting. 

For a demonstration of the method of LU-factorization,
we will content ourselves here with trying the method on the Vandermonde
determinant. That is, let $M_n$ be the determinant in
\eqref{eq:Vandermonde}. We go to the computer and crank out the
matrices $L_n$ and $U_n$ for small values of $n$. For the purpose of
guessing, it suffices that I just display the matrices $L_5$ and $U_5$.
They are

{
\scriptsize
\begin{multline} \notag
L_5=\left(\matrix
1 & 0 & 0 \\
1 & (X_2-X_1) & 0 \\
1 & (X_3-X_1) & ( X_3-X_1 )  ( X_3-X_2 )  \\
1 & (X_4-X_1) & ( X_4-X_1 )  ( X_4-X_2 )  \\
1 & (X_5-X_1) &  ( X_5-X_1 )  ( X_5-X_2 )  \endmatrix\right.\\
\left.\matrix 
0&0\\
0&0\\
0&0\\
( X_4-X_1 )  ( X_4-X_2
  )  ( X_4-X_3 )  & 0\\
 (
  X_5-X_1 )  ( X_5-X_2 )  ( X_5-X_3 )  & ( X_5-X_1 )  ( X_5-X_2
  )  ( X_5-X_3 )  ( X_5-X_4 )
\endmatrix
\right)
\end{multline}
}

\noindent
(in this display, the first line contains columns $1,2,3$ of $L_5$,
while the second line contains the remaining columns), and

{\scriptsize
\begin{equation} \notag
U_5=\left(\matrix
1 & -e_1(X_1) & \hphantom{-}e_2(X_1, X_2) & - e_3(X_1 ,X_2, X_3)   &
  \hphantom{-}e_4(X_1, X_2 ,X_3 ,X_4)\\
0 & 1 & -e_1(X_1, X_2) & \hphantom{-}e_2(X_1, X_2, X_3) & - e_3(X_1, X_2, X_3,X_4)\\
0 & 0 & 1
  & -e_1(X_1, X_2, X_3) & \hphantom{-}e_2(X_1, X_2, X_3,X_4)\\
0 & 0 & 0 & 1 & -e_1(X_1,X_2,X_3,X_4)\\
0 & 0 & 0 & 0 & 1
\endmatrix
\right),
\end{equation}
}

\noindent
where $e_m(X_1,X_2,\dots,X_s)=\sum _{1\le i_1< \dots< i_m\le s}
^{}X_{i_1}\cdots X_{i_m}$ denotes the {\em $m$-th elementary
symmetric function\/}.

Having seen that, it will not take you for long to guess that,
apparently, $L_n$ is given by
$$L_n=\bigg(\prod _{k=1} ^{j-1}(X_i-X_k)\bigg)_{1\le i,j\le n},$$
and that $U_n$ is given by
$$U_n=\big((-1)^{j-i}e_{j-i}(X_1,\dots,X_{j-1})\big)_{1\le i,j\le n},$$
where, of course, $e_{m}(X_1,\dots):=0$ if $m<0$.
That \eqref{eq:LU} holds with these choices of $L_n$ and $U_n$ is
easy to verify. Thus, the Vandermonde determinant equals the product
of diagonal entries of $L_n$, which is exactly the product on the
right-hand side of \eqref{eq:Vandermonde}.

Applications of LU-factorization are abundant in the
work of George Andrews
\machSeite{AndrAK}%
\machSeite{AndrAN}%
\machSeite{AndrAO}%
\machSeite{AndrAS}%
\machSeite{AndrAW}%
\machSeite{AnStAA}%
\cite{AndrAK,AndrAN,AndrAO,AndrAS,AndrAW,AnStAA}. All of them concern
solutions to difficult enumeration problems on various types of plane
partitions. To mention another example, Aomoto and Kato
\machSeite{AoKaAA}%
\cite[Theorem~3]{AoKaAA}
computed the LU-factorization of a matrix which arose in the
theory of $q$-difference equations, thus proving a conjecture by
Mimachi 
\machSeite{MimaAD}%
\cite{MimaAD}.

\medskip
Needless to say that this allows for variations. You may try to guess
\eqref{eq:LU-echt} directly (and not its variation \eqref{eq:LU1}),
or you may try to guess the {\em U}(pper triangular){\em L}(ower
triangular) {\em
factorization}, or its variation in the style of \eqref{eq:LU1}. I am
saying this because it may be easy to guess the form of one
of these variations, while it can be very difficult to guess the form of
another.

\medskip
It should be observed that the way LU-factorization is used here in
order to evaluate determinants is very much in the same spirit as
``identification of factors" as described in the previous section. In
both cases, the essential steps are to {\em first guess} something, and
{\em then prove} the guess. Therefore, the remarks from the previous
section about guessing and proving binomial (hypergeometric) identities
apply here as well. In particular, for guessing you are once more
referred to Appendix~A.

\medskip
It is important to note that, as opposed to ``condensation" or
``identification of factors," LU-factorization does not require any
parameter. So, in principle, {\em it is applicable to any
determinant} (which satisfies the aforementioned conditions). If
there are limitations, then, from my experience, it is that the
coefficients which have to be guessed in LU-factorization tend to be
more complicated than in ``identification of factors". That is, 
guessing \eqref{eq:LU} (or one of its variations) may sometimes be
not so easy.

\subsection{Hankel determinants} \label{sec:Hankel}

A {\em Hankel determinant\/} is a determinant of a matrix which has
constant entries along antidiagonals, i.e., it is a determinant of
the form
$$\det_{1\le i,j,\le n}(c_{i+j}).$$
If you encounter a Hankel determinant, which you think evaluates
nicely, then expect the evaluation of your Hankel determinant 
to be found within the domain of
{\em continued fractions} and {\em orthogonal polynomials}. In this
section I explain what this connection is.

To make things concrete, let us suppose that we want to evaluate
\begin{equation} \label{eq:Berneinf}
\det_{0\le i,j\le n-1}(B_{i+j+2}),
\end{equation}
where $B_k$ denotes the $k$-th {\em Bernoulli number}. (The Bernoulli
numbers are defined via their generating function, $\sum _{k=0}
^{\infty}B_kz^k/k!=z/(e^z-1)$.) You have to try hard if you want to find an  
evaluation of \eqref{eq:Berneinf} {\em explicitly} in the literature.
Indeed, you can find it, hidden in
Appendix~A.5 of 
\machSeite{MehtAB}%
\cite{MehtAB}. However, even if you are not able to discover this
reference (which I would not have as well, unless the author of 
\machSeite{MehtAB}%
\cite{MehtAB} would not have drawn my attention to it), 
there is a rather straight-forward way to find an
evaluation of \eqref{eq:Berneinf}, which I outline below. It is
based on the fact, and this is the main point of this section, that
evaluations of Hankel determinants like \eqref{eq:Berneinf}
are, at least {\em implicitly}, in the literature on the theory of orthogonal
polynomials and continued fractions, which is very accessible today.

\medskip
So, let us review the relevant facts about orthogonal
polynomials and continued fractions (see
\machSeite{JoThAA}%
\machSeite{KoSwAA}%
\machSeite{PerronKB}%
\machSeite{SzegoOP}%
\machSeite{VienAE}%
\machSeite{WallCF}%
\cite{JoThAA,KoSwAA,PerronKB,SzegoOP,VienAE,WallCF} for more
information on these topics). 

We begin by citing the result, due to Heilermann, 
which makes the connection between
Hankel determinants and continued fractions.
\begin{Theorem}
\label{cor:cfracHankel}
{\em (Cf\@. 
\machSeite{WallCF}%
\cite[Theorem 51.1]{WallCF} or
\machSeite{VienAE}%
\cite[Corollaire~6, (19), on p.~IV-17]{VienAE})}.
Let $(\mu_k)_{k\ge0}$ be a sequence of numbers with generating
function $\sum_{k=0}^\infty{\mu_k}x^k$ written in the form
\begin{equation}
\label{eq:momentgf}
\sum_{k=0}^\infty{\mu_k}x^k=\cfrac{
	\mu_0}
		{1+a_0x-\cfrac{
			b_1x^2}
				{1+a_1x-\cfrac{
					b_2x^2}
						{1+a_2x-\cdots}}}\quad .
\end{equation}
Then the Hankel determinant $\det_{0\le i,j\le n-1}(\mu_{i+j})$
equals $\mu_0^nb_1^{n-1}b_2^{n-2}\cdots b_{n-2}^2b_{n-1}$.\quad \quad \qed
\end{Theorem}
(We remark that a continued fraction of the type as in
\eqref{eq:momentgf} is called a {\em J-fraction}.)

Okay, that means we would have evaluated \eqref{eq:Berneinf} once we
are able to explicitly expand the generating function $\sum _{k=0}
^{\infty}B_{k+2}x^k$ in terms of a continued fraction of the form of
the right-hand side of \eqref{eq:momentgf}. Using the tools explained
in Appendix~A, it is easy to work out a conjecture,
\begin{equation} \label{eq:Bern-cont}
\sum _{k=0} ^{\infty} B_{k+2}x^k=\cfrac{
	1/6}
		{1-\cfrac{
			b_1x^2}
				{1-\cfrac{
					b_2x^2}
						{1-\cdots}}}\quad ,
\end{equation}
where $b_i=-i(i+1)^2(i+2)/4(2i+1)(2i+3)$, $i=1,2,\dots$.
If we would find this expansion in the literature then we would be
done. But if not (which is the case here),
how to prove such an expansion? The key is {\em orthogonal
polynomials}.

A sequence $(p_n(x))_{n\ge0}$ of polynomials
is called (formally) {\em orthogonal\/} if $p_n(x)$ has degree $n$,
$n=0,1,\dots$, and if there
exists a linear functional $L$ such that
$L(p_n(x)p_m(x))=\delta_{mn}c_n$ for some sequence $(c_n)_{n\ge0}$ of
nonzero numbers, with $\delta_{m,n}$ denoting the Kronecker delta (i.e.,
$\delta_{m,n}=1$ if $m=n$ and $\delta_{m,n}=0$ otherwise).

The first important theorem in the theory of orthogonal polynomials
is Favard's Theorem, which gives an unexpected characterization for
sequences of orthogonal polynomials, in that it completely avoids the
mention of the functional $L$.

\begin{Theorem} \label{thm:Favard}
{\em (Cf\@.  
\machSeite{VienAE}%
\cite[Th\'eor\`eme~9 on
p.~I-4]{VienAE} or 
\machSeite{WallCF}%
\cite[Theorem 50.1]{WallCF}).}
Let $(p_n(x))_{n\ge0}$ be a sequence of monic polynomials, the
polynomial $p_n(x)$ having degree $n$, $n=0,1,\dots$. 
Then the sequence $(p_n(x))$ is
(formally) orthogonal if and only if
there exist sequences $(a_n)_{n\ge1}$ and $(b_n)_{n\ge1}$, with
$b_n\ne0$ for all $n\ge1$, such that the three-term recurrence
\begin{equation}
p_{n+1}(x)=(a_{n}+x)p_{n}(x)-b_{n}p_{n-1}(x),\quad \quad 
\text{ for } n\geq 1,
\label{eq:three-term}
\end{equation}
holds, with initial conditions $p_0(x)=1$ and $p_1(x)=x+a_0$.\quad \quad \qed
\end{Theorem}

What is the connection between orthogonal polynomials and continued
fractions? This question is answered by the next theorem,
the link being the {\em generating function of the moments}.
\begin{Theorem} \label{thm:momentgf}
{\em (Cf\@. 
\machSeite{WallCF}%
\cite[Theorem 51.1]{WallCF} or 
\machSeite{VienAE}%
\cite[Proposition~1, (7), on p.~V-5]{VienAE})}.
Let $(p_n(x))_{n\ge0}$ be a sequence of monic polynomials, the
polynomial $p_n(x)$ having degree $n$, which is orthogonal with
respect to some functional $L$. Let
\begin{equation}
p_{n+1}(x)=(a_{n}+x)p_{n}(x)-b_{n}p_{n-1}(x)
\label{eq:three-term2}
\end{equation}
be the corresponding three-term recurrence which is guaranteed by
Favard's theorem. Then the generating function $\sum _{k=0}
^{\infty}\mu_kx^k$ for the moments
$\mu_k=L(x^k)$ satisfies \eqref{eq:momentgf} with the $a_i$'s and
$b_i$'s being the coefficients in the three-term recurrence
\eqref{eq:three-term2}.\quad \quad \qed
\end{Theorem}

Thus, what we have to do is to find orthogonal polynomials
$(p_n(x))_{n\ge0}$, the three-term recurrence of which is explicitly known, 
and which are orthogonal with respect to some linear
functional $L$ whose moments $L(x^k)$ are exactly equal to $B_{k+2}$.
So, what would be very helpful at this point 
is some sort of table of orthogonal polynomials.
Indeed, there is such a {\em table for hypergeometric and basic
hypergeometric orthogonal polynomials}, proposed by Richard Askey
(therefore called the {\em``Askey table"}), and 
compiled by Koekoek and Swarttouw 
\machSeite{KoSwAA}%
\cite{KoSwAA}.

Indeed, in Section~1.4 of 
\machSeite{KoSwAA}%
\cite{KoSwAA}, we find
the family of orthogonal polynomials that is of relevance here, the
{\em continuous Hahn polynomials}, first studied by
Atakishiyev and Suslov 
\machSeite{AtSuAA}%
\cite{AtSuAA} and Askey 
\machSeite{AskeyCHP}%
\cite{AskeyCHP}. 
These polynomials depend on four parameters, $a,b,c,d$. It is
just the special choice $a=b=c=d=1$ which is of interest to us.
The theorem below lists the
relevant facts about these special polynomials.
\begin{Theorem} \label{thm:cHahn}
The continuous Hahn polynomials with parameters $a=b=c=d=1$, 
$(p_n(x))_{n\ge0}$, are
the monic polynomials defined by
\begin{equation}
p_n(x)
=\complexikl^n \frac {(n+1)!^2\,(n+2)!} {(2n+2)!}
\sum _{k=0} ^{\infty}\frac {(-n)_k\,(n+3)_k\,(1+x\complexi )_k}
{k!\,(k+1)!^2},
\end{equation}
with the shifted factorial $(a)_k$ defined as previously {\em(}see
{\em\eqref{eq:MRR2})}.
These polynomials satisfy the three-term recurrence
\begin{equation} \label{eq:cHahnrec}
p_{n+1}(x)  
=	xp_{n}(x)+\frac {n(n+1)^2(n+2)} {4(2n+1)(2n+3)}p_{n-1}(x).
\end{equation}
They are orthogonal with respect to the functional $L$ which is given
by
\begin{equation} \label{eq:cHahnortho}
L(p(x))=\frac{\pi}{2}\int_{-\infty}^{\infty}
\frac {x^2} {\sinh^2( \pi x)} p(x)\,dx\,.
\end{equation}
Explicitly, the orthogonality relation is
\begin{equation}
L(p_m(x)p_n(x)) 
=\frac{n!\,(n+1)!^4\,(n+2)!}
{(2n+2)!\,(2n+3)!}\delta_{m,n}.
\end{equation}
In particular, $L(1)=1/6$.\quad \quad \qed
\end{Theorem}

\smallskip
Now, by combining Theorems~\ref{cor:cfracHankel}, \ref{thm:momentgf}, 
and \ref{thm:cHahn}, and by using an integral representation of
Bernoulli numbers (see 
\machSeite{Noerlund}%
\cite[p.~75]{Noerlund}),
$$B_\nu=\frac{1}{2\pi \complexi }
\int_{-\infty\complexi }^{\infty\complexi }z^\nu
\left(\frac{\pi}{\sin\pi z}\right)^2 dz
$$
(if $\nu=0$ or $\nu=1$ then the path of integration is indented so
that it avoids the singularity $z=0$, passing it on the negative side)
 we obtain without difficulty the desired
determinant evaluation,
\begin{align} 
\notag
\det_{0\le i,j,\le n-1}(B_{i+j+2})&=
(-1)^{\binom n2}\(\frac {1} {6}\)^n\,\,\prod _{i=1} ^{n-1}\(\frac
{i(i+1)^2(i+2)} {4(2i+1)(2i+3)}\)^{n-i}\\
\label{eq:Bernoulli-det}
&=
(-1)^{\binom n2}\frac {1} {6}\prod _{i=1} ^{n-1}\frac
{i!\,(i+1)!^4\,(i+2)!} {(2i+2)!\,(2i+3)!}.
\end{align}
The general determinant evaluation which results from using continuous
Hahn polynomials with generic nonnegative integers
$a,b,c,d$ is worked out in
\machSeite{FuKrAD}%
\cite[Sec.~5]{FuKrAD}.

Let me mention that, given a Hankel determinant evaluation such as
\eqref{eq:Bernoulli-det}, one has automatically proved a more general
one, by means of the following simple fact (see
for example 
\machSeite{MuirAD}%
\cite[p.~419]{MuirAD}):

\begin{Lemma} \label{thm:Hankel-x}
Let $x$ be an indeterminate. For any nonnegative integer $n$ there
holds
\begin{equation} \label{eq:Hankel-x}
\det_{0\le i,j\le n-1}(A_{i+j})=
\det_{0\le i,j\le n-1}\(\sum _{k=0} ^{i+j}\binom
{i+j}kA_{k}x^{i+j-k}\).
\end{equation}
\quad \quad \qed
\end{Lemma}

The idea of using continued fractions and/or
orthogonal polynomials for the evaluation of Hankel determinants has
been also exploited in
\machSeite{AlCaAA}%
\machSeite{DelsAA}%
\machSeite{MilnAN}%
\machSeite{MilnAO}%
\machSeite{MilnAP}%
\machSeite{MilnAQ}%
\cite{AlCaAA,DelsAA,MilnAN,MilnAO,MilnAP,MilnAQ}. Some of these
results are exhibited in Theorem~\ref{thm:Hankel}. See the
remarks after Theorem~\ref{thm:Hankel} for pointers to further Hankel
determinant evaluations.

\subsection{Miscellaneous} \label{sec:misc}

This section is a collection of various further results on
determinant evaluation of the
general sort, which I personally like, regardless whether
they may be more or less useful.

Let me begin with a result by Strehl and Wilf
\machSeite{StWiAA}%
\cite[Sec.~II]{StWiAA}, a special case of which was already in the
seventies advertised by van der Poorten 
\machSeite{PoorAA}%
\cite[Sec.~4]{PoorAA} as 
`a determinant evaluation that should be better known'. (For a
generalization see 
\machSeite{KedlAA}%
\cite{KedlAA}.)

\begin{Lemma} \label{thm:StWi}
Let $f(x)$ be a formal power series. Then for any positive integer $n$ 
there holds
\begin{equation} \label{eq:StWi}
\det_{1\le i,j\le n}\(\(\frac {d} {dx}\)^{i-1}f(x)^{a_j}\)
=\(\frac {f'(x)} {f(x)}\)^{\binom n2} f(x)^{a_1+\dots+a_n}
\prod _{1\le i<j\le n} ^{}(a_j-a_i).
\end{equation}
\quad \quad \qed
\end{Lemma}
By specializing, this result allows for the quick proof of various,
sometimes surprising, determinant evaluations, see 
Theorems~\ref{thm:StWi-cor} and \ref{thm:StWi-cor2}.

\medskip
An extremely beautiful determinant evaluation is the evaluation of
the determinant of the {\em circulant matrix}.
\begin{Theorem} \label{thm:circulant}
Let $n$ by a fixed positive integer, and let 
$a_0,a_1,\dots,a_{n-1}$ be indeterminates. Then
\begin{equation} \label{eq:circulant}
\det\pmatrix a_0&a_1&a_2&\dots&a_{n-1}\\
a_{n-1}&a_0&a_1&\dots&a_{n-2}\\
a_{n-2}&a_{n-1}&a_0&\dots&a_{n-3}\\
\hdotsfor5\\
a_{1}&a_2&a_3&\dots&a_{0}
\endpmatrix=\prod _{i=0}
^{n-1}(a_0+\om^{i}a_1+\om^{2i}a_2+\dots+\om^{(n-1)i}a_{n-1}),
\end{equation}
where $\om$ is a primitive $n$-th root of unity.\quad \quad \qed
\end{Theorem}
Actually, the {\em circulant determinant} is just a very special case in a whole
family of determinants, called {\em group determinants}. This would
bring us into the vast territory of group representation theory, and
is therefore beyond the scope of this article. It must suffice to
mention that the group determinants were in fact the {\em cause of birth
of group representation theory} (see 
\machSeite{LamTAA}%
\cite{LamTAA} for a
beautiful introduction into these matters).

\medskip
The next theorem does actually not give the evaluation of a
determinant, but of a {\em Pfaffian}. The Pfaffian $\Pf(A)$ of a
skew-symmetric $(2n)\times(2n)$ matrix $A$ is defined by
$$\Pf(A)=\sum _{\pi} ^{}(-1)^{c(\pi)}\prod _{(ij)\in \pi} ^{}A_{ij},$$
where the sum is over all perfect matchings $\pi$ of the complete
graph on $2n$ vertices, where $c(\pi)$ is the {\em crossing
number} of $\pi$, and where the product is over all edges $(ij)$,
$i<j$, in the matching $\pi$ (see e.g\@. 
\machSeite{StemAE}%
\cite[Sec.~2]{StemAE}). 
What links Pfaffians
so closely to determinants is (aside from similarity of definitions)
the fact that the Pfaffian of a skew-symmetric matrix is, up to sign,
the square root of its determinant. That is,
$\det(A)=\Pf(A)^2$ for any skew-symmetric $(2n)\times(2n)$ matrix $A$
(cf\@. 
\machSeite{StemAE}%
\cite[Prop.~2.2]{StemAE}).\footnote{Another point of view, 
beautifully set forth in 
\machSeite{KnutAF}%
\cite{KnutAF}, is that ``Pfaffians are more fundamental than
determinants, in the sense that determinants are merely the bipartite
special case of a general sum over matchings."} 
Pfaffians play an important role, for example,
in the enumeration of plane partitions, due to the results by Laksov,
Thorup and Lascoux 
\machSeite{LaLTAA}%
\cite[Appendix, Lemma~(A.11)]{LaLTAA} and Okada
\machSeite{OkadAA}%
\cite[Theorems~3 and 4]{OkadAA}
on sums of minors of a given matrix (a combinatorial view as enumerating 
nonintersecting lattice paths with varying starting and/or ending
points has been given by Stembridge
\machSeite{StemAE}%
\cite[Theorems~3.1, 3.2, and 4.1]{StemAE}), 
and their generalization in form 
of the powerful {\em minor summation formulas} due to Ishikawa and
Wakayama 
\machSeite{IsWaAA}%
\cite[Theorems~2 and 3]{IsWaAA}.

Exactly in this context, the context of enumeration of plane
partitions, Gordon 
\machSeite{GordAB}%
\cite[implicitly in Sec.~4, 5]{GordAB} 
(see also 
\machSeite{StemAE}%
\cite[proof of Theorem~7.1]{StemAE}) proved two extremely
useful reductions of Pfaffians to determinants.
\begin{Lemma} \label{lem:Gordon}
Let $(g_i)$ be a sequence with the property $g_{-i}=g_i$, and let $N$
be a positive integer. Then
\begin{equation} \label{eq:Gordon-even}
\Pf\limits_{1\le i<j\le 2N}\bigg(\sum _{-(j-i)<\al\le j-i} ^{}g_\al\bigg)=
\det\limits_{1\le i,j\le N}(g_{i-j}+g_{i+j-1}),
\end{equation}
and
\begin{equation} \label{eq:Gordon-odd}
\Pf\limits_{1\le i<j\le 2N+2}\(\left\{\begin{matrix} 
\sum _{-(j-i)<\al\le j-i} ^{}g_\al&j\le 2N+1\\
X&j=2N+2
\end{matrix}\right\}\)=
X\cdot\det\limits_{1\le i,j\le N}(g_{i-j}-g_{i+j}).
\end{equation}
{\em(}In these statements only one half of the entries of the
Pfaffian is given, the other half being uniquely determined by
skew-symmetry{\em)}.
\quad \quad \qed
\end{Lemma}
This result looks somehow technical, but its usefulness was
sufficiently proved by its applications in the enumeration of plane
partitions and tableaux in 
\machSeite{GordAB}%
\cite{GordAB} and
\machSeite{StemAE}%
\cite[Sec.~7]{StemAE}.

\medskip
Another technical, but useful result is due to Goulden and Jackson
\machSeite{GoJaAJ}%
\cite[Theorem~2.1]{GoJaAJ}.
\begin{Lemma} \label{lem:GoJa}
Let $F_m(t)$, $G_m(t)$ and $H_m(t)$ by formal power series, with
$H_m(0)=0$, $m=0,1,\dots,n-1$. Then for any positive integer $n$
there holds
\begin{equation} \label{eq:GoJa}
\det_{0\le i,j,\le n-1}\(\CT\(\frac {F_j(t)}
{H_j(t)^{i}}G_i(H_j(t))\)\)=
\det_{0\le i,j\le n-1}\(\CT\(\frac {F_j(t)}
{H_j(t)^{i}}G_i(0)\)\),
\end{equation}
where $\CT(f(t))$ stands for the constant term of the Laurent series
$f(t)$.\quad \quad \qed
\end{Lemma}
What is the value of this theorem? In some cases, out of a given
determinant evaluation, it immediately
implies a more general one, containing (at least) one more parameter. 
For example, consider the determinant
evaluation \eqref{eq:MRR}. Choose $F_j(t)=t^j(1+t)^{\mu+j}$,
$H_j(t)=t^2/(1+t)$, and $G_i(t)$ such that
$G_i(t^2/(1+t))=(1+t)^k+(1+t)^{-k}$ for a fixed $k$ 
(such a choice does indeed exist; see 
\machSeite{GoJaAJ}%
\cite[proof of Cor.~2.2]{GoJaAJ}) 
in Lemma~\ref{lem:GoJa}. This yields
$$\det_{0\le i,j\le n-1}\(\binom {\mu+k+i+j}{2i-j}+\binom
{\mu-k+i+j}{2i-j}\)=
\det_{0\le i,j\le n-1}\(2\binom {\mu+i+j}{2i-j}\).$$
Thus, out of the validity of \eqref{eq:MRR}, this enables to establish the
validity of \eqref{eq:AB}, and even of \eqref{eq:Chu1}, by choosing
$F_j(t)$ and $H_j(t)$ as above, but $G_i(t)$ such that 
$G_i(t^2/(1+t))=(1+t)^{x_i}+(1+t)^{-x_i}$, $i=0,1,\dots,n-1$.

\section{A list of determinant evaluations}
\label{sec3}

In this section I provide a list of determinant evaluations, some of
which are very frequently met, others maybe not so often. In any case, I
believe that all of them are useful or attractive, or even both.
However, this is not intended to be, and cannot possibly be, an
exhaustive list of known determinant evaluations. The selection
depends totally on my taste. This may
explain that many of these determinants arose in the enumeration of 
plane partitions and rhombus tilings. On the other hand, it is exactly
this field (see 
\machSeite{PropAA}%
\machSeite{RobbAA}%
\machSeite{StanAA}%
\machSeite{StanAI}%
\cite{PropAA,RobbAA,StanAA,StanAI} for more information on
these topics) which is a particular rich source of nontrivial
determinant evaluations. If you do not find ``your" determinant here,
then, at least, the many references given in this section or the general
results and methods from Section~\ref{sec2} may turn out to be helpful.

Throughout this section we use the standard hypergeometric and basic
hypergeometric notations. To wit, for nonnegative integers $k$
the {\em shifted factorial\/} $(a)_k$ is defined (as already before) by 
$$(a)_k:=a(a+1)\cdots(a+k-1),$$ 
so that in particular $(a)_0:=1$.
Similarly, for nonnegative integers $k$ the {\em shifted $q$-factorial\/} 
$(a;q)_k$ is given by 
$$(a;q)_k:=(1-a)(1-aq)\cdots(1-aq^{k-1}),$$ 
so that $(a;q)_0:=1$. 
Sometimes we make use of the notations $[\al]_q:=(1-q^\al)/(1-q)$,
$[n]_q!:=[n]_q[n-1]_q\cdots [1]_q$, $[0]_q!:=1$. 
The {\em $q$-binomial coefficient\/} is defined by
$$\begin{bmatrix} \al\\k\end{bmatrix}_q:=\frac
{[\al]_q[\al-1]_q\cdots[\al-k+1]_q}
{[k]_q!}=\frac {(1-q^\al)(1-q^{\al-1})\cdots(1-q^{\al-k+1})}
{(1-q^k)(1-q^{k-1})\cdots(1-q)}.$$
Clearly we have $\lim_{q\to1}\[\smallmatrix \al\\k\endsmallmatrix\]_q=
\binom \al k$.

Occasionally shifted ($q$-)factorials will appear which contain a
subscript which is a negative integer. By convention, a shifted
factorial $(a)_k$, where $k$ is a negative integer, is interpreted as
$(a)_k:= {1}/ {(a-1)(a-2)\cdots (a+k)}$, whereas 
a shifted
$q$-factorial $(a;q)_k$, where $k$ is a negative integer, is interpreted as
$(a;q)_k:={1}/ (1-q^{a-1})(1-q^{a-2})\cdots \break(1-q^{a+k})$.
(A uniform way to define the shifted factorial, for positive {\em and}
negative $k$, is by
$(a)_k:=\Ga(a+k)/\Ga(a)$, respectively by an appropriate limit in
case that $a$ or $a+k$ is a nonpositive integer, see 
\machSeite{GrKPAA}%
\cite[Sec.~5.5, p.~211f]{GrKPAA}. A uniform way to define the shifted
$q$-factorial is by means of $(a;q)_k:=(a;q)_\infty/(aq^k;q)_\infty$,
see 
\machSeite{GaRaAA}%
\cite[(1.2.30)]{GaRaAA}.)

\medskip
We begin our list with two determinant evaluations which generalize
the Vandermonde determinant evaluation \eqref{eq:Vandermonde} in a
nonstandard way.
The determinants appearing in these evaluations 
can be considered as ``augmentations" of the
Vandermonde determinant by columns which are formed by
differentiating ``Van\-der\-monde-type" columns. 
(Thus, these determinants can also be considered as
certain {\em generalized Wronskians}.)
Occurences of the first determinant can be found e.g\@. in
\machSeite{FlHaAA}%
\cite{FlHaAA}, 
\machSeite{MehtAA}%
\cite[App.~A.16]{MehtAA}, 
\machSeite{MehtAB}%
\cite[(7.1.3)]{MehtAB}, 
\machSeite{ScheAA}%
\cite{ScheAA}, 
\machSeite{VoigAA}%
\cite{VoigAA}. (It is called ``confluent alternant" in 
\machSeite{MehtAA}%
\machSeite{MehtAB}%
\cite{MehtAA,MehtAB}.)
The motivation in
\machSeite{FlHaAA}%
\cite{FlHaAA} to study these determinants came from {\em Hermite
interpolation} and the {\em analysis of linear recursion relations}. 
In \machSeite{MehtAA}%
\cite[App.~A.16]{MehtAA}, special cases of these determinants
are used in the context of {\em random matrices}. Special cases arose
also in the context of {\em transcendental number theory} (see
\machSeite{PoorAA}%
\cite[Sec.~4]{PoorAA}).

\begin{Theorem} \label{thm:FlHa1}
Let $n$ be a nonnegative integer, and
let $A_m(X)$ denote the $n\times m$ matrix
$$\begin{pmatrix} 1&0&0&0&\dots&0\\
X&1&0&0&\dots&0\\
X^2&2X&2&0&\dots&0\\
X^3&3X^2&6X&6&\dots&0\\
\hdotsfor6\\
X^{n-1}&(n-1)X^{n-2}&(n-1)(n-2)X^{n-3}&\hdotsfor2&
(n-1)\cdots(n-m+1)X^{n-m}
\end{pmatrix},$$
i.e., any next column is formed by differentiating the previous
column with respect to $X$. 
Given a composition of $n$, $n=m_1+\dots+m_\ell$, there holds
\begin{equation} \label{eq:FlHa1}
\det_{1\le i,j,\le n}\big(A_{m_1}(X_1)\,A_{m_2}(X_2)\dots
A_{m_\ell}(X_\ell)\big)=
\bigg(\prod _{i=1} ^{\ell}\prod _{j=1} ^{m_i-1}j!\bigg)
\prod _{1\le i<j\le \ell} ^{}(X_j-X_i)^{m_im_j}.
\end{equation}
\quad \quad \qed
\end{Theorem}
The paper 
\machSeite{FlHaAA}%
\cite{FlHaAA} has as well an ``Abel-type" variation of this result.

\begin{Theorem} \label{thm:FlHa2}
Let $n$ be a nonnegative integer, and
let $B_m(X)$ denote the $n\times m$ matrix
$$\begin{pmatrix} 1&0&0&0&\dots&0\\
X&X&X&X&\dots&X\\
X^2&2X^2&4X^2&\hphantom{3}8X^2&\dots&2^{m-1}X^2\\
X^3&3X^3&9X^3&27X^3&\dots&3^{m-1}X^3\\
\hdotsfor6\\
X^{n-1}&(n-1)X^{n-1}&(n-1)^2X^{n-1}&\hdotsfor2&
(n-1)^{m-1}X^{n-1}
\end{pmatrix},$$
i.e., any next column is formed by applying the operator $X(d/dX)$.
Given a composition of $n$, $n=m_1+\dots+m_\ell$, there holds
\begin{equation} \label{eq:FlHa2}
\det_{1\le i,j,\le n}\big(B_{m_1}(X_1)\,B_{m_2}(X_2)\dots
B_{m_\ell}(X_\ell)\big)=
\bigg(\prod _{i=1} ^{\ell}X_i^{\binom {m_i}2}\prod _{j=1}
^{m_i-1}j!\bigg)
\prod _{1\le i<j\le \ell} ^{}(X_j-X_i)^{m_im_j}.
\end{equation}
\quad \quad \qed
\end{Theorem}

As Alain Lascoux taught me, the natural environment for this type of
determinants is {\em divided differences} and (generalized) {\em discrete
Wronskians}. The divided difference $\partial_{x,y}$ is a linear
operator which maps polynomials in $x$ and $y$ to polynomials
symmetric in $x$ and $y$, and is defined by
$$\partial_{x,y}f(x,y)=\frac {f(x,y)-f(y,x)} {x-y}.$$
Divided differences have been introduced by Newton to solve the
interpolation problem in one variable. (See 
\machSeite{LascAF}%
\cite{LascAF} for an
excellent introduction to interpolation, divided differences, and
related matters, such as {\em Schur functions} and
{\em Schubert polynomials}.) In fact, given a polynomial $g(x)$ in
$x$, whose coefficients do not depend on $a_1,a_2,\dots,a_m$, 
{\em Newton's interpolation formula} reads as follows (cf\@. e.g\@.
\machSeite{LascAF}%
\cite[(Ni2)]{LascAF}),
\begin{multline} \label{eq:Newton}
g(x)=g(a_1)+(x-a_1)\partial_{a_1,a_2}g(a_1)
+(x-a_1)(x-a_2)\partial_{a_2,a_3}\partial_{a_1,a_2}g(a_1)\\
+(x-a_1)(x-a_2)(x-a_3)\partial_{a_3,a_4}\partial_{a_2,a_3}
\partial_{a_1,a_2}g(a_1)+\cdots.
\end{multline}

Now suppose that $f_1(x),f_2(x),\dots,f_n(x)$ are polynomials in one
variable $x$, whose coefficients do not depend on $a_1,a_2,\dots,a_n$, 
and consider the determinant 
\begin{equation} \label{eq:det-fi}
\det_{1\le i,j,\le n}(f_i(a_j)).
\end{equation}
Let us for the moment concentrate on the first $m_1$ columns of this
determinant. We may apply \eqref{eq:Newton}, and write
\begin{multline*} 
f_i(a_j)=f_i(a_1)+(a_j-a_1)\partial_{a_1,a_2}f_i(a_1)
+(a_j-a_1)(a_j-a_2)\partial_{a_2,a_3}\partial_{a_1,a_2}f_i(a_1)\\
+\dots +(a_j-a_1)(a_j-a_2)\cdots(a_j-a_{j-1})
\partial_{a_{j-1},a_j}\cdots\partial_{a_2,a_3}
\partial_{a_1,a_2}f_i(a_1),
\end{multline*}
$j=1,2,\dots,m_1$. Following 
\machSeite{LascAF}%
\cite[Proof of Lemma~(Ni5)]{LascAF}, 
we may perform column reductions to the effect that the
determinant \eqref{eq:det-fi}, with column $j$ replaced by
$$(a_j-a_1)(a_j-a_2)\cdots(a_j-a_{j-1})
\partial_{a_{j-1},a_j}\cdots\partial_{a_2,a_3}
\partial_{a_1,a_2}f_i(a_1),$$
$j=1,2,\dots,m_1$, has the same value as the original determinant. 
Clearly, the product $\prod
_{k=1} ^{j-1}(a_j-a_k)$ can be taken out of column $j$,
$j=1,2,\dots,m_1$. Similar reductions can be applied to the next
$m_2$ columns, then to the next $m_3$ columns, etc.

This proves the following fact about generalized discrete Wronskians: 

\begin{Lemma} \label{lem:Wronski-allg}
Let $n$ be a nonnegative integer, and
let $W_m(x_1,x_2,\dots,x_m)$ denote the $n\times m$ matrix
$\pmatrix \partial_{x_{j-1},x_j}\cdots\partial_{x_2,x_3}
\partial_{x_1,x_2}f_i(x_1)\endpmatrix_{1\le i\le n,\ 1\le j\le m}$.
Given a composition of $n$, $n=m_1+\dots+m_\ell$, there holds
\begin{multline} \label{eq:Wronski-allg}
\det_{1\le i,j,\le n}\big(W_{m_1}(a_1,\dots,a_{m_1})\,
W_{m_2}(a_{m_1+1},\dots,a_{m_1+m_2})\dots
W_{m_\ell}(a_{m_1+\dots+m_{\ell-1}+1},\dots,a_n)\big)\\
= {\det_{1\le i,j,\le n}(f_i(a_j))}\bigg/ {\prod _{k=1} ^{\ell}\bigg(\prod
_{m_1+\dots+m_{k-1}+1\le i<j\le m_1+\dots+m_k} ^{}(a_j-a_i)}\bigg).
\end{multline}
\quad \quad \qed
\end{Lemma}

If we now choose $f_i(x):=x^{i-1}$, so that 
${\det_{1\le i,j,\le n}(f_i(a_j))}$ 
is a Vandermonde determinant, then the right-hand
side of \eqref{eq:Wronski-allg} factors completely by
\eqref{eq:Vandermonde}. The final step to
obtain Theorem~\ref{thm:FlHa1} is to let $a_1\to X_1$, $a_2\to X_1$, \dots,
$a_{m_1}\to X_1$, $a_{m_1+1}\to X_2$, \dots,
$a_{m_1+m_2}\to X_2$, etc., in \eqref{eq:Wronski-allg}. This does
indeed yield \eqref{eq:FlHa1}, because
$$\lim_{x_j\to x}\dots\lim_{x_2\to x}\lim_{ x_1\to x}
\partial_{x_{j-1},x_j}\cdots\partial_{x_2,x_3}
\partial_{x_1,x_2}g(x_1)=\frac {1} {(j-1)!}\(\frac {d}
{dx}\)^{j-1}g(x),$$
as is easily verified.

The Abel-type variation in Theorem~\ref{thm:FlHa2} follows from
Theorem~\ref{thm:FlHa1} by multiplying column $j$ in \eqref{eq:FlHa1}
by $X_1^{j-1}$ for $j=1,2,\dots,m_1$, 
by $X_2^{j-m_1-1}$ for $j=m_1+1,m_1+2,\dots,m_2$, etc., and by then using
the relation
$$X\frac {d} {dX}g(X)=\frac {d} {dX}Xg(X)-g(X)$$
many times, so that a typical entry 
$X_k^{j-1}(d/dX_k)^{j-1}X_k^{i-1}$ in row $i$ and column $j$ of the
$k$-th submatrix is expressed as $(X_k(d/dX_k))^{j-1}X_k^{i-1}$
plus a linear combination of terms $(X_k(d/dX_k))^{s}X_k^{i-1}$ with
$s<j-1$. Simple column reductions then yield \eqref{eq:FlHa2}.

\medskip
It is now not very difficult to adapt this analysis 
to derive, for example, $q$-analogues of Theorems~\ref{thm:FlHa1} and
\ref{thm:FlHa2}. The results
below do actually contain $q$-analogues of extensions of
Theorems~\ref{thm:FlHa1} and \ref{thm:FlHa2}.

\begin{Theorem} \label{thm:qFlHa1}
Let $n$ be a nonnegative integer, and
let $A_m(X)$ denote the $n\times m$ matrix
\begin{multline} \notag
\left(\begin{matrix} 1&[C]_q\,X^{-1}&[C]_q\,[C-1]_q\,X^{-2}\\
X&[C+1]_q&[C+1]_q\,[C]_q\,X^{-1}\\
X^2&[C+2]_q\,X&[C+2]_q\,[C+1]_q\\
\hdotsfor3\\
X^{n-1}&[C+n-1]_q\,X^{n-2}&[C+n-1]_q\,[C+n-2]_q\,X^{n-3}
\end{matrix}\right.\\
\left.\begin{matrix} 
\dots&[C]_q\cdots [C-m+2]_q\,X^{1-m}\\
\dots&[C+1]_q\cdots
[C-m+3]_q\,X^{2-m}\\
\dots&[C+2]_q\cdots
[C-m+4]_q\,X^{3-m}\\
\hdotsfor2\\
\dots&[C+n-1]_q\cdots
[C+n-m+1]_q\,X^{n-m}
\end{matrix}\right),
\end{multline} 
i.e., any next column is formed by applying the operator
$X^{-C}D_q\,X^C$, with $D_q$ denoting the usual $q$-derivative, $D_q
f(X):=(f(qX)-f(X))/(q-1)X$.
Given a composition of $n$, $n=m_1+\dots+m_\ell$, there holds
\begin{multline} \label{eq:qFlHa1}
\det_{1\le i,j,\le n}\big(A_{m_1}(X_1)\,A_{m_2}(X_2)\dots
A_{m_\ell}(X_\ell)\big)\\
=
q^{N_1}\bigg(\prod _{i=1} ^{\ell}\prod _{j=1} ^{m_i-1}[j]_q!\bigg)
\prod _{1\le i<j\le \ell} ^{}\prod _{s=0} ^{m_i-1}\prod _{t=0} ^{m_j-1}
(q^{t-s}X_j-X_i),
\end{multline}
where $N_1$ is the quantity
$$\textstyle\sum\limits _{i=1} ^{\ell}\sum\limits _{j=1}
^{m_i}\((C+j+m_1+\dots+m_{i-1}-1)(m_i-j)-\binom {m_i}3\)-
\sum\limits _{1\le i<j\le \ell} ^{}\(m_i\binom {m_j}2-m_j\binom {m_i}2\).$$
\quad \quad \qed
\end{Theorem}

To derive \eqref{eq:qFlHa1} one would choose strings of geometric
sequences for the variables $a_j$ in Lemma~\ref{lem:Wronski-allg},
i.e., $a_1=X_1$, $a_2=qX_1$, $a_3=q^2X_1$, \dots, $a_{m_1+1}=X_2$,
$a_{m_1+2}=qX_2$, etc., and, in addition, use the relation
\begin{equation} \label{eq:div-comm}
y^C\partial_{x,y}f(x,y)=\partial_{x,y}(x^Cf(x,y))-(\partial_{x,y}x^C)f
(x,y)
\end{equation}
repeatedly.

\medskip
A ``$q$-Abel-type" variation of this result reads as follows.

\begin{Theorem} \label{thm:qFlHa2}
Let $n$ be a nonnegative integer, and
let $B_m(X)$ denote the $n\times m$ matrix
$$\begin{pmatrix} 1&[C]_q&[C]_q^2&\dots&[C]_q^{m-1}\\
X&[C+1]_q\,X&[C+1]_q^2\,X&\dots&[C+1]_q^{m-1}\,X\\
X^2&[C+2]_q\,X^2&[C+2]_q^2\,X^2&\dots&[C+2]_q^{m-1}\,X^2\\
\hdotsfor5\\
X^{n-1}&[C+n-1]_q\,X^{n-1}&[C+n-1]_q^2\,X^{n-1}&\dots&[C+n-1]_q^{m-1}\,X^
{n-1}
\end{pmatrix},$$
i.e., any next column is formed by applying the operator
$X^{1-C}D_qX^C$, with $D_q$ denoting the $q$-derivative as in
Theorem~\ref{thm:qFlHa1}.
Given a composition of $n$, $n=m_1+\dots+m_\ell$, there holds
\begin{multline} \label{eq:qFlHa2}
\det_{1\le i,j,\le n}\big(B_{m_1}(X_1)\,B_{m_2}(X_2)\dots
B_{m_\ell}(X_\ell)\big)\\
=
q^{N_2}\bigg(\prod _{i=1} ^{\ell}X_i^{\binom {m_i}2}
\prod _{j=1} ^{m_i-1}[j]_q!\bigg)
\prod _{1\le i<j\le \ell} ^{}\prod _{s=0} ^{m_i-1}\prod _{t=0} ^{m_j-1}
(q^{t-s}X_j-X_i),
\end{multline}
where $N_2$ is the quantity
$$\textstyle\sum\limits _{i=1} ^{\ell}\sum\limits _{j=1}
^{m_i}\((C+j+m_1+\dots+m_{i-1}-1)(m_i-j)\)-
\sum\limits _{1\le i<j\le \ell} ^{}\(m_i\binom {m_j}2-m_j\binom {m_i}2\).$$
\quad \quad \qed
\end{Theorem}

Yet another generalization of the Vandermonde determinant evaluation
is found in 
\machSeite{BuCRAA}%
\cite{BuCRAA}. Multidimensional analogues are contained in
\machSeite{TaVaAB}%
\cite[Theorem~A.7, Eq.~(A.14), Theorem~B.8, Eq.~(B.11)]{TaVaAB}
and 
\machSeite{VarcAA}%
\cite[Part I, p.~547]{VarcAA}.

\medskip
Extensions of Cauchy's double alternant \eqref{eq:Cauchy} 
can also be found in the
literature (see e.g\@. 
\machSeite{MilnAR}%
\machSeite{RobbAF}%
\cite{MilnAR,RobbAF}). I want to mention here particularly 
Borchardt's variation 
\machSeite{BorcAA}%
\cite{BorcAA} in
which the $(i,j)$-entry in Cauchy's double alternant is replaced by
its square,
\begin{equation} \label{eq:Borchardt}
\det_{1\le i,j\le n}\(\frac {1} {(X_i - Y_j)^2}\) =
\frac {\prod_{1\le i<j\le n} 
(X_i - X_j)(Y_i - Y_j)} { \prod_{1\le i,j\le n} (X_i - Y_j)}
\underset{1\le i,j\le n}\Per\(\frac {1} {X_i - Y_j}\) ,
\end{equation}
where $\Per M$ denotes the {\em permanent\/} of the matrix $M$. Thus,
there is no closed form expression such as in \eqref{eq:Cauchy}. This
may not look that useful. However, most remarkably, there is a
($q$-)deformation of this identity which did indeed lead to a
``closed form evaluation," thus solving a famous enumeration problem
in an unexpected way,
the problem of enumerating {\em alternating sign
matrices}.\footnote{\label{foot:ASM}%
An alternating sign matrix is a square matrix with
entries $0,1,-1$, with all row and column sums equal to 1, and such
that, on disregarding the $0$s, 
in each row and column the $1$s and $(-1)$s alternate. Alternating
sign matrix are currently the most fascinating, and most mysterious, 
objects in enumerative combinatorics. The reader is referred to
\machSeite{BrPrAA}%
\machSeite{RobbAA}%
\machSeite{BresAO}%
\machSeite{KupeAD}%
\machSeite{MiRRAB}%
\machSeite{ZeilBD}%
\machSeite{ZeilAO}%
\cite{BresAO,BrPrAA,MiRRAB,RobbAA,KupeAD,ZeilBD,ZeilAO} 
for more detailed material.
Incidentally, the ``birth" of alternating sign matrices came through
--- determinants, see 
\machSeite{RoRuAA}%
\cite{RoRuAA}.}
This $q$-deformation is equivalent to 
Izergin's evaluation 
\machSeite{IzerAA}%
\cite[Eq.~(5)]{IzerAA} (building on results by Korepin 
\machSeite{KoreAA}%
\cite{KoreAA}) of the partition function of the {\em 
six-vertex model\/} under certain boundary conditions (see also
\machSeite{KupeAD}%
\cite[Theorem~8]{KupeAD} and 
\machSeite{KoBIAA}%
\cite[Ch.~VII, (10.1)/(10.2)]{KoBIAA}).
\begin{Theorem} \label{thm:IzKor}
For any nonnegative integer $n$ there holds
\begin{multline} \label{eq:IzKor}
\det_{1\le i,j\le n}\(\frac {1} {(X_i - Y_j)(q X_i - Y_j)}\) 
   = \frac {\prod_{1\le i<j\le n} 
(X_i - X_j)(Y_i - Y_j)} { \prod_{1\le i,j\le n} (X_i - Y_j)(q X_i - Y_j)    }\\
\times
       \sum_{A} (1-q)^{2 N(A)} \prod_{i=1}^{n} X_i^{N_i(A)} Y_i^{N^i(A)}
         \prod_{i,j \text{ \em such that }A_{ij} = 0} (\alpha_{i,j} X_i - Y_j),
\end{multline}
where the sum is over all $n \times n$ alternating sign matrices
$A=(A_{ij})_{1\le i,j\le n}$,
$N(A)$ is the number of $(-1)$s in $A$, $N_i(A)$ {\em(}respectively 
$N^i(A)${\em)} is the
number of $(-1)$s in the $i$-th row {\em(}respectively column{\em)} 
of $A$, and $\alpha_{ij} = q$ if $\sum_{k=1}^{j}
A_{ik} = \sum_{k=1}^{i} A_{kj}$, and $\alpha_{ij} = 1$ otherwise.
\end{Theorem}
Clearly, equation \eqref{eq:Borchardt} results immediately from
\eqref{eq:IzKor} by setting $q=1$. Roughly, Kuperberg's 
solution 
\machSeite{KupeAD}%
\cite{KupeAD} of the enumeration of alternating sign
matrices consisted of suitably specializing the $x_i$'s, $y_i$'s
and $q$
in \eqref{eq:IzKor}, so that each summand on the right-hand side
would reduce to the same quantity, 
and, thus, the sum would basically count $n\times n$
alternating sign matrices, and in evaluating the left-hand side
determinant for that special choice of the $x_i$'s, $y_i$'s and $q$.
The resulting number of $n\times n$ alternating sign matrices is
given in \eqref{eq:ASM} in the Appendix.
(The first, very different, solution is due to Zeilberger 
\machSeite{ZeilBD}%
\cite{ZeilBD}.)
Subsequently, Zeilberger 
\machSeite{ZeilAO}%
\cite{ZeilAO} improved on Kuperberg's approach and succeeded in
proving the refined alternating sign matrix conjecture from 
\machSeite{MiRRAB}%
\cite[Conj.~2]{MiRRAB}.
For a different expansion of the determinant of
Izergin, in terms of {\em Schur functions},
and a variation, see 
\machSeite{LascAG}%
\cite[Theorem~q, Theorem~$\ga$]{LascAG}.

\medskip
Next we turn to typical applications of Lemma~\ref{lem:Krat1}. They
are listed in the following theorem. 

\begin{Theorem} \label{thm:PP}
Let $n$ be a nonnegative integer, and let
$L_1,L_2,\dots,L_n$ and $A,B$ be indeterminates. Then there hold
\begin{equation} \label{eq:PP1}
\det_{1\le i,j\le n}\(\bmatrix
L_i+A+j\\L_i+j\endbmatrix_q\)
=q^{\sum _{i=1} ^{n}(i-1)(L_i+i)}\frac {\prod _{1\le i<j\le n}
^{}[L_i-L_j]_q} {\prod _{i=1} ^{n}[L_i+n]_q!}
\frac {\prod _{i=1} ^{n}[L_i+A+1]_q!} {\prod _{i=1} ^{n}[A+1-i]_q!},
\end{equation}
and
\begin{equation} \label{eq:PP2}
\det_{1\le i,j\le n}\(q^{jL_i}\bmatrix
A\\L_i+j\endbmatrix_q\)
=q^{\sum _{i=1} ^{n}iL_i}\frac {\prod _{1\le i<j\le n}
^{}[L_i-L_j]_q} {\prod _{i=1} ^{n}[L_i+n]_q!}
\frac {\prod _{i=1} ^{n}[A+i-1]_q!} {\prod _{i=1} ^{n}[A-L_i-1]_q!},
\end{equation}
and
\begin{multline} \label{eq:PP3}
\det_{1\le i,j\le n}\bigg(\pmatrix
BL_i+A\\L_i+j\endpmatrix\bigg)\\
=\frac {\prod _{1\le i<j\le n} ^{}(L_i-L_j)} {\prod _{i=1}
^{n}(L_i+n)!} \prod _{i=1} ^{n}\frac {(BL_i+A)!} {((B-1)L_i+A-1)!}
\prod _{i=1} ^{n}(A-Bi+1)_{i-1}\ ,
\end{multline}
and
\begin{equation} \label{eq:Abel}
\det_{1\le i,j\le n}\left(\frac {(A+BL_i)^{j-1}}
{(j-L_i)!}\right)
=\prod _{i=1} ^{n}\frac {(A+Bi)^{i-1}}
{(n-L_i)!}\prod _{1\le i<j\le n} ^{}(L_j-L_i).
\end{equation}
\quad \quad \qed
\end{Theorem}
(For derivations of \eqref{eq:PP1} and \eqref{eq:PP2} using
Lemma~\ref{lem:Krat1} see the proofs of Theorems~6.5 and 6.6 in
\machSeite{KratAM}%
\cite{KratAM}. For a derivation of \eqref{eq:PP3} using
Lemma~\ref{lem:Krat1} see the proof of Theorem~5 in 
\machSeite{KratAK}%
\cite{KratAK}.)

Actually, the evaluations \eqref{eq:PP1} and \eqref{eq:PP2} are
equivalent. This is seen by observing that 
$$\[\matrix 
L_i+A+j\\L_i+j\endmatrix\]_q=(-1)^{L_i+j}q^{\binom {L_i}2+\binom
j2+jL_i+(A+1)(L_i+j)}
\[\matrix 
-A-1\\L_i+j\endmatrix\]_q.$$
Hence, replacement of $A$ by $-A-1$ in \eqref{eq:PP1} leads to
\eqref{eq:PP2} after little manipulation.

The determinant evaluations \eqref{eq:PP1} and \eqref{eq:PP2}, 
and special cases thereof, are rediscovered and reproved 
in the literature over and over. 
(This phenomenon will probably persist.) To the best of my
knowledge, the evaluation \eqref{eq:PP1} appeared in print explicitly 
for the first time in 
\machSeite{CarlAH}%
\cite{CarlAH}, although it was (implicitly)
known earlier to people in group representation theory, as it also results
from the principal specialization (i.e., set $x_i=q^i$,
$i=1,2,\dots,N$) 
of a Schur function of arbitrary
shape, by comparing the Jacobi--Trudi identity with the
bideterminantal form (Weyl character formula) of the Schur function 
(cf\@. 
\machSeite{MacdAC}%
\cite[Ch.~I, (3.4),
Ex.~3 in Sec.~2, Ex.~1 in Sec.~3]{MacdAC}; the determinants arising in the
bideterminantal form are Vandermonde determinants and therefore
easily evaluated).

The main applications of \eqref{eq:PP1}--\eqref{eq:PP3} are in the
enumeration of {\em tableaux, plane partitions and rhombus tilings}.
For example, the {\em hook-content
formula} 
\machSeite{StanAA}%
\cite[Theorem~15.3]{StanAA} for tableaux of a given shape
with bounded entries follows immediately from the theory of
nonintersecting lattice paths (cf\@. 
\machSeite{GeViAB}%
\cite[Cor.~2]{GeViAB}
and 
\machSeite{StemAE}%
\cite[Theorem~1.2]{StemAE}) and the determinant evaluation
\eqref{eq:PP1} (see 
\machSeite{GeViAB}%
\cite[Theorem~14]{GeViAB} and 
\machSeite{KratAM}%
\cite[proof of 
Theorem~6.5]{KratAM}). MacMahon's {\em ``box formula"}
\machSeite{MacMahon}%
\cite[Sec.~429; proof in Sec.~494]{MacMahon} for the generating
function of plane partitions which are contained inside a given box
follows from nonintersecting lattice paths and the determinant
evaluation \eqref{eq:PP2} (see 
\machSeite{GeViAB}%
\cite[Theorem~15]{GeViAB} and 
\machSeite{KratAM}%
\cite[proof of 
Theorem~6.6]{KratAM}).
The $q=1$ special case of the determinant which is relevant here is
the one in \eqref{eq:MacMahon} (which is the one which was evaluated as an
illustration in Section~\ref{sec:general}).
To the best of my knowledge, the evaluation \eqref{eq:PP3} is due to
Proctor 
\machSeite{ProcAA}%
\cite{ProcAA} who used it 
for enumerating {\em plane partitions of staircase shape} (see also
\machSeite{KratAK}%
\cite{KratAK}). 
The determinant evaluation \eqref{eq:Abel} can be used to
give closed form expressions in the enumeration of {\em $\la$-parking
functions} (an extension of the notion of $k$-parking functions such
as in 
\machSeite{StanAZ}%
\cite{StanAZ}), if one starts with determinantal expressions due
to Gessel (private communication). 
Further applications of \eqref{eq:PP1}, in the domain of {\em
multiple} ({\em basic}) {\em hypergeometric series}, 
are found in 
\machSeite{GuKrAB}%
\cite{GuKrAB}. Applications of these determinant
evaluations in statistics are contained in 
\machSeite{HaMoAB}%
\cite{HaMoAB} and
\machSeite{StecAB}%
\cite{StecAB}.

It was pointed out in 
\machSeite{DT}%
\cite{DT} that plane partitions in a given box
are in bijection with {\em rhombus tilings of a ``semiregular"
hexagon}.
Therefore, the determinant \eqref{eq:MacMahon} counts as well rhombus
tilings in a hexagon with side lengths $a,b,n,a,b,n$. In this regard,
generalizations of the evaluation of this determinant, and of a special
case of \eqref{eq:PP3}, appear in
\machSeite{CiucAH}%
\cite{CiucAH} and 
\machSeite{CiEKAA}%
\cite{CiEKAA}. The theme of these papers is to
enumerate {\em rhombus tilings of a hexagon with triangular holes}.


\medskip
The next theorem provides a typical application of
Lemma~\ref{lem:Krat2}. For a derivation of this determinant
evaluation using this lemma see 
\machSeite{KratAO}%
\cite[proofs of Theorems~8 and
9]{KratAO}.
\begin{Theorem} \label{thm:shifted}
Let $n$ be a nonnegative integer, and let
$L_1,L_2,\dots,L_n$ and $A$ be indeterminates. Then there holds
\begin{multline} \label{eq:shifted}
\det_{1\le i,j\le n}\(q^{jL_i}\bmatrix
L_i+A-j\\L_i+j\endbmatrix_q\)\\
=q^{\sum _{i=1} ^{n}iL_i}\prod _{i=1} ^{n}\frac {[L_i+A-n]_q!} {[L_i+n]_q!\,[A-2i]_q!}
\prod _{1\le i<j\le n} ^{}\big([L_i-L_j]_q\,[L_i+L_j+A+1]_q\big) .
\end{multline}
\quad \quad \qed
\end{Theorem} 

This result was used to compute generating
functions for {\em shifted plane partitions of trapezoidal shape} (see
\machSeite{KratAO}%
\cite[Theorems~8 and 9]{KratAO}, 
\machSeite{ProcAF}%
\cite[Prop.~4.1]{ProcAF} and
\machSeite{ProcAB}%
\cite[Theorem~1]{ProcAB}).  

\medskip
Now we turn to typical applications of Lemma~\ref{lem:Krat3}, given
in Theorems~\ref{thm:PP5}--\ref{thm:PP6} below.
All of them can be derived
in just the same way as we evaluated the determinant \eqref{eq:MacMahon} 
in Section~\ref{sec:general} (the only difference being that
Lemma~\ref{lem:Krat3} is invoked instead of Lemma~\ref{lem:Krat1}).

The first application is the evaluation of a determinant 
whose entries are a product of two $q$-binomial coefficients.

\begin{Theorem} \label{thm:PP5}
Let $n$ be a nonnegative integer, and let
$L_1,L_2,\dots,L_n$ and $A,B$ be indeterminates. Then there holds
\begin{multline} \label{eq:PP5}
\det_{1\le i,j\le n}\(\bmatrix
L_i+j\\B\endbmatrix_q\cdot
\bmatrix L_i+A-j\\B\endbmatrix_q\)\\
=q^{\sum _{i=1} ^{n}(i-1)L_i-B\binom n2+2\binom {n+1}3}
\prod _{1\le i<j\le n} ^{}\big([L_i-L_j]_q\,[L_i+L_j+A-B+1]_q\big)
\\
\times
\prod _{i=1} ^{n}\frac {[L_i+1]_q!\,[L_i+A-n]_q!} 
{[L_i-B+n]_q!\,[L_i+A-B-1]_q!}
\frac {[A-2i-1]_q!} {[A-i-n-1]_q!\,[B+i-n]_q!\,[B]_q!}.
\end{multline}
\quad \quad \qed
\end{Theorem}
As is not difficult to verify, this determinant evaluation contains
\eqref{eq:PP1}, \eqref{eq:PP2}, as well as \eqref{eq:shifted} as special,
respectively limiting cases.

This determinant evaluation found applications in basic
hypergeometric functions theory. In 
\machSeite{WilJAA}%
\cite[Sec.~3]{WilJAA}, Wilson used
a special case
to construct {\em biorthogonal rational functions}. On the other hand,
Schlosser applied it in 
\machSeite{SchlAD}%
\cite{SchlAD} to find several new summation
theorems for {\em multidimensional basic hypergeometric series}.

In fact, as Joris Van der Jeugt pointed out to me, there is a
generalization of Theorem~\ref{thm:PP5} of the following form
(which can be also proved by means of Lemma~\ref{lem:Krat3}).
\begin{Theorem} \label{thm:PP5a}
Let $n$ be a nonnegative integer, and let
$X_0,X_1,\dots,X_{n-1}$, $Y_0,Y_1,\dots,\break Y_{n-1}$, $A$ and $B$ be 
indeterminates. Then there holds
\begin{multline} \label{eq:PP5a}
\det_{0\le i,j\le n-1}\(\frac {\bmatrix X_i+Y_j\\j\endbmatrix_q
\bmatrix Y_j+A-X_i\\j\endbmatrix_q}
{\bmatrix X_i+B\\j\endbmatrix_q \bmatrix A+B-X_i\\j\endbmatrix_q}\)\\
=q^{2\binom n3+\sum_{i=0}^{n-1}i(X_i+Y_i-A-2B)}
   \prod_{0\le i<j\le n-1}^{}[X_i-X_j]_q\,[X_i+X_j-A]_q\\
\times
   \prod_{i=0}^{n-1}\frac {(q^{B-Y_i-i+1})_i\,(q^{Y_i+A+B+2-2i})_i} 
     {  (q^{X_i-A-B})_{n-1}\,(q^{X_i+B-n+2})_{n-1}}.
\end{multline}
\quad \quad \qed
\end{Theorem}

As another application of Lemma~\ref{lem:Krat3} we list two
evaluations of determinants (see below) where the entries are, up to
some powers of $q$,
a difference of two $q$-binomial coefficients. 
A proof of the first evaluation which uses 
Lemma~\ref{lem:Krat3} can be found in
\machSeite{KratAP}%
\cite[proof of Theorem~7]{KratAP}, a proof of the second
evaluation using Lemma~\ref{lem:Krat3} can be found in
\machSeite{SchlAA}%
\cite[Ch.~VI, \S3]{SchlAA}. Once more, the second
evaluation was always (implicitly) known to 
people in group representation theory, as it also results
from a principal specialization (set $x_i=q^{i-1/2}$,
$i=1,2,\dots$) 
of a symplectic character of arbitrary
shape, by comparing the symplectic dual Jacobi--Trudi identity with the
bideterminantal form (Weyl character formula) 
of the symplectic character (cf\@. 
\machSeite{FuHaAA}%
\cite[Cor.~24.24 and (24.18)]{FuHaAA}; the determinants arising in the
bideterminantal form are easily evaluated by means of \eqref{eq:Bn}).

\begin{Theorem} \label{thm:PP4}
Let $n$ be a nonnegative integer, and let
$L_1,L_2,\dots,L_n$ and $A$ be indeterminates. Then there hold
\begin{multline} \label{eq:PP4}
\det_{1\le i,j\le n}\bigg(q^{j(L_{j}-L_{i})}
\bigg(\bmatrix A\\j-L_{i}\endbmatrix_q
-q^{j(2L_{i}+A-1)}\bmatrix 
A\\-j-L_{i}+1\endbmatrix_q\bigg)\bigg)\\
=
\prod _{i=1} ^{n}\frac {[A+2i-2]_q!}
{[n-L_{i}]_q!\,[A+n-1+L_{i}]_q!}
\prod _{1\le i<j\le n}
^{}[L_{j}-L_{i}]_q\prod _{1\le i\le j\le n}
^{}[L_{i}+L_{j}+A-1]_q
\end{multline}
and
\begin{multline} \label{eq:PP4a}
\det_{1\le i,j\le n}\bigg(q^{j(L_{j}-L_{i})}
\bigg(\bmatrix A\\j-L_{i}\endbmatrix_q
-q^{j(2L_{i}+A)}\bmatrix 
A\\-j-L_{i}\endbmatrix_q\bigg)\bigg)\\
=
\prod _{i=1} ^{n}\frac {[A+2i-1]_q!}
{[n-L_{i}]_q!\,[A+n+L_{i}]_q!}
\prod _{1\le i<j\le n}
^{}[L_{j}-L_{i}]_q\prod _{1\le i\le j\le n}
^{}[L_{i}+L_{j}+A]_q\ .
\end{multline}
\quad \quad \qed
\end{Theorem}

A special case of \eqref{eq:PP4a} was the second determinant
evaluation which
Andrews needed in 
\machSeite{AndrAK}%
\cite[(1.4)]{AndrAK} in order to prove the 
{\em MacMahon Conjecture} (since then, ex-Conjecture) about the
$q$-enumeration of {\em symmetric plane partitions}. Of course, Andrews'
evaluation proceeded by LU-factorization, while Schlosser 
\machSeite{SchlAA}%
\cite[Ch.~VI, \S3]{SchlAA} simplified Andrews' proof significantly by
making use of Lemma~\ref{lem:Krat3}.
The determinant evaluation \eqref{eq:PP4}
was used in 
\machSeite{KratAP}%
\cite{KratAP} in the proof of
{\em refinements} of the 
{\em MacMahon (ex-)Conjecture} and the {\em Bender--Knuth
(ex-)Conjecture}.
(The latter makes an assertion about the generating function for tableaux
with bounded entries and a bounded number of columns. The first proof
is due to Gordon 
\machSeite{GordAC}%
\cite{GordAC}, the first published proof
\machSeite{AndrAJ}%
\cite{AndrAJ} is due to Andrews.)

Next, in the theorem below, 
we list two very similar determinant evaluations. This time, 
the entries of the determinants are, up to
some powers of $q$,
a sum of two $q$-binomial coefficients. 
A proof of the first evaluation which uses 
Lemma~\ref{lem:Krat3} can be found in
\machSeite{SchlAA}%
\cite[Ch.~VI, \S3]{SchlAA}. A proof of the second evaluation can be
established analogously. Again, the second
evaluation was always (implicitly) known to 
people in group representation theory, as it also results
from a principal specialization (set $x_i=q^{i}$,
$i=1,2,\dots$) 
of an odd orthogonal character of arbitrary
shape, by comparing the orthogonal dual Jacobi--Trudi identity with the
bideterminantal form (Weyl character formula) 
of the orthogonal character (cf\@. 
\machSeite{FuHaAA}%
\cite[Cor.~24.35 and (24.28)]{FuHaAA}; the determinants arising in the
bideterminantal form are easily evaluated by means of \eqref{eq:Cn}).

\begin{Theorem} \label{thm:PP6}
Let $n$ be a nonnegative integer, and let
$L_1,L_2,\dots,L_n$ and $A$ be indeterminates. Then there hold
\begin{multline} \label{eq:PP6}
\det_{1\le i,j\le n}\bigg(q^{(j-1/2)(L_{j}-L_{i})}
\bigg(\bmatrix A\\j-L_{i}\endbmatrix_q
+q^{(j-1/2)(2L_{i}+A-1)}\bmatrix 
A\\-j-L_{i}+1\endbmatrix_q\bigg)\bigg)\\
=
\prod _{i=1} ^{n}\frac {(1+q^{L_i+A/2-1/2})} {(1+q^{i+A/2-1/2})}
\frac {[A+2i-1]_q!} {[n-L_{i}]_q!\,[A+n+L_{i}-1]_q!}\\
\times
\prod _{1\le i<j\le n}
^{}[L_{j}-L_{i}]_q\,[L_{i}+L_{j}+A-1]_q
\end{multline}
and
\begin{multline} \label{eq:PP6a}
\det_{1\le i,j\le n}\bigg(q^{(j-1/2)(L_{j}-L_{i})}
\bigg(\bmatrix A\\j-L_{i}\endbmatrix_q
+q^{(j-1/2)(2L_{i}+A-2)}\bmatrix 
A\\-j-L_{i}+2\endbmatrix_q\bigg)\bigg)\\
=
\frac {\prod _{i=1} ^{n}(1+q^{L_i+A/2-1})} 
{\prod _{i=2} ^{n}(1+q^{i+A/2-1})}
\prod _{i=1} ^{n}\frac {[A+2i-2]_q!} {[n-L_{i}]_q!\,[A+n+L_{i}-2]_q!}\\
\times
\prod _{1\le i<j\le n}
^{}[L_{j}-L_{i}]_q\,[L_{i}+L_{j}+A-2]_q\ .
\end{multline}
\quad \quad \qed
\end{Theorem}
A special case of \eqref{eq:PP6} was the first determinant
evaluation which
Andrews needed in 
\machSeite{AndrAK}%
\cite[(1.3)]{AndrAK} in order to prove the 
{\em MacMahon Conjecture} on {\em symmetric plane partitions}. 
Again, Andrews'
evaluation proceeded by LU-factorization, while Schlosser 
\machSeite{SchlAA}%
\cite[Ch.~VI, \S3]{SchlAA} simplified Andrews' proof significantly by
making use of Lemma~\ref{lem:Krat3}.

\medskip
Now we come to determinants which belong to a different category what
regards difficulty of evaluation, as it is
not possible to introduce more parameters in a substantial way.

\medskip
The first determinant evaluation in this category that we list here
is a determinant evaluation due to Andrews 
\machSeite{AndrAN}%
\machSeite{AndrAO}%
\cite{AndrAN,AndrAO}.
It solved, at the same time, Macdonald's problem of enumerating
{\em cyclically symmetric plane partitions}  
and Andrews' own conjecture about the enumeration of
{\em descending plane partitions}.

\begin{Theorem} \label{thm:Andrews}
Let $\mu$ be an indeterminate. For nonnegative integers $n$ there
holds
\begin{multline} \label{eq:Andrews}
\det_{0\le i,j\le n-1}\left(\delta_{ij}+\binom {2\mu+i+j}{j}\right)\\
=\cases   \displaystyle
2^{\cl{n/2}}\prod _{i=1} ^{n-2}(\mu+\cl{i/2}+1)_{\fl{(i+3)/4}}\\
\displaystyle\times
\frac {\prod _{i=1} ^{n/2}\(\mu+\frac {3n}2-\cl{\frac {3i}2}+\frac 32\)_{\cl{i/2}-1}
                \(\mu+\frac {3n}2-\cl{\frac {3i}2}+\frac 32\)_{\cl{i/2}}}
        {\prod _{i=1} ^{n/2-1}(2i-1)!!\,(2i+1)!!}&\text {if $n$ is
even,}\\
\displaystyle
2^{\cl{n/2}}\prod _{i=1} ^{n-2}(\mu+\cl{i/2}+1)_{\cl{(i+3)/4}}\\
\displaystyle\times
\frac {\prod _{i=1} ^{(n-1)/2}\(\mu+\frac {3n}2-\cl{\frac {3i-1}2}+1\)_{\cl{(i-1)/2}}
           \(\mu+\frac {3n}2-\cl{\frac {3i}2}\)_{\cl{i/2}}}
        {\prod _{i=1} ^{(n-1)/2}(2i-1)!!^2}&\text {if $n$ is odd.}
\endcases
\end{multline}
\quad \quad \qed
\end{Theorem} 

The specializations of this determinant evaluation which are of
relevance for the enumeration of cyclically symmetric plane
partitions and descending plane partitions are the cases $\mu=0$ and
$\mu=1$, respectively. In these cases, Macdonald, respectively Andrews,
actually had conjectures about $q$-enumeration. These were proved by Mills,
Robbins and Rumsey 
\machSeite{MiRRAA}%
\cite{MiRRAA}. Their theorem which solves the
$q$-enumeration of {\em cyclically symmetric plane
partitions} is the following.
\begin{Theorem} \label{thm:CSP}
For nonnegative integers $n$ there holds
\begin{equation} \label{eq:CSP}
\det_{0\le i,j\le n-1}\left(\delta_{ij}+q^{3i+1}\begin{bmatrix}
i+j\\j\end{bmatrix} _{q^3}\right)=
\prod _{i=1} ^{n}\frac {1-q^{3i-1}} {1-q^{3i-2}}
\prod _{1\le i\le j\le n} ^{}\frac {1-q^{3(n+i+j-1)}}
{1-q^{3(2i+j-1)}}.
\end{equation}
\quad \quad \qed
\end{Theorem} 

The theorem by Mills, Robbins and Rumsey in 
\machSeite{MiRRAA}%
\cite{MiRRAA} which
concerns the enumeration of {\em descending plane partitions} is the
subject of the next theorem.
\begin{Theorem} \label{thm:descPP}
For nonnegative integers $n$ there holds
\begin{equation} \label{eq:descPP}
\det_{0\le i,j\le n-1}\left(\delta_{ij}+q^{i+2}\begin{bmatrix}
i+j+2\\j\end{bmatrix}_q\right)=
\prod _{1\le i\le j\le n+1} ^{}\frac {1-q^{n+i+j}}
{1-q^{2i+j-1}}.
\end{equation}
\quad \quad \qed
\end{Theorem} 

It is somehow annoying that so far nobody was able to come
up with a full $q$-analogue of the Andrews determinant
\eqref{eq:Andrews} (i.e., not
just in the cases $\mu=0$ and $\mu=1$). This issue is already addressed
in 
\machSeite{AndrAO}%
\cite[Sec.~3]{AndrAO}. In particular, it is shown there that the
result for a
natural $q$-enumeration of a parametric family of descending plane
partitions does not factor nicely in general, and thus does not lead to a
$q$-analogue of \eqref{eq:Andrews}. Yet, such a $q$-analogue should
exist. Probably the binomial coefficient in \eqref{eq:Andrews} has to be
replaced by something more complicated than just a $q$-binomial times
some power of $q$.

On the other hand, there are surprising variations of
the Andrews determinant \eqref{eq:Andrews}, discovered by Douglas Zare.
These can be interpreted as certain {\em weighted enumerations} of
{\em cyclically
symmetric plane partitions} and of {\em rhombus tilings of a hexagon
with a triangular hole} (see 
\machSeite{CiEKAA}%
\cite{CiEKAA}).

\begin{Theorem} \label{thm:Zare1}
Let $\mu$ be an indeterminate. For nonnegative integers $n$ there
holds
\begin{multline} \label{eq:Zare1}
\det_{0\le i,j\le n-1}\left(-\delta_{ij}+\binom {2\mu+i+j}{j}\right)
\\=
\cases 0,&\text {if $n$ is odd,}\\
(-1)^{n/2}\prod _{i=0} ^{n/2-1}\frac {i!^2\,(\mu+i)!^2\,
(\mu+3i+1)!^2\,(2\mu+3i+1)!^2} {(2i)!\,(2i+1)!\,(\mu+2i)!^2\,
(\mu+2i+1)!^2\,(2\mu+2i)!\,(2\mu+2i+1)!},
&\text {if $n$ is even.}\endcases
\end{multline}

If $\om$ is a primitive 3rd root of unity,
then for nonnegative integers $n$ there holds
\begin{multline} \label{eq:Zare3}
\det_{0\le i,j\le n-1}\left(\om\delta_{ij}+\binom {2\mu+i+j}{j}\right)
=\frac {(1+\om)^n2^{\fl{n/2}}}
        {\prod _{i=1} ^{\fl{n/2}}(2i-1)!!
        \prod _{i=1} ^{\fl{(n-1)/2}}(2i-1)!!}\\
\times
        \prod _{i\ge0} ^{}\(\mu+3i+1\)_{\fl{(n-4i)/2}}
                \(\mu+3i+3\)_{\fl{(n-4i-3)/2}}\\
\cdot
                \(\mu+n-i+\tfrac {1} {2}\)_{\fl{(n-4i-1)/2}}
                \(\mu+n-i-\tfrac {1} {2}\)_{\fl{(n-4i-2)/2}},
\end{multline}
where, in abuse of notation, by $\fl{\al}$ we mean the usual floor
function if $\al\ge0$, however, if $\al<0$ then $\fl{\al}$ must be
read as $0$, so that the product over $i$ in \eqref{eq:Zare3} is indeed a
finite product.

If $\om$ is a primitive 6th root of unity,
then for nonnegative integers $n$ there holds
\begin{multline} \label{eq:Zare2}
\det_{0\le i,j\le n-1}\left(\om\delta_{ij}+\binom {2\mu+i+j}{j}\right)
=\frac {(1+\om)^n\(\frac {2} {3}\)^{\fl{n/2}}}
        {\prod _{i=1} ^{\fl{n/2}}(2i-1)!!
        \prod _{i=1} ^{\fl{(n-1)/2}}(2i-1)!!}\\
\times
        \prod _{i\ge0} ^{}\(\mu+3i+\tfrac {3} {2}\)_{\fl{(n-4i-1)/2}}
                \(\mu+3i+\tfrac {5} {2}\)_{\fl{(n-4i-2)/2}}\\
\cdot
                \(\mu+n-i\)_{\fl{(n-4i)/2}}
                \(\mu+n-i\)_{\fl{(n-4i-3)/2}},
\end{multline}
where again, in abuse of notation, by $\fl{\al}$ we mean the usual floor
function if $\al\ge0$, however, if $\al<0$ then $\fl{\al}$ must be
read as $0$, so that the product over $i$ in \eqref{eq:Zare2} is indeed a
finite product.\quad \quad \qed
\end{Theorem}

There are no really simple proofs of
Theorems~\ref{thm:Andrews}--\ref{thm:Zare1}. Let me just address the
issue of proofs of the evaluation of the Andrews determinant,
Theorem~\ref{thm:Andrews}. The only {\em direct\/}
proof of Theorem~\ref{thm:Andrews} is the original proof of Andrews
\machSeite{AndrAN}%
\cite{AndrAN}, who worked out the LU-factorization of the determinant. 
Today one agrees that the ``easiest" way of evaluating
the determinant \eqref{eq:Andrews} is by first employing a
magnificent factorization theorem 
\machSeite{MiRRAD}%
\cite[Theorem~5]{MiRRAD} due to
Mills, Robbins and Rumsey, and then evaluating each of the two resulting
determinants. For these, for some reason, 
more elementary evaluations exist (see in particular
\machSeite{AnStAA}%
\cite{AnStAA} for such a derivation). What I state below is
a (straightforward) generalization of this factorization theorem from
\machSeite{KratBH}%
\cite[Lemma~2]{KratBH}. 

\begin{Theorem} \label{thm:MRR-Factor}
Let $Z_n(x;\mu,\nu)$ be defined by
$$Z_n(x;\mu,\nu):=\det_{0\le i,j\le
n-1}\(\de_{ij}+\sum _{t=0} ^{n-1}\sum _{k=0} ^{n-1}\binom
{i+\mu}t\binom {k+\nu}{k-t}\binom {j-k+\mu-1}{j-k}x^{k-t}\),$$
let $T_n(x;\mu,\nu)$ be defined by
$$T_n(x;\mu,\nu):=
\det_{0\le i,j\le n-1}\(\sum _{t=i} ^{2j}\binom {i+\mu}{t-i}\binom
{j+\nu}{2j-t}x^{2j-t}\),$$
and let $R_n(x;\mu,\nu)$ be defined by
\begin{multline} \notag
R_n(x;\mu,\nu):=
\det_{0\le i,j\le n-1}\Bigg(\sum _{t=i} ^{2j+1}
\(\binom {i+\mu}{t-i-1}+\binom {i+\mu+1}{t-i}\)\\
\cdot\(\binom {j+\nu}{2j+1-t}+\binom
{j+\nu+1}{2j+1-t}\)x^{2j+1-t}\Bigg).
\end{multline}
Then for all positive integers $n$ there hold
\begin{equation} \label{eq:ZTR1}
Z_{2n}(x;\mu,\nu)=T_n(x;\mu,\nu/2)\,R_n(x;\mu,\nu/2)
\end{equation}
and
\begin{equation} \label{eq:ZTR2}
Z_{2n-1}(x;\mu,\nu)=2\,T_{n}(x;\mu,\nu/2)\,R_{n-1}(x;\mu,\nu/2).
\end{equation}
\quad \quad \qed
\end{Theorem}

The reader should observe that $Z_n(1;\mu,0)$ is identical with the
determinant in \eqref{eq:Andrews}, as the sums in the entries
simplify by means of Chu--Vandermonde summation (see e.g\@. 
\machSeite{GrKPAA}%
\cite[Sec.~5.1, (5.27)]{GrKPAA}). However, also the entries in
the determinants $T_n(1;\mu,0)$ and $R_n(1;\mu,0)$ simplify. The
respective evaluations read as follows (see
\machSeite{MiRRAD}%
\cite[Theorem~7]{MiRRAD} and 
\machSeite{AnBuAA}%
\cite[(5.2)/(5.3)]{AnBuAA}).

\begin{Theorem} \label{thm:MRR}
Let $\mu$ be an indeterminate. For nonnegative integers $n$ there holds
\begin{multline} \label{eq:MRR}
\det_{0\le i,j\le n-1}\(\binom {\mu+i+j}{2i-j} \)\\
=(-1)^{\chi(n\equiv 3\mod 4)}2^{\binom {n-1}2}\prod _{i=1} ^{n-1}
\frac {\(\mu+i+1\)_{\fl{(i+1)/2}}\,
\(-\mu-3n+i+\frac {3} {2}\)_{\fl{i/2}}} {(i)_i},
\end{multline}
where $\chi(\mathcal A)=1$ if $\mathcal A$ is
true and $\chi(\mathcal A)=0$ otherwise, and
\begin{multline} \label{eq:R_n}
\det_{0\le i,j\le n-1}\(\binom {\mu+i+j}{2i-j} +2\binom
{\mu+i+j+2}{2i-j+1}\)\\
=2^n\prod _{i=1} ^{n}\frac {(\mu+i)_{\fl{i/2}}\,(\mu+3n-\fl{\frac
{3i-1} {2}}+\frac {1} {2})_{\fl{(i+1)/2}}} {(2i-1)!!}.
\end{multline}
\quad \quad \qed
\end{Theorem} 
The reader should notice that the determinant in \eqref{eq:MRR} is
the third determinant from the Introduction, \eqref{eq:MRR1}.
Originally, in 
\machSeite{MiRRAD}%
\cite[Theorem~7]{MiRRAD}, Mills, Robbins and Rumsey proved
\eqref{eq:MRR}
by applying their factorization theorem
(Theorem~\ref{thm:MRR-Factor})
the other way round, relying on Andrews' Theorem~\ref{thm:Andrews}. 
However, in the meantime there exist short {\em direct\/} proofs of
\eqref{eq:MRR}, see 
\machSeite{AnStAA}%
\machSeite{KratBI}%
\machSeite{PeWiAA}%
\cite{AnStAA,KratBI,PeWiAA}, either by
LU-factorization, or by ``identification of factors". A proof based
on the determinant evaluation \eqref{eq:Krat} and some combinatorial
considerations is given in 
\machSeite{CiKrAB}%
\cite[Remark~4.4]{CiKrAB}, see the remarks after
Theorem~\ref{thm:Krat}. As shown in
\machSeite{AnBuAA}%
\machSeite{AnStAA}%
\cite{AnBuAA,AnStAA}, the determinant \eqref{eq:R_n} can easily be
transformed into a special case of the determinant in \eqref{eq:Krat}
(whose evaluation is easily proved using condensation, see the
corresponding remarks there). Altogether, this gives an alternative,
and simpler,
proof of Theorem~\ref{thm:Andrews}. 

Mills, Robbins and Rumsey needed the evaluation of \eqref{eq:MRR}
because it allowed them to prove the (at that time) conjectured enumeration
of {\em cyclically symmetric transpose-complementary plane
partitions} (see 
\machSeite{MiRRAD}%
\cite{MiRRAD}). The unspecialized determinants $Z_n(x;\mu,\nu)$
and
$T_n(x;\mu,\nu)$ have combinatorial meanings as well (see
\machSeite{MiRRAA}%
\cite[Sec.~4]{MiRRAA}, respectively 
\machSeite{KratBH}%
\cite[Sec.~3]{KratBH}), as the
{\em weighted enumeration of certain descending plane partitions and
triangularly shaped plane partitions}.

\medskip
It must be mentioned that the determinants $Z_n(x;\mu,\nu)$,
$T_n(x;\mu,\nu)$, $R_n(x;\mu,\nu)$ do also factor nicely for $x=2$.
This was proved by Andrews 
\machSeite{AndrAS}%
\cite{AndrAS} using LU-factorization, 
thus confirming a
conjecture by Mills, Robbins and Rumsey (see 
\machSeite{KratBH}%
\cite{KratBH} for an
alternative proof by ``identification of factors").

\medskip
It was already mentioned in Section~\ref{sec:misc} that there is a
general theorem by Goulden and Jackson 
\machSeite{GoJaAJ}%
\cite[Theorem~2.1]{GoJaAJ}
(see Lemma~\ref{lem:GoJa} and the remarks thereafter)
which, given the evaluation \eqref{eq:MRR}, immediately implies a
generalization containing one more parameter. (This property of the
determinant \eqref{eq:MRR} is called by Goulden and Jackson the {\em
averaging property}.) The resulting determinant evaluation had been
earlier
found by Andrews and Burge 
\machSeite{AnBuAA}%
\cite[Theorem~1]{AnBuAA}. They derived it
by showing that it can be obtained by 
multiplying the matrix underlying the determinant \eqref{eq:MRR} by 
a suitable triangular matrix.

\begin{Theorem} \label{thm:AB}
Let $x$ and $y$ be indeterminates. For nonnegative integers $n$ there holds
\begin{multline} \label{eq:AB}
\det_{0\le i,j\le n-1}\(\binom {x+i+j}{2i-j} + \binom
{y+i+j}{2i-j}\)\\
=(-1)^{\chi(n\equiv 3\mod 4)}2^{\binom n2+1}\prod _{i=1} ^{n-1}
\frac {\(\frac {x+y} {2}+i+1\)_{\fl{(i+1)/2}}\,
\(-\frac {x+y} {2}-3n+i+\frac {3} {2}\)_{\fl{i/2}}} {(i)_i},
\end{multline}
where $\chi(\mathcal A)=1$ if $\mathcal A$ is
true and $\chi(\mathcal A)=0$ otherwise.\quad \quad \qed
\end{Theorem} 
(The evaluation \eqref{eq:AB} does indeed reduce to \eqref{eq:MRR} by
setting $x=y$.)

The above described procedure of Andrews and Burge to multiply a
matrix, whose determinant is known, by an appropriate triangular
matrix, and thus obtain a new determinant evaluation, was
systematically exploited by
Chu 
\machSeite{ChuWAZ}%
\cite{ChuWAZ}. He derives numerous variations of \eqref{eq:AB},
\eqref{eq:R_n}, and
special cases of \eqref{eq:PP3}. We content ourselves with displaying
two typical identities from 
\machSeite{ChuWAZ}%
\cite[(3.1a), (3.5a)]{ChuWAZ}, 
just enough to get an idea of the character of these.

\begin{Theorem} \label{thm:Chu}
Let $x_0,x_1,\dots,x_{n-1}$ and $c$ 
be indeterminates. For nonnegative integers $n$ there hold
\begin{multline} \label{eq:Chu1}
\det_{0\le i,j\le n-1}\(\binom {c+x_i+i+j}{2i-j} + \binom
{c-x_i+i+j}{2i-j}\)\\
=(-1)^{\chi(n\equiv 3\mod 4)}2^{\binom n2+1}\prod _{i=1} ^{n-1}
\frac {\(c+i+1\)_{\fl{(i+1)/2}}\,
\(-c-3n+i+\frac {3} {2}\)_{\fl{i/2}}} {(i)_i}
\end{multline}
and
\begin{multline} \label{eq:Chu2}
\det_{0\le i,j\le n-1}\(
\frac {(2i-j)+(2c+3j+1)(2c+3j-1)} {(c+i+j+\frac {1} {2})(c+i+j-\frac
{1} {2})}
\binom {c+i+j+\frac {1} {2}}{2i-j} \)\\
=(-1)^{\chi(n\equiv 3\mod 4)}2^{\binom {n+1}2+1}\prod _{i=1} ^{n-1}
\frac {\(c+i+\frac {1} {2}\)_{\fl{(i+1)/2}}\,
\(-c-3n+i+2\)_{\fl{i/2}}} {(i)_i},
\end{multline}
where $\chi(\mathcal A)=1$ if $\mathcal A$ is
true and $\chi(\mathcal A)=0$ otherwise.\quad \quad \qed
\end{Theorem} 

The next determinant (to be precise, the special case $y=0$), whose
evaluation is stated in the theorem below,
seems to be closely related to the
Mills--Robbins--Rumsey determinant \eqref{eq:MRR}, although it is in
fact a lot easier to evaluate. Indications that the evaluation
\eqref{eq:MRR} is much deeper than the following evaluation are,
first, that it does not seem to be possible to
introduce a second parameter into the Mills--Robbins--Rumsey
determinant \eqref{eq:MRR} in a similar way, and, second, the much
more irregular form of the right-hand side of \eqref{eq:MRR} (it
contains many floor functions!), as
opposed to the right-hand side of \eqref{eq:Krat}.

\begin{Theorem} \label{thm:Krat}
Let $x,y,n$ be nonnegative integers. Then there
holds
\begin{multline} \label{eq:Krat}
\det_{0\le i,j\le n-1}\(\frac {(x+y+i+j-1)!}
{(x+2i-j)!\,(y+2j-i)!}\)\\
=\prod _{i=0} ^{n-1}\frac {i!\,(x+y+i-1)!\,(2x+y+2i)_i\,(x+2y+2i)_i}
{(x+2i)!\,(y+2i)!}.
\end{multline}
\quad \quad \qed
\end{Theorem}
This determinant evaluation is due to the author, who proved it in
\machSeite{KratBD}%
\cite[(5.3)]{KratBD} as an aside to the (much more difficult) 
determinant evaluations 
which were needed there to settle a conjecture by Robbins and
Zeilberger about a generalization of the enumeration of {\em totally
symmetric self-complementary plane partitions}.
(These are the determinant evaluations of
Theorems~\ref{thm:TSSCPP1} and \ref{thm:TSSCPP2} 
below.) It was proved there
by ``identification of factors". However, Amdeberhan 
\machSeite{AmdeAD}%
\cite{AmdeAD}
observed that it can be easily proved by ``condensation".

Originally there was no application for \eqref{eq:Krat}. However, not much
later, Ciucu 
\machSeite{CiKrAB}%
\cite{CiKrAB} found not just one application. He
observed that if the determinant evaluation \eqref{eq:Krat} is
suitably combined with his beautiful {\em Matchings Factorization
Theorem} 
\machSeite{CiucAB}%
\cite[Theorem~1.2]{CiucAB} (and some combinatorial
considerations), then not only does one obtain
simple proofs for the evaluation of the Andrews determinant
\eqref{eq:Andrews} and the Mills--Robbins--Rumsey determinant
\eqref{eq:MRR}, but also simple proofs for the enumeration of four
different symmetry classes of plane partitions,
{\em cyclically symmetric plane partitions, cyclically symmetric
self-complementary plane partitions} (first proved by Kuperberg
\machSeite{KupeAA}%
\cite{KupeAA}), {\em cyclically symmetric
transpose-complementary plane partitions} (first proved by Mills,
Robbins and Rumsey 
\machSeite{MiRRAD}%
\cite{MiRRAD}), and {\em totally symmetric
self-complementary plane partitions} (first proved by Andrews
\machSeite{AndrAW}%
\cite{AndrAW}).

\medskip
A $q$-analogue of the previous determinant evaluation is contained in
\machSeite{KratBG}%
\cite[Theorem~1]{KratBG}. Again, Amdeberhan 
\machSeite{AmdeAD}%
\cite{AmdeAD} observed
that it can be easily proved by means of ``condensation".

\begin{Theorem}
Let $x,y,n$ be nonnegative integers. Then there
holds
\begin{multline} \label{eq:qKrat}
\det_{0\le i,j\le n-1} \bigg(
{\frac {(q;q)_{x+y+i+j-1}} {(q;q)_{x+2i-j} \,(q;q)_{y+2j-i}}}
{\frac{q^{-2 i j}} {(-q^{x+y+1};q)_{i+j}}}
\bigg)\\
=
\prod _{i=0} ^{n-1}q^{-2 i^2} 
{\frac{(q^2;q^2)_i \,(q;q)_{x+y+i-1} \,(q^{2x+y+2i};q)_i \,(q^{x+2y+2i};q)_i}
{(q;q)_{x+2i} \,(q;q)_{y+2i} \,(-q^{x+y+1};q)_{n-1+i}}}
.
\end{multline}
\quad \quad \qed
\end{Theorem}

The reader should observe that this is not a straightforward $q$-analogue
of \eqref{eq:Krat} as it does contain the terms
$(-q^{x+y+1};q)_{i+j}$ in the determinant, respectively 
$(-q^{x+y+1};q)_{n-1+i}$ in the denominator of the right-hand side
product, which can be cleared only if $q=1$.

A similar determinant evaluation, with some overlap with
\eqref{eq:qKrat}, was found by Andrews and Stanton 
\machSeite{AnStAA}%
\cite[Theorem~8]{AnStAA}
by making use of LU-factorization, in their ``\'etude" on the Andrews and the
Mills--Robbins--Rumsey determinant.

\begin{Theorem} \label{thm:AnSt}
Let $x$ and $E$ be indeterminates and $n$ be a nonnegative integer. 
Then there holds
\begin{multline} \label{eq:AnSt}
\det_{0\le i,j\le n-1}\(\frac {(E/xq^{i};q^2)_{i-j}\,
(q/Exq^i;q^2)_{i-j}\, (1/x^2q^{2+4i};q^2)_{i-j}} {(q;q)_{2i+1-j}\,
(1/Exq^{2i};q)_{i-j}\, (E/xq^{1+2i};q)_{i-j}}\)\\
=\prod _{i=0} ^{n-1}\frac {(x^2q^{2i+1};q)_i\, (xq^{3+i}/E;q^2)_i\,
(Exq^{2+i};q^2)_i} {(x^2q^{2i+2};q^2)_i\, (q;q^2)_{i+1}\,
(Exq^{1+i};q)_i\, (xq^{2+i}/E;q)_i}.
\end{multline}
\quad \quad \qed
\end{Theorem}

The next group of determinants is (with one exception) 
from 
\machSeite{KratBD}%
\cite{KratBD}. These
determinants were needed in the proof of a conjecture by
Robbins and
Zeilberger about a generalization of the enumeration of {\em totally
symmetric self-complementary plane partitions}.

\begin{Theorem} \label{thm:TSSCPP1}
Let $x,y,n$ be nonnegative integers.
Then
\begin{align} \notag
\det_{0\le i,j\le n-1}&\(
\frac
{(x+y+i+j-1)!\,(y-x+3j-3i)}
{(x+2i-j+1)!\,(y+2j-i+1)!}\)\\
\notag
&=\prod _{i=0} ^{n-1}\(\frac {i!\,(x+y+i-1)!\,(2x+y+2i+1)_i\,(x+2y+2i+1)_i}
{(x+2i+1)!\,(y+2i+1)!}\)\\
&\hskip2cm\cdot \sum _{k=0} ^{n}(-1)^k\binom nk (x)_k\,(y)_{n-k}.
\label{eq:TSSCPP1}
\end{align}
\quad \quad \qed
\end{Theorem}

This is Theorem~8 from 
\machSeite{KratBD}%
\cite{KratBD}. A $q$-analogue, provided in
\machSeite{KratBG}%
\cite[Theorem~2]{KratBG}, is the following theorem. 
\begin{Theorem} \label{thm:qTSSCPP1}
Let $x,y,n$ be nonnegative integers. Then there
holds
\begin{multline} \label{eq:qTSSCPP1}
\det_{0\le i,j\le n-1} \bigg(
{\frac{(q;q)_{x+y+i+j-1}\,
(1-q^{y+2j-i}-q^{y+2j-i+1}+q^{x+y+i+j+1})}
 {(q;q)_{x+2i-j+1} \,(q;q)_{y+2j-i+1}}}\\
\hskip5cm \cdot{\frac{q^{-2 i j}} {(-q^{x+y+2};q)_{i+j}}}
\bigg)\\
=
\prod _{i=0} ^{n-1}\(q^{-2 i^2} 
{\frac{(q^2;q^2)_i \,(q;q)_{x+y+i-1} \,(q^{2x+y+2i+1};q)_i \,(q^{x+2y+2i+1};q)_i}
{(q;q)_{x+2i+1} \,(q;q)_{y+2i+1} \,(-q^{x+y+2};q)_{n-1+i}}}\)\\
\times
\sum _{k=0} ^{n} (-1)^k q^{n k} \bmatrix n\\k\endbmatrix_q q^{y k}\, (q^x;q)_{k}
\,(q^y;q)_{n-k}
.
\end{multline}
\quad \quad \qed
\end{Theorem}

Once more, Amdeberhan
observed that, in principle, Theorem~\ref{thm:TSSCPP1} as well as
Theorem~\ref{thm:qTSSCPP1} could be proved by means of
``condensation". However, as of now, nobody provided a proof of the
double sum identities which would establish \eqref{eq:cond} in these cases.

We continue with Theorems~2 and Corollary~3 from 
\machSeite{KratBD}%
\cite{KratBD}.

\begin{Theorem} \label{thm:TSSCPP2}
Let $x,m,n$ be nonnegative integers with $m\le
n$. Under the convention that sums are interpreted by
$$\sum _{r=A+1} ^{B}\text {\rm Expr}(r)=\begin{cases} \hphantom{-}
\sum _{r=A+1} ^{B} \text {\rm Expr}(r)&A<B\\
\hphantom{-}0&A=B\\
-\sum _{r=B+1} ^{A}\text {\rm Expr}(r)&A>B,\end{cases}$$
there holds
\begin{multline} \label{eq:TSSCPP2}
\det_{0\le i,j\le n-1}\bigg(\sum
 _{x+2i-j<r\le x+m+2j-i}
^{}\binom {2x+m+i+j}r\bigg)\\
=\prod _{i=1} ^{n-1}\(\frac {(2x+m+i)!\,(3x+m+2i+2)_i\,(3x+2m+2i+2)_i}
{(x+2i)!\,(x+m+2i)!}\)\hskip2.5cm\\
\times \frac {(2x+m)!} {(x+\fl{m/2})!\,(x+m)!}\cdot\prod _{i=0}
^{\fl{n/2}-1}(2x+2\cl{m/2}+2i+1)\cdot P_1(x;m,n),
\end{multline}
where $P_1(x;m,n)$ is a polynomial in $x$ of degree $\le \fl{m/2}$.

In particular,
for $m=0$ the determinant equals
\begin{equation} \label{eq:tsscpp1}
\cases \displaystyle\prod _{i=0} ^{n-1}\(\frac {i!\,(2x+i)!\,(3x+2i+2)_i^2}
{(x+2i)!^2}\)
\frac{\prodl _{i=0}
^{n/2-1}(2x+2i+1)} {(n-1)!!}&n\text { even}\\
0&n\text { odd,}
\endcases
\end{equation}
for $m=1$, $n\ge1$, it equals
\begin{equation} \label{eq:tsscpp2}
\prod _{i=0} ^{n-1}\(\frac
{i!\,(2x+i+1)!\,(3x+2i+3)_i\,(3x+2i+4)_i}
{(x+2i)!\,(x+2i+1)!}\)
\frac{\prodl _{i=0}
^{\fl{n/2}-1}(2x+2i+3)} {(2\fl{n/2}-1)!!},
\end{equation}
for $m=2$, $n\ge2$, it equals
\begin{multline} \label{eq:tsscpp3}
\prod _{i=0} ^{n-1}\(\frac
{i!\,(2x+i+2)!\,(3x+2i+4)_i\,(3x+2i+6)_i}
{(x+2i)!\,(x+2i+2)!}\)
\frac{\prodl _{i=0}
^{\fl{n/2}-1}(2x+2i+3)} {(2\fl{n/2}-1)!!}\\
\times\frac {1} {(x+1)}\cdot\cases (x+n+1)&n\text { even}\\
(2x+n+2)&n\text { odd,}\endcases
\end{multline}
for $m=3$, $n\ge3$, it equals
\begin{multline} \label{eq:tsscpp4}
\prod _{i=0} ^{n-1}\(\frac
{i!\,(2x+i+3)!\,(3x+2i+5)_i\,(3x+2i+8)_i}
{(x+2i)!\,(x+2i+3)!}\)
\frac{\prodl _{i=0}
^{\fl{n/2}-1}(2x+2i+5)} {(2\fl{n/2}-1)!!}\\
\times\frac {1} {(x+1)}\cdot\cases (x+2n+1)&n\text { even}\\
(3x+2n+5)&n\text { odd,}\endcases
\end{multline}
and for $m=4$, $n\ge4$, it equals
\vskip3pt
\vbox{{}
\begin{multline} \label{eq:tsscpp5}
\prod _{i=0} ^{n-1}\(\frac
{i!\,(2x+i+4)!\,(3x+2i+6)_i\,(3x+2i+10)_i}
{(x+2i)!\,(x+2i+4)!}\)
\frac{\prodl _{i=0}
^{\fl{n/2}-1}(2x+2i+5)} {(2\fl{n/2}-1)!!}\\
\times\frac {1} {(x+1)(x+2)}\cdot\cases (x^2+(4n+3)x+2(n^2+4n+1))
&n\text { even}\\
(2x+n+4)(2x+2n+4)&n\text { odd.}\endcases
\end{multline}
}
\quad \quad \qed
\end{Theorem}

\medskip
One of the most embarrassing failures of ``identification of factors,"
respectively of LU-factorization, is the problem of $q$-enumeration of
{\em totally symmetric plane partitions}, as stated for example in
\machSeite{StanAH}%
\cite[p.~289]{StanAH} or 
\machSeite{StanAI}%
\cite[p.~106]{StanAI}. 
It is now known for quite a while that also
this problem can be reduced to the evaluation of a certain
determinant, by means of Okada's result 
\machSeite{OkadAA}%
\cite[Theorem~4]{OkadAA}
about the sum of all minors of a given matrix, that was already
mentioned in Section~\ref{sec:misc}. In fact, in 
\machSeite{OkadAA}%
\cite[Theorem~5]{OkadAA},
Okada succeeded to
transform the resulting determinant into a reasonably simple one, so
that the problem of $q$-enumerating totally symmetric plane
partitions reduces to resolving the following conjecture.

\begin{Conjecture} \label{thm:qTSP}
For any nonnegative integer $n$ there holds
\begin{equation} \label{eq:Okada}
\det_{1\le i,j\le n}\(T_n^{(1)}+T_n^{(2)}\)=\prod _{1\le i\le j\le
k\le n} ^{}\(\frac {1-q^{i+j+k-1}} {1-q^{i+j+k-2}}\)^2,
\end{equation}
where 
$$T_n^{(1)}=\(q^{i+j-1}\(\begin{bmatrix} i+j-2\\i-1\end{bmatrix}_q
+q\begin{bmatrix} i+j-1\\i\end{bmatrix}_q\)\)_{1\le i,j\le n}$$
and
$$T_n^{(2)}=\begin{pmatrix} 1+q\\
-1&1+q^2&&&0\\
&-1&1+q^3\\
&&-1&1+q^4\\
&0&&\ddots&\ddots\\
&&&&-1&1+q^n
\end{pmatrix}.$$
\end{Conjecture}
While the problem of (plain) enumeration of
totally symmetric plane partitions was solved a few years ago 
by Stembridge 
\machSeite{StemAG}%
\cite{StemAG} (by some ingenious transformations of the
determinant which results directly from Okada's result on the sum of
all minors of a matrix), the problem of $q$-enumeration is still wide
open. ``Identification of factors" cannot even get started because
so far nobody came up with a way of introducing a
parameter in \eqref{eq:Okada} or any equivalent determinant 
(as it turns out, the parameter
$q$ cannot serve as a parameter in the sense of
Section~\ref{sec:ident}), and, apparently, guessing the
LU-factorization is too difficult.

\medskip
Let us proceed by giving a few more determinants which arise in the
enumeration of rhombus tilings.

Our next determinant evaluation is the evaluation of a determinant
which, on disregarding the second binomial coefficient, 
would be just a special
case of \eqref{eq:PP3}, and which, on the other hand, resembles very
much the $q=1$ case of \eqref{eq:PP4}. (It is the determinant that
was shown as \eqref{eq:pentagon1} in the Introduction.) However, neither
Lemma~\ref{lem:Krat1} nor Lemma~\ref{lem:Krat3} suffice to give a
proof. The proof in 
\machSeite{CiKrAC}%
\cite{CiKrAC} by means of ``identification of
factors" is unexpectedly difficult. 

\begin{Theorem} \label{thm:pentagon}
Let $n$ be a positive integer, and let $x$ and
$y$ be nonnegative integers. Then the
following determinant evaluation holds:
\begin{multline} \label{eq:pentagon}
\det_{1\le i,j\le n}\(\binom {x+y+j}{x-i+2j}-\binom
{x+y+j}{x+i+2j}\)\\
=\prod _{j=1} ^{n}\frac {(j-1)!\,(x+y+2j)!\,(x-y+2j+1)_j\,
(x+2y+3j+1)_{n-j}} {(x+n+2j)!\,(y+n-j)!}.
\end{multline}
\quad \quad \qed
\end{Theorem}
This determinant evaluation is used in 
\machSeite{CiKrAC}%
\cite{CiKrAC} to enumerate
rhombus tilings of a certain pentagonal subregion of a hexagon.

\medskip
To see an example of different nature, I present a determinant
evaluation from 
\machSeite{FuKrAC}%
\cite[Lemma~2.2]{FuKrAC}, which can be considered as a
determinant of a mixture of two matrices, 
out of one we take all rows except the
$l$-th, while out of the other we take just the $l$-th row. The
determinants of both of
these matrices could be straightforwardly evaluated by means of
Lemma~\ref{lem:Krat1}. (They are in fact equivalent to special cases of
\eqref{eq:PP3}.) However, to evaluate this mixture is much more
challenging. In particular, the mixture does not anymore factor
completely into ``nice" factors, i.e., the result is of the form
\eqref{eq:UGLY}, so that for a proof one has to resort to the
extension of ``identification of factors" described there.
\begin{Theorem} \label{thm:FuKr2}
Let $n,m,l$ be positive integers such that $1\le l\le n$. Then there
holds
\begin{multline}
\label{eq:FuKr2}
	\det_{1\leq i,j\leq n}\left(
		\begin{cases}
\binom {n+m-i}{m+i-j}\frac {(m+\frac{n-j+1}{2})} {(n+j-2i+1)}
&\text{ if }i\neq l\\
\vrule height17pt width0pt
\binom {n+m-i}{m+i-j}&\text{ if }i=l
		\end{cases}
	\right)\\
= \prod_{i=1}^{n}\frac{(n+m-i)!}{(m+i-1)!\,(2n-2i+1)!}
\prod _{i=1}^{\fl{n/2}}
	\left(
		\po{m+i}{n-2i+1}\,\po{m+i+\frac{1}{2}}{n-2i}
	\right)\kern1.5cm
\\
\times
2^{\frac{(n-1)(n-2)}{2}}\frac{
	\po{m}{n+1}\prod_{j=1}^{n}(2j-1)!
}{
	n!\prod_{i=1}^{\fl{n/2}}\po{2i}{2n-4i+1}
}
\sum_{e=0}^{l-1}
(-1)^{e}\binom{n}{e}\frac{
	(n-2e)\,\po{\frac{1}{2}}{e}
}{
	(m+e)\,(m+n-e)\,\po{\frac{1}{2}-n}{e}
}.
\end{multline}
\quad \quad \qed
\end{Theorem}
In 
\machSeite{FuKrAC}%
\cite{FuKrAC}, this result was applied to 
enumerate all rhombus tilings of a symmetric hexagon that contain a
fixed rhombus. In Section~4 of 
\machSeite{FuKrAC}%
\cite{FuKrAC} there can be found
several conjectures about the enumeration of rhombus tilings with
more than just one fixed rhombus, all of which amount to evaluating
other mixtures of the above-mentioned two determinants.

\medskip
As last binomial determinants, 
I cannot resist to show the, so far, weirdest determinant
evaluations that I am aware of. 
They arose in an attempt
\machSeite{BoHPAA}%
\cite{BoHPAA} by Bombieri, Hunt and van der Poorten 
to improve on Thue's method of approximating an
algebraic number. In their paper, they conjectured the following
determinant evaluation, which, thanks to van der Poorten
\machSeite{PoorAB}%
\cite{PoorAB}, has recently become a theorem (see the subsequent
paragraphs and, in particular, Theorem~\ref{thm:Poorten} and the
remark following it).
\begin{Theorem} \label{conj:Bombieri}
Let $N$ and $l$ be positive integers. Let $M$ be the matrix with rows
labelled by pairs $(i_1,i_2)$ with $0\le i_1\le 2l(N-i_2)-1$
{\em(}the intuition is that the points $(i_1,i_2)$ are 
the lattice points in a right-angled
triangle{\em)}, with columns labelled by 
pairs $(j_1,j_2)$ with $0\le j_2\le N$ and $2l(N-j_2)\le j_1\le l(3N-2j_2)-1$
{\em(}the intuition is that the points $(j_1,j_2)$ are 
the lattice points in a lozenge{\em)}, and entry in row $(i_1,i_2)$
and column $(j_1,j_2)$ equal to
$$\binom {j_1}{i_1}\binom {j_2}{i_2}.$$
Then the determinant of $M$ is given by
$$\pm\(\frac {\prod _{k=0} ^{l-1}k!\prod _{k=2l} ^{3l-1}k!} {\prod
_{k=l} ^{2l-1}k!}\)^{\binom {N+2}3}.$$
\end{Theorem}
This determinant evaluation is just one in a whole series of
conjectured determinant evaluations and greatest common divisors of
minors of a certain matrix, many of them reported in 
\machSeite{BoHPAA}%
\cite{BoHPAA}. 
These conjectures being settled, the
authors of 
\machSeite{BoHPAA}%
\cite{BoHPAA} expect important implications in the
approximation of algebraic numbers.

The case $N=1$ of Theorem~\ref{conj:Bombieri} 
is a special case of \eqref{eq:PP1}, and, thus, on a
shallow level. On the other hand, the next case, $N=2$, is 
already on a considerably deeper level. It was first proved in
\machSeite{KrZeAA}%
\cite{KrZeAA}, by establishing, in fact, a much more general result, given in
the next theorem. It reduces
to the $N=2$ case of Theorem~\ref{conj:Bombieri} for $x=0$, $b=4l$,
and $c=2l$. (In fact, the $x=0$ case of Theorem~\ref{thm:Bombieri}
had already been conjectured in
\machSeite{BoHPAA}%
\cite{BoHPAA}.)

\begin{Theorem} \label{thm:Bombieri}
Let $b,c$ be nonnegative integers, $c\le b$, and
let $\De(x;b,c)$ be the determinant of the $(b+c)\times(b+c)$ matrix
\vskip10pt
\begin{equation} \label{det:Bombieri}
\(
\PfadDicke{.3pt}
\SPfad(0,1),111111111111\endSPfad
\SPfad(0,-1),111111111111\endSPfad
\SPfad(4,-3),222222\endSPfad
\SPfad(8,-3),222222\endSPfad
\Label\o{\raise25pt\hbox{$\displaystyle\binom {x+j}i$}}(2,1)
\Label\o{\raise25pt\hbox{$\displaystyle\binom {x+j}i$}}(6,1)
\Label\o{\raise25pt\hbox{$\displaystyle\binom {2x+j}i$}}(10,1)
\Label\o{\raise25pt\hbox{$\displaystyle\binom {x+j}i$}}(6,-1)
\Label\o{\raise25pt\hbox{$\displaystyle\binom {2x+j}i$}}(10,-1)
\Label\o{\raise25pt\hbox{$\displaystyle 2\binom {x+j}{i-b}$}}(2,-3)
\Label\o{\raise25pt\hbox{$\displaystyle\binom {x+j}{i-b}$}}(6,-3)
\Label\o{\raise15pt\hbox{\LARGE$0$}}(2,-1)
\Label\o{\raise15pt\hbox{\LARGE$0$}}(10,-3)
\Label\o{\raise20pt\hbox{$0\le
i<c\hphantom{b+{}}$\hskip10pt}}(-3,1)
\Label\o{\raise20pt\hbox{$c\le
i<b\hphantom{c+{}}$\hskip10pt}}(-3,-1)
\Label\o{\raise20pt\hbox{$b\le i<b+c$\hskip10pt}}(-3,-3)
\hskip7.2cm
\)
\hbox{\hskip-7.7cm}
\Label\o{\raise15pt\hbox{$0\le j<c$}}(2,3)
\Label\o{\raise15pt\hbox{$c\le j<b$}}(6,3)
\Label\o{\raise15pt\hbox{$b\le j<b+c$}}(10,3)
\hskip7.7cm.\hskip-2cm
\end{equation}
Then
\begin{enumerate}
\item[{\rm (i)}] $\De(x;b,c)=0$ if $b$ is even and $c$ is odd;
\item[{\rm (ii)}] if any of these conditions does not hold, then
\begin{multline} \label{eq:Bombieri}
\De(x;b,c)=(-1)^{c}2^{c}\prod _{i=1} ^{b-c}\frac {\(i+\frac {1}
{2}-\cl{\frac {b} {2}}\)_c} {(i)_c}\\
\times\prod _{i=1} ^{c}\frac {\(x+\cl{\frac {c+i}
{2}}\)_{b-c+\cl{i/2}-\cl{(c+i)/2}} \, \(x+\cl{\frac {b-c+i}
{2}}\)_{\cl{(b+i)/2}-\cl{(b-c+i)/2}}} {\(\frac {1} {2}-\cl{\frac {b} {2}}
+\cl{\frac {c+i}
{2}}\)_{b-c+\cl{i/2}-\cl{(c+i)/2}} \, \(\frac {1} {2}-\cl{\frac {b} {2}}
+\cl{\frac {b-c+i} {2}}\)_{\cl{(b+i)/2}-\cl{(b-c+i)/2}}}.
\end{multline}
\quad \quad \qed
\end{enumerate}
\end{Theorem}
The proof of this result in
\machSeite{KrZeAA}%
\cite{KrZeAA} could be called ``heavy". It proceeded by
``identification of factors". Thus, it was only 
the introduction of the parameter
$x$ in the determinant in \eqref{det:Bombieri} that allowed the
attack on this special case of the conjecture of Bombieri, Hunt and
van der Poorten. However, the authors of
\machSeite{KrZeAA}%
\cite{KrZeAA} (this includes myself) failed to find a way to
introduce a parameter into the determinant in Theorem~\ref{conj:Bombieri}
for generic $N$
(that is, in a way such the determinant would still factor nicely). This
was accomplished by van der Poorten
\machSeite{PoorAB}%
\cite{PoorAB}. He first changed the entries in 
the determinant slightly, without changing
the value of the determinant, and then was able to introduce a
parameter. 
I state his result, 
\machSeite{PoorAB}%
\cite[Sec.~5, Main Theorem]{PoorAB}, 
in the theorem below. For the proof of his result
he used ``identification of factors" as well, thereby considerably
simplifying and illuminating arguments from
\machSeite{KrZeAA}%
\cite{KrZeAA}.

\begin{Theorem} \label{thm:Poorten}
Let $N$ and $l$ be positive integers. Let $M$ be the matrix with rows
labelled by pairs $(i_1,i_2)$ with $0\le i_1\le 2l(N-i_2)-1$, 
with columns labelled by 
pairs $(j_1,j_2)$ with $0\le j_2\le N$ and $0\le j_1\le lN-1$,
and entry in row $(i_1,i_2)$
and column $(j_1,j_2)$ equal to
\begin{equation} \label{eq:Poorten}
(-1)^{i_1-j_1}\binom {-x(N-j_2)}{i_1-j_1}\binom {j_2}{i_2}.
\end{equation}
Then the determinant of $M$ is given by
$$\pm\(\prod _{i=1} ^{l}\binom {x+i-1}{2i-1}\bigg/
\binom {l+i-1}{2i-1}\)^{\binom {N+2}3}.$$
\end{Theorem}
Although not completely obvious, the special case $x=-2l$ establishes
Theorem~\ref{conj:Bombieri}, see 
\machSeite{PoorAB}%
\cite{PoorAB}. Van der Poorten proves as well an evaluation that
overlaps with the $x=0$ case of Theorem~\ref{thm:Bombieri}, see
\machSeite{PoorAB}%
\cite[Sec.~6, Example Application]{PoorAB}.

\medskip
Let us now turn to a few remarkable Hankel determinant
evaluations.
\begin{Theorem} \label{thm:Hankel}
Let $n$ be a positive integer. Then there hold
\begin{equation} \label{eq:Hankel1}
\det_{0\le i,j\le n-1}(E_{2i+2j})=\prod _{i=0} ^{n-1}(2i)!^2,
\end{equation}
where $E_{2k}$ is the $(2k)$-th {\em(}signless{\em)} 
{\em Euler number}, defined through
the generating function\break $1/\cos z=\sum _{k=0}
^{\infty}E_{2k}z^{2k}/(2k)!$,
and
\begin{equation} \label{eq:Hankel1a}
\det_{0\le i,j\le n-1}(E_{2i+2j+2})=\prod _{i=0} ^{n-1}(2i+1)!^2.
\end{equation}

Furthermore, define the {\em Bell polynomials} $\mathcal B_m(x)$ by 
$\mathcal B_{m}(x)=\sum _{k=1} ^{m}S(m,k)x^k$, where $S(m,k)$ is a
{\em Stirling number of the second kind} {\em(}the number of
partitions of an $m$-element set into $k$ blocks; cf\@.
\machSeite{StanAP}%
\cite[p.~33]{StanAP}{\em)}. Then
\begin{equation} \label{eq:Hankel3}
\det_{0\le i,j\le n-1}(\mathcal B_{i+j}(x))=x^{n(n-1)/2}\prod _{i=0} ^{n-1}i!.
\end{equation}

Next, there holds
\begin{equation} \label{eq:Hankel4}
\det_{0\le i,j\le n-1}(H_{i+j}(x))=(-1)^{n(n-1)/2}\prod _{i=0}
^{n-1}i!,
\end{equation}
where $H_m(x)=\sum _{k\ge0} ^{}\frac {m!} {k!\,(m-2k)!}\(-\frac {1}
{2}\)^kx^{m-2k}$ is the $m$-th {\em Hermite polynomial}.

Finally, the following Hankel determinant evaluations featuring {\em
Bernoulli numbers} hold,
\begin{equation} \label{eq:Hankel5}
\det_{0\le i,j,\le n-1}(B_{i+j})=
(-1)^{\binom {n}2}\,\,\prod _{i=1} ^{n-1}
\frac {i!^6} {(2i)!\,(2i+1)!},
\end{equation}
and
\begin{equation} \label{eq:Hankel6}
\det_{0\le i,j,\le n-1}(B_{i+j+1})=
(-1)^{\binom {n+1}2}\frac {1} {2}\,\,\prod _{i=1} ^{n-1}
\frac {i!^3\,(i+1)!^3} {(2i+1)!\,(2i+2)!},
\end{equation}
and
\begin{equation} \label{eq:Hankel7}
\det_{0\le i,j,\le n-1}(B_{i+j+2})=
(-1)^{\binom n2}\frac {1} {6}\,\,\prod _{i=1} ^{n-1}
\frac {i!\,(i+1)!^4\,(i+2)!} {(2i+2)!\,(2i+3)!},
\end{equation}
and
\begin{equation} \label{eq:Hankel8}
\det_{0\le i,j,\le n-1}(B_{2i+2j+2})=
\prod _{i=0} ^{n-1}
\frac {(2i)!\,(2i+1)!^4\,(2i+2)!} {(4i+2)!\,(4i+3)!},
\end{equation}
and
\begin{equation} \label{eq:Hankel9}
\det_{0\le i,j,\le n-1}(B_{2i+2j+4})=
(-1)^n\prod _{i=1} ^{n}
\frac {(2i-1)!\,(2i)!^4\,(2i+1)!} {(4i)!\,(4i+1)!}.
\end{equation}
\quad \quad \qed
\end{Theorem}
All these evaluations can be deduced from continued fractions and
orthogonal polynomials, in the way that was described in
Section~\ref{sec:Hankel}. 
To prove \eqref{eq:Hankel1} and \eqref{eq:Hankel1a} 
one would resort to suitable special
cases of {\em Meixner--Pollaczek polynomials} (see
\machSeite{KoSwAA}%
\cite[Sec.~1.7]{KoSwAA}), and use an integral representation for Euler
numbers, given for example in 
\machSeite{Noerlund}%
\cite[p.~75]{Noerlund},
$$E_{2\nu}=(-1)^{\nu+1}\sqrt{-1}
\int_{-\infty\complexi }^{\infty\complexi }
\frac{(2z)^{2\nu}}{\cos\pi z} dz.
$$
Slightly different proofs of \eqref{eq:Hankel1} can be found in
\machSeite{AlCaAA}%
\cite{AlCaAA} and 
\machSeite{MehtAB}%
\cite[App.~A.5]{MehtAB}, together with more Hankel determinant evaluations
(among which are also \eqref{eq:Hankel5} and \eqref{eq:Hankel7},
respectively).
The evaluation \eqref{eq:Hankel3} can be derived by
considering {\em Charlier polynomials} (see 
\machSeite{DelsAA}%
\cite{DelsAA} for such a
derivation in a special case).
The evaluation \eqref{eq:Hankel4} follows from the fact that Hermite
polynomials are moments of slightly shifted Hermite polynomials, as
explained in 
\machSeite{IsStAB}%
\cite{IsStAB}. In fact, the papers 
\machSeite{IsStAB}%
\cite{IsStAB} and
\machSeite{IsStAC}%
\cite{IsStAC} contain more examples of orthogonal polynomials which
are moments, thus in particular 
implying Hankel determinant evaluations whose entries are
{\em Laguerre polynomials}, {\em Meixner polynomials}, and
{\em Al-Salam--Chihara polynomials}.
Hankel determinants where the entries are (classical)
orthogonal polynomials are also considered in 
\machSeite{KaSzAA}%
\cite{KaSzAA}, where
they are related to {\em Wronskians} of orthogonal polynomials. In
particular, there result Hankel determinant evaluations with entries
being {\em Legendre}, {\em ultraspherical}, and {\em Laguerre
polynomials}
\machSeite{KaSzAA}%
\cite[(12.3), (14.3), (16.5), \S~28]{KaSzAA}, respectively. The
reader is also referred to 
\machSeite{LeclAB}%
\cite{LeclAB}, where illuminating proofs
of these identities
between Hankel determinants and Wronskians are given, by using the
fact that Hankel determinants can be seen as certain {\em Schur
functions} of rectangular shape, and by applying a `master identity'
of Turnbull 
\machSeite{TurnAA}%
\cite[p.~48]{TurnAA} on minors of a matrix.
(The evaluations \eqref{eq:Hankel1}, \eqref{eq:Hankel4} and
\eqref{eq:Hankel5} can be found in
\machSeite{LeclAB}%
\cite{LeclAB} as well, as corollaries to more 
general results.)
Alternative proofs of \eqref{eq:Hankel1}, \eqref{eq:Hankel3} and
\eqref{eq:Hankel4} can be
found in 
\machSeite{RadoAB}%
\cite{RadoAB}, see also 
\machSeite{RadoAD}%
\cite{RadoAD} and
\machSeite{RadoAA}%
\cite{RadoAA}.

Clearly, to prove \eqref{eq:Hankel5}--\eqref{eq:Hankel7} one 
would proceed in the same way as in Section~\ref{sec:Hankel}. (Identity
\eqref{eq:Hankel7} is in fact the evaluation \eqref{eq:Bernoulli-det}
that we derived in Section~\ref{sec:Hankel}.) The evaluations
\eqref{eq:Hankel8} and \eqref{eq:Hankel9} are equivalent to 
\eqref{eq:Hankel7}, because the matrix underlying 
the determinant in \eqref{eq:Hankel7} has
a checkerboard pattern (recall that
Bernoulli numbers with odd indices are
zero, except for $B_1$), and therefore decomposes into the product of a
determinant of the form \eqref{eq:Hankel8} and a determinant of the
form \eqref{eq:Hankel9}.
Very interestingly,
variations of \eqref{eq:Hankel5}--\eqref{eq:Hankel9} arise as
normalization constants in {\em statistical mechanics models}, see
e.g\@.
\machSeite{AuPeAA}%
\cite[(4.36)]{AuPeAA},
\machSeite{CuThAA}%
\cite[(4.19)]{CuThAA}, and
\machSeite{MehtAB}%
\cite[App.~A.5]{MehtAB}.
A common generalization of \eqref{eq:Hankel5}--\eqref{eq:Hankel7} can
be found in 
\machSeite{FuKrAD}%
\cite[Sec.~5]{FuKrAD}. Strangely enough, it was needed
there in the enumeration of rhombus tilings.

In view of Section~\ref{sec:Hankel}, any continued fraction expansion
of the form \eqref{eq:momentgf} gives rise to a Hankel determinant
evaluation. Thus, many more Hankel determinant evaluations follow
e.g\@.
from work by Rogers 
\machSeite{RogeAA}%
\cite{RogeAA}, 
Stieltjes 
\machSeite{StieAA}%
\machSeite{StieAB}%
\cite{StieAA,StieAB}, Flajolet 
\machSeite{FlajAA}%
\cite{FlajAA}, 
Han, Randrianarivony and Zeng
\machSeite{HaZeAA}%
\machSeite{HaRZAA}%
\machSeite{RandAA}%
\machSeite{RandAB}%
\machSeite{RandAC}%
\machSeite{RaZeAA}%
\machSeite{RaZeAB}%
\machSeite{ZengAA}%
\cite{HaZeAA,HaRZAA,RandAA,RandAB,RandAC,RaZeAA,RaZeAB,ZengAA},
Ismail, Masson and Valent 
\machSeite{IsMaAA}%
\machSeite{IsVaAA}%
\cite{IsMaAA,IsVaAA} or
Milne 
\machSeite{MilnAN}%
\machSeite{MilnAO}%
\machSeite{MilnAP}%
\machSeite{MilnAQ}%
\cite{MilnAN,MilnAO,MilnAP,MilnAQ}, 
in particular, evaluations of Hankel determinant
featuring {\em Euler numbers} with odd indices (these are given
through the generating function $\tan z=\sum _{k=0}
^{\infty}E_{2k+1}z^{2k+1}/(2k+1)!$), {\em Genocchi numbers}, 
$q$- and other extensions of {\em Catalan, Euler} and {\em Genocchi
numbers}, and
coefficients in the power series expansion of
{\em Jacobi elliptic functions}. Evaluations of the latter type
played an important role in Milne's recent beautiful results
\machSeite{MilnAN}%
\machSeite{MilnAO}%
\cite{MilnAN,MilnAO} on the number of representations of integers as sums of
$m$-th powers (see also 
\machSeite{MehtAB}%
\cite[App.~A.5]{MehtAB}).

For further evaluations of Hankel determinants, which apparently do
not follow from known results about continued fractions or orthogonal
polynomials, see 
\machSeite{HeGeAA}%
\cite[Prop.~14]{HeGeAA} and 
\machSeite{FuKrAD}%
\cite[Sec.~4]{FuKrAD}.

\medskip
Next we state two charming applications of Lemma~\ref{thm:StWi} (see
\machSeite{WilfAE}%
\cite{WilfAE}).
\begin{Theorem} \label{thm:StWi-cor}
Let $x$ be a nonnegative integer. For any nonnegative integer $n$ there hold
\begin{equation} \label{eq:StWi-cor1}
\det_{0\le i,j\le n}\(\frac {(xi)!} {(xi+j)!}S(xi+j,xi)\)=
\(\frac {x} {2}\)^{\binom {n+1}2}
\end{equation}
where $S(m,k)$ is a {\em Stirling number of the second kind}
{\em(}the number of
partitions of an $m$-element set into $k$ blocks; 
cf\@. 
\machSeite{StanAP}%
\cite[p.~33]{StanAP}{\em)},
and 
\begin{equation} \label{eq:StWi-cor2}
\det_{0\le i,j\le n}\(\frac {(xi)!} {(xi+j)!}s(xi+j,xi)\)=
\(-\frac {x} {2}\)^{\binom {n+1}2},
\end{equation}
where $s(m,k)$ is a {\em Stirling number of the first kind}
{\em(}up to sign, the number of
permutations of $m$ elements with exactly $k$ cycles;
cf\@. 
\machSeite{StanAP}%
\cite[p.~18]{StanAP}{\em)}.\quad \quad \qed
\end{Theorem}

\begin{Theorem} \label{thm:StWi-cor2}
Let $A_{ij}$ denote the number of representations of $j$ as a sum of
$i$ squares of nonnegative integers. Then $\det_{0\le i,j\le
n}(A_{ij})=1$ for any nonnegative integer $n$. 
The same is true if ``squares" is replaced by
``cubes," etc.\quad \quad \qed
\end{Theorem}

\medskip
After having seen so many determinants where rows and columns are
indexed by integers, it is time for a change. There are quite a few
interesting determinants 
whose rows and columns are indexed by (other) combinatorial
objects. (Actually, we already encountered one in
Conjecture~\ref{conj:Bombieri}.)

We start by a determinant where rows and columns are indexed by
permutations. Its beautiful evaluation was obtained at roughly the same time by
Varchenko 
\machSeite{VarcAC}%
\cite{VarcAC} and Zagier 
\machSeite{ZagiAA}%
\cite{ZagiAA}.
\begin{Theorem} \label{thm:Zagier}
For any positive integer $n$ there holds
\begin{equation} \label{eq:Zagier}
\det_{\si,\pi\in \frak S_n}\(q^{\inv(\si\pi^{-1})}\)=
\prod _{i=2} ^{n}(1-q^{i(i-1)})^{\binom ni(i-2)!\,(n-i+1)!},
\end{equation}
where $\frak S_n$ denotes the symmetric group on $n$ elements.\quad \quad \qed
\end{Theorem}
This determinant evaluation appears in
\machSeite{ZagiAA}%
\cite{ZagiAA}
in the study of
certain models in {\em infinite statistics}.
However, as Varchenko et al\@.  
\machSeite{BrVaAA}%
\machSeite{ScVaAA}%
\machSeite{VarcAC}%
\cite{BrVaAA,ScVaAA,VarcAC} show, 
this determinant evaluation is in fact just a special instance in a whole
series of determinant evaluations.
The latter papers give 
evaluations of determinants corresponding to certain bilinear forms
associated to {\em hyperplane arrangements} and {\em matroids}. Some
of these bilinear forms are relevant to the study of {\em hypergeometric
functions} and the {\em representation theory of quantum groups} (see
also 
\machSeite{VarcAD}%
\cite{VarcAD}). In
particular, these results contain analogues of \eqref{eq:Zagier} for all
{\em finite Coxeter groups} as special cases. 
For other developments related to
Theorem~\ref{thm:Zagier} (and different proofs) see 
\machSeite{HaStAA}%
\machSeite{DeHaAA}%
\machSeite{DeHaAB}%
\machSeite{DuKKAA}%
\cite{DeHaAA,DeHaAB,DuKKAA,HaStAA}, tying the subject also to
the {\em representation theory of the symmetric group}, to {\em
noncommutative symmetric functions}, and to {\em free Lie
algebras}, and 
\machSeite{MeSvAA}%
\cite{MeSvAA}. For more remarkable 
determinant evaluations
related to hyperplane arrangements see 
\machSeite{DoTeAA}%
\machSeite{VarcAA}%
\machSeite{VarcAB}%
\cite{DoTeAA,VarcAA,VarcAB}. For more determinant evaluations related
to hypergeometric functions and quantum groups and algebras, see
\machSeite{TaVaAA}%
\machSeite{TaVaAB}%
\cite{TaVaAA,TaVaAB}, where determinants arising in the context of 
{\em Knizhnik-Zamolodchikov equations} are computed.

The results in 
\machSeite{BrVaAA}%
\machSeite{ScVaAA}%
\cite{BrVaAA,ScVaAA} may be considered as a
generalization of the {\em Shapovalov determinant} evaluation
\machSeite{ShapAA}%
\cite{ShapAA}, associated to the {\em Shapovalov form}
in {\em Lie theory}. The latter has since been extended to {\em Kac--Moody
algebras} (although not yet in full generality), see 
\machSeite{CoCaAA}%
\cite{CoCaAA}.

There is a result similar to Theorem~\ref{thm:Zagier} 
for another prominent permutation statistics, 
MacMahon's {\em major index}. (The major
index $\maj(\pi)$ is defined as the sum of all positions of descents
in the permutation $\pi$, see e.g\@. 
\machSeite{FoatAB}%
\cite{FoatAB}.) 
\begin{Theorem} \label{thm:Zagiermaj}
For any positive integer $n$ there holds
\begin{equation} \label{eq:Zagiermaj}
\det_{\si,\pi\in \frak S_n}\(q^{\maj(\si\pi^{-1})}\)=
\prod _{i=2} ^{n}(1-q^{i})^{n!\,(i-1)/i}.
\end{equation}
\end{Theorem}
As Jean--Yves Thibon explained to me, this determinant evaluation follows
from results about the {\em descent algebra} of the symmetric group given in
\machSeite{KrLTAA}%
\cite{KrLTAA}, presented there in an equivalent form, in terms of
{\em noncommutative symmetric functions}. For the details of
Thibon's argument see Appendix~C.
Also the bivariate determinant $\det_{\si,\pi\in
\frak S_n}\(x^{\des(\si\pi^{-1})}q^{\maj(\si\pi^{-1})}\)$ seems to possess
an interesting factorization.

\medskip
The next set of determinant evaluations shows determinants where the
rows and columns are indexed by {\em set partitions}. In what follows,
the set of all partitions of $\{1,2,\dots,n\}$ is denoted by
$\Pi_n$. The number of blocks of a partition $\pi$ is denoted by
$\bk(\pi)$. A partition $\pi$ is called {\em noncrossing}, if there do
not exist $i<j<k<l$ such that both $i$ and $k$ belong to one block, $B_1$
say, while both $j$ and $l$ belong to another block which is
different from $B_1$.
The set of all noncrossing partitions of $\{1,2,\dots,n\}$ is denoted by
$\NC_n$. (For more information about noncrossing partitions see
\machSeite{SimiAD}%
\cite{SimiAD}.)

Further, poset-theoretic, notations which are needed in the following
theorem are:  Given a poset $P$, the {\em join} of two elements $x$
and $y$ in $P$
is denoted by $x\lor_P y$, while the {\em meet} of $x$ and $y$ is
denoted by $x\land_P y$.
The {\em characteristic polynomial} of a poset
$P$ is written as $\chi_P(q)$
(that is, if the maximum element of $P$ has rank $h$ and $\mu$ is
the {\em M\"obius function} of $P$, then $\chi_P(q):=\sum _{p\in P}
^{}\mu(\hat 0,p)q^{h-\rank(p)}$, where $\hat 0$ stands for the minimal
element of $P$).
The symbol $\tilde\chi_P(q)$ denotes the
{\em reciprocal polynomial} $q^h\chi_P(1/q)$ of $\chi_P(q)$.
Finally, $P^*$ is the {\em order-dual} of $P$.
\begin{Theorem} \label{thm:NC}
Let $n$ be a positive integer. Then
\begin{equation} \label{eq:NC1}
\det_{\pi,\ga\in\Pi_n}\(q^{\bk(\pi\land_{\Pi_n}\ga)}\)=
\prod _{i=1} ^{n}\Big(q\,\tilde\chi_{\Pi^*_i}(q)\Big)^{\binom
niB(n-i)},
\end{equation}
where $B(k)$ denotes the $k$-th {\em Bell number} {\em(}the total
number of partitions of a $k$-element set; cf\@. 
\machSeite{StanAP}%
\cite[p.~33]{StanAP}{\em)}. Furthermore,
\begin{equation} \label{eq:NC2}
\det_{\pi,\ga\in\Pi_n}\(q^{\bk(\pi\lor_{\Pi_n}\ga)}\)=
\prod _{i=1} ^{n}\Big(q\,\chi_{\Pi_i}(q)\Big)^{S(n,i)},
\end{equation}
where $S(m,k)$ is a {\em Stirling number of the second kind}
{\em(}the number of partitions of an $m$-element set into $k$ blocks; 
cf\@. 
\machSeite{StanAP}%
\cite[p.~33]{StanAP}{\em)}. Next,
\begin{equation} \label{eq:NC3}
\det_{\pi,\ga\in\NC_n}\(q^{\bk(\pi\land_{\NC_n}\ga)}\)=
q^{\binom {2n-1}n}\prod _{i=1} ^{n}\Big(\tilde\chi_{\NC_i}(q)\Big)^{\binom
{2n-1-i}{n-1}},
\end{equation}
and
\begin{equation} \label{eq:NC4}
\det_{\pi,\ga\in\NC_n}\(q^{\bk(\pi\lor_{\NC_n}\ga)}\)=
q^{\frac {1} {n+1}\binom {2n}n}
\prod _{i=1} ^{n}\Big(\chi_{\NC_i}(q)\Big)^{\binom
{2n-1-i}{n-1}},
\end{equation}
Finally,
\begin{equation} \label{eq:Tutte}
\det_{\pi,\ga\in\NC_n}\(q^{\bk(\pi\lor_{\Pi_n}\ga)}\)=
q^{\binom {2n-1}n}\prod _{i=1} ^{n-1}\(\frac {U_{i+1}(\sqrt q/2)}
{qU_{i-1}(\sqrt q/2)}\)^{\frac {i+1} {n}\binom {2n}{n-1-i}},
\end{equation}
where $U_{m}(x):=\sum _{j\ge0} ^{}(-1)^j\binom
{m-j}j(2x)^{m-2j}$ is the $m$-th
{\em Chebyshev polynomials of the second kind}.\quad \quad \qed
\end{Theorem}

The evaluations \eqref{eq:NC1}--\eqref{eq:NC4} are due to Jackson
\machSeite{JackAD}%
\cite{JackAD}. The last determinant evaluation, \eqref{eq:Tutte}, is
the hardest among those. It was proved independently by Dahab
\machSeite{DahaAA}%
\cite{DahaAA} and Tutte 
\machSeite{TuttAC}%
\cite{TuttAC}. All these determinants are
related to the so-called {\em Birkhoff--Lewis equation} from
chromatic graph theory (see 
\machSeite{DahaAA}%
\machSeite{TuttAB}%
\cite{DahaAA,TuttAB} for more information).

A determinant of somewhat similar type appeared in work by Lickorish
\machSeite{LickAA}%
\cite{LickAA} on {\em $3$-manifold invariants}. Let $\match(2n)$ denote the
set of all noncrossing perfect matchings of $2n$ elements.
Equivalently, $\match(2n)$ can be considered as 
the set of all noncrossing partitions of $2n$ elements
with all blocks containing exactly 2 elements. Lickorish considered a
bilinear form on the linear space spanned by $\match(2n)$. The
corresponding determinant was evaluated by Ko and Smolinsky
\machSeite{KoSmAA}%
\cite{KoSmAA} using an elegant recursive approach, 
and independently by Di Francesco 
\machSeite{DiFrAA}%
\cite{DiFrAA}, whose calculations
are done within the framework of the {\em Temperley--Lieb algebra} (see also
\machSeite{DiGGAA}%
\cite{DiGGAA}).
\begin{Theorem} \label{thm:meander}
For $\al,\be\in\match(2n)$, let $c(\al,\be)$ denote the number of
connected components when the matchings $\al$ and $\be$ are
superimposed. Equivalently, $c(\al,\be)=\bk(\al\lor_{\Pi_{2n}}\be)$. For
any positive integer $n$ there holds
\begin{equation} \label{eq:meander}
\det_{\al,\be\in\match(2n)}\(q^{c(\al,\be)}\)=
\prod _{i=1} ^{n}U_i(q/2)^{a_{2n,2i}},
\end{equation}
where $U_m(q)$ is the {\em Chebyshev polynomial of the second kind}
as given in Theorem~\ref{thm:NC}, and
where $a_{2n,2i}=c_{2n,2i}-c_{2n,2i+2}$ with $c_{n,h}=\binom
{n}{(n-h)/2}-\binom {n}{(n-h)/2-1}$.\quad \quad \qed
\end{Theorem}
Di Francesco 
\machSeite{DiFrAA}%
\cite[Theorem~2]{DiFrAA} does also provide a generalization to
partial matchings, and in 
\machSeite{DiFrAB}%
\cite{DiFrAB} a generalization in an $SU(n)$
setting, the previously mentioned results being situated in the $SU(2)$
setting. While the derivations in 
\machSeite{DiFrAA}%
\cite{DiFrAA} are mostly
combinatorial, the derivations in 
\machSeite{DiFrAB}%
\cite{DiFrAB} are based on
computations in quotients of type $A$ Hecke algebras.

\medskip
There is also an interesting determinant evaluation,
which comes to my mind, where rows and
columns of the determinant are indexed by {\em integer partitions}. 
It is a result due to Reinhart 
\machSeite{ReiGAA}%
\cite{ReiGAA}.
Interestingly, it arose in the analysis of {\em algebraic differential
equations}.

\medskip
In concluding, let me attract your attention to other determinant
evaluations which I like, but which would take too much space to
state and introduce properly. 

For example, there is a determinant evaluation, 
conjectured by Good, and proved by 
Goulden and Jackson 
\machSeite{GoJaAC}%
\cite{GoJaAC}, which arose
in the calculation of cumulants of a statistic analogous to
{\em Pearson's chi-squared for a multinomial sample}. Their method of
derivation is very combinatorial, in particular making use of
generalized {\em ballot sequences}. 

Determinants arising from certain {\em raising operators of $\mathfrak
s\mathfrak
l(2)$-representations} are presented in 
\machSeite{ProcAH}%
\cite{ProcAH}. As special
cases, there result beautiful determinant evaluations where rows and
columns are indexed by integer partitions and the entries are numbers
of {\em standard Young tableaux} of skew shapes.

In 
\machSeite{KraHAA}%
\cite[p.~4]{KraHAA} (see also 
\machSeite{WilRAA}%
\cite{WilRAA}),
an interesting mixture of linear
algebra and combinatorial matrix theory yields, as a by-product, the
evaluation of the determinant of certain 
{\em incidence mappings}. There, rows and columns of
the relevant matrix are indexed by all subsets of an $n$-element set 
of a fixed size. 

As a by-product of
the analysis of an interesting matrix in {\em quantum information
theory}
\machSeite{KrSlAA}%
\cite[Theorem~6]{KrSlAA}, the evaluation of a determinant of a matrix
whose rows and columns are indexed by all subsets of an $n$-element
set is obtained.

Determinant evaluations of {\em $q$-hypergeometric functions} are used in
\machSeite{TerTAA}%
\cite{TerTAA} to compute {\em $q$-Selberg integrals}.

And last, but not least, let me once more mention the remarkable
determinant evaluation, arising in connection with {\em holonomic
$q$-difference equations}, due to Aomoto and Kato 
\machSeite{AoKaAA}%
\cite[Theorem~3]{AoKaAA}, who
thus proved a conjecture by Mimachi 
\machSeite{MimaAD}%
\cite{MimaAD}.

\section*{Appendix A: A word about guessing}

\global\def\theequation{\mbox{A.\arabic{equation}}}
\setcounter{equation}{0}

The problem of
guessing a formula for the generic element $a_n$ of a sequence
$(a_n)_{n\ge0}$ out of the first few elements 
was present at many places, in particular this is crucial for a
successful application of the ``identification of factors" method
(see Section~\ref{sec:ident}) or
of LU-factorization (see Section~\ref{sec:LU}). Therefore some
elaboration on guessing is in order.

First of all, as I already claimed, {\em guessing can be largely
automatized}. This is due to the following tools\footnote{In
addition, one has to mention Frank Garvan's {\tt qseries}
\machSeite{GarvAG}%
\cite{GarvAG}, which is
designed for guessing and computing within the territory of $q$-series,
$q$-products, eta and theta functions, and the like. Procedures like
{\tt prodmake} or {\tt qfactor}, however, might also be helpful for the
evaluation of ``$q$-determinants". The package is available from {\tt
http://www.math.ufl.edu/\~{}frank/qmaple.html}.}:

\begin{enumerate}
\item {\tt Superseeker}, the electronic version of the
``Encyclopedia of Integer Sequences" 
\machSeite{SlPlAA}%
\machSeite{SloaAA}%
\cite{SlPlAA,SloaAA} by
Neil Sloane and Simon Plouffe (see Footnote~\ref{foot:Integer} in the Introduction),
\item {\tt gfun} by Bruno Salvy
and Paul Zimmermann and {\tt Mgfun} by Frederic Chyzak 
(see Footnote~\ref{foot:gfun} in the Introduction),
\item {\tt Rate} by the author (see Footnote~\ref{foot:Rate} in the Introduction).
\end{enumerate}

If you send the first few elements of your sequence to {\tt
Superseeker} then, if it overlaps with a sequence that is stored
there, you will receive information about your sequence such as where
your sequence already appeared in the literature, a formula,
generating function, or a recurrence for your sequence. 

The {\sl Maple} package {\tt
gfun} provides tools for finding a generating function and/or a recurrence 
for your sequence. (In fact, {\tt Superseeker} does also
automatically invoke features from {\tt gfun}.)
{\tt Mgfun} does the same in a multidimensional setting. 

Within the ``hypergeometric paradigm," the most
useful is the {\sl Mathematica} program {\tt Rate} (``Rate!" is
German for ``Guess!"), respectively its {\sl Maple} equivalent {\tt
GUESS}. 
Roughly speaking, it allows to automatically
guess ``closed forms".\footnote{\label{foot:NICE}Commonly, by ``closed form" (``NICE"
in Zeilberger's ``terminology") one means an expression which is
built by forming products and quotients of factorials. A strong
indication that you encounter a sequence $(a_n)_{n\ge0}$ 
for which a ``closed form"
exists is that the prime factors in the prime factorization of $a_n$ do not
grow rapidly as $n$ becomes larger. (In fact, they should grow
linearly.)} 
The program is based on the observation that
any ``closed form" sequence $(a_n)_{n\ge0}$ that appears within the 
``hypergeometric paradigm" is either given by a rational expression,
like $a_n=n/(n+1)$, or the sequence of successive quotients
$(a_{n+1}/a_n)_{n\ge0}$ is given by a rational expression, like
in the case of central binomial coefficients
$a_n=\binom {2n}n$, 
or the sequence of successive quotients of successive quotients
$((a_{n+2}/a_{n+1})/(a_{n+1}/a_n))_{n\ge0}$ 
is given by a rational expression, like in the case of the famous sequence of
numbers of {\em alternating sign matrices} (cf\@. the paragraphs
following \eqref{eq:Borchardt}, and
\machSeite{BresAO}%
\machSeite{BrPrAA}%
\machSeite{MiRRAB}%
\machSeite{RobbAA}%
\machSeite{KupeAD}%
\machSeite{ZeilBD}%
\machSeite{ZeilAO}%
\cite{BresAO,BrPrAA,MiRRAB,RobbAA,KupeAD,ZeilBD,ZeilAO} for information on
alternating sign matrices),
\begin{equation} \label{eq:ASM}
a_n=\prod _{i=0} ^{n-1}\frac {(3i+1)!}
{(n+i)!},
\end{equation}
etc.
Given enough special values,
a rational expression is easily found by rational interpolation. 

This is implemented in {\tt Rate}. Given the first $m$ terms of a
sequence, it takes the first $m-1$ terms and applies rational
interpolation to these, then it applies rational interpolation to
the successive quotients of these $m-1$ terms, etc. For each of the 
obtained results it is tested if it does also give the $m$-th term correctly. 
If it does, then the corresponding result is added to the output, 
if it does not, then nothing is added to the output. 

Here is a short demonstration of the {\sl Mathematica} program {\tt Rate}.
The output shows guesses for the $i0$-th element of the sequence.
\MATH

In[1]:= \MATHkleiner \MATHkleiner rate.m
\goodbreakpoint%
In[2]:= Rate[1,2,3]
\goodbreakpoint%
Out[2]= %
\MATHlbrace i0%
\MATHrbrace 
\goodbreakpoint%
In[3]:= Rate[2/3,3/4,4/5,5/6]
\goodbreakpoint%
         1 + i0
\vskip-5pt%
Out[3]= %
\MATHlbrace ------%
\MATHrbrace 
\vskip-5pt%
\leavevmode         2 + i0

\endMATH
\noindent
Now we try the central binomial coefficients:
\MATH

In[4]:= Rate[1,2,6,20,70]
\goodbreakpoint%
         -1 + i0
\vskip-5pt%
\leavevmode          %
\MATHhStrich %
\MATHtStueck %
\MATHhStrich %
\MATHtStueck %
\MATHhStrich  2 (-1 + 2 i1)
\vskip-5pt%
Out[4]= %
\MATHlbrace   %
\MATHvStrich  %
\MATHvStrich   -------------%
\MATHrbrace 
\vskip-5pt%
\leavevmode           %
\MATHvStrich  %
\MATHvStrich        i1
\leavevmode          i1=1

\endMATH
\noindent
It needs the first 8 values to guess the formula \eqref{eq:ASM} for
the numbers of alternating sign matrices:
\MATH

In[5]:= Rate[1,2,7,42,429,7436,218348,10850216]
\goodbreakpoint%
         -1 + i0   -1 + i1
\vskip-5pt%
\leavevmode          %
\MATHhStrich %
\MATHtStueck %
\MATHhStrich %
\MATHtStueck %
\MATHhStrich      %
\MATHhStrich %
\MATHtStueck %
\MATHhStrich %
\MATHtStueck %
\MATHhStrich  3 (2 + 3 i2) (4 + 3 i2)
\vskip-5pt%
Out[5]= %
\MATHlbrace   %
\MATHvStrich  %
\MATHvStrich   2 (  %
\MATHvStrich  %
\MATHvStrich   -----------------------)%
\MATHrbrace 
\vskip-5pt%
\leavevmode           %
\MATHvStrich  %
\MATHvStrich        %
\MATHvStrich  %
\MATHvStrich   4 (1 + 2 i2) (3 + 2 i2)
\leavevmode          i1=1      i2=1

\endMATH

\medskip
However, what if we encounter a sequence
where all these nice automatic tools fail? Here are a few hints. 
First of all, it is not uncommon to encounter a sequence
$(a_n)_{n\ge0}$ which has actually a split definition. For example,
it may be the case that the subsequence $(a_{2n})_{n\ge0}$ of
even-numbered terms follows a ``nice" formula, and that the
subsequence $(a_{2n+1})_{n\ge0}$ of
odd-numbered terms follows as well a ``nice," but different, formula.
Then {\tt Rate} will fail on any number of first terms of
$(a_n)_{n\ge0}$, while it will give you something for sufficiently
many first terms of $(a_{2n})_{n\ge0}$, and 
it will give you something else for sufficiently
many first terms of $(a_{2n+1})_{n\ge0}$.

Most of the subsequent hints 
apply to a situation where you encounter a sequence
$p_0(x),p_1(x),p_2(x),\dots$
of polynomials $p_n(x)$ in $x$ for which you want to find (i.e., guess) a formula.
This is indeed the situation in which you are generally during the
guessing for ``identification of factors," and also usually when you
perform a guessing where a parameter is involved.

To make things concrete, let us suppose that the first 10 elements of
your sequence of polynomials are

{
\scriptsize
\begin{multline} \label{eq:Folge}
1,\  1 + 2 x,\  1 + x + 3 {x^2},\  
\frac {1} {6}\big(6 + 31 x - 15 {x^2} + 20 {x^3}\big),\  
  \frac {1} {12}\big(12 - 58 x + 217 {x^2} - 98 {x^3} + 35 {x^4}\big),\
\\
  \frac {1} {20}\big(20 + 508 x - 925 {x^2} + 820 {x^3} - 245 {x^4} + 42 {x^5}\big),\
 \frac {1} {120}\big(120 - 8042 x 
+ 20581 {x^2} - 17380 {x^3} + 7645 {x^4} - 
      1518 {x^5} + 154 {x^6}\big),\ \\
  \frac {1} {1680}\big(1680 + 386012 x - 958048 {x^2} + 943761 {x^3} - 455455 {x^4} + 
      123123 {x^5} - 17017 {x^6} + 1144 {x^7}\big),\ \\
  \frac {1} {20160}\big(20160 - 15076944 x + 40499716 {x^2} 
- 42247940 {x^3} + 
    23174515 {x^4} - 7234136 {x^5} + 1335334 {x^6} - 134420 {x^7} + 
      6435 {x^8}\big),\ \\
  \frac {1} {181440}\big(181440 + 462101904 x - 1283316876 {x^2} + 1433031524 {x^3} - 
      853620201 {x^4} + 303063726 {x^5} \\
- 66245634 {x^6} + 
      8905416 {x^7} - 678249 {x^8} + 24310 {x^9}\big),\ \dots
\normalsize
\end{multline}
}

You may of course try to guess the coefficients of powers of $x$ in
these polynomials. But within the ``hypergeometric paradigm" this
does usually not work. In particular, that does not work with the above
sequence of polynomials.

A first very useful idea is {\em to guess through interpolation}.
(For example, this is what helped to guess coefficients in 
\machSeite{FiscAA}%
\cite{FiscAA}.)
What this means is that, for each $p_n(x)$ you try to find enough
values of $x$ for which $p_n(x)$ appears to be ``nice" (the prime
factorization of $p_n(x)$ has small prime
factors, see Footnote~\ref{foot:NICE}). Then you guess these special evaluations of
$p_n(x)$ (by, possibly, using {\tt Rate} or {\tt GUESS}), 
and, by interpolation, are able to write down a guess for
$p_n(x)$ itself.

Let us see how this works for our sequence \eqref{eq:Folge}. A few
experiments will convince you that $p_n(x)$, $0\le n\le 9$ (this is all
we have), appears to be ``nice" for
$x=0,1,\dots,n$. Furthermore, using {\tt Rate} or {\tt GUESS}, 
you will quickly find
that, apparently, $p_n(e)=\binom {2n+e}e$ for $e=0,1,\dots,n$.
Therefore, as it also appears to be the case that $p_n(x)$ is of
degree $n$, our sequence of polynomials should be given (using
Lagrange interpolation) by 
\begin{equation} \label{eq:p_n}
p_n(x)=\sum _{e=0} ^{n}\binom {2n+e}e\frac {x(x-1)\cdots
(x-e+1)(x-e-1)\cdots (x-n)} {e(e-1)\cdots
1\cdot(-1)\cdots (e-n)} .
\end{equation}

\medskip
Another useful idea is to try to {\em expand your polynomials with respect
to a ``suitable" basis}. (For example, 
this is what helped to guess coefficients
in 
\machSeite{CiKrAC}%
\cite{CiKrAC} or 
\machSeite{KrZeAA}%
\cite[e.g., (3.15), (3.33)]{KrZeAA}.) 
Now, of course, you do not know beforehand
what ``suitable" could be in your situation. 
Within the ``hypergeometric paradigm"
candidates for a suitable basis are always more or less sophisticated
shifted factorials. So,
let us suppose that we know that we were working within the
``hypergeometric paradigm" when we
came across the example \eqref{eq:Folge}. Then the simplest possible
bases are $(x)_k$, $k=0,1,\dots$, or $(-x)_k$, $k=0,1,\dots$. It
is just a matter of taste, which of these to try first. Let us try the
latter. Here are the expansions of $p_3(x)$ and $p_4(x)$ 
in terms of this basis
(actually, in terms of the equivalent basis $\binom xk$,
$k=0,1,\dots$):
\begin{align} \notag
\textstyle
\frac {1} {6}\big(6 + 31 x - 15 {x^2} + 20 {x^3}\big)
&\textstyle
=1+6\binom x1
+15\binom x2
+20\binom x3,\\
\notag
\textstyle
  \frac {1} {12}\big(12 - 58 x + 217 {x^2} - 98 {x^3} + 35
{x^4}\big)
&\textstyle
=1+8\binom x1
+28\binom x2
+56\binom x3
+70\binom x4.
\end{align}
I do not know how you feel. For me this is enough to guess that,
apparently,
$$p_n(x)=\sum _{k=0} ^{n}\binom {2n}k\binom xk.$$
(Although this is not the same expression as in \eqref{eq:p_n}, it is
identical by means of a $_3F_2$-transformation due to Thomae, see
\machSeite{GaRaAA}%
\cite[(3.1.1)]{GaRaAA}).

As was said before, we do not know beforehand what a ``suitable"
basis is. Therefore you are advised to get as much a priori
information about your polynomials as possible. For example, in
\machSeite{CiKrAA}%
\cite{CiKrAA} it was ``known" to the authors that the result which they
wanted to guess (before being able to think about a proof) is of the
form (NICE PRODUCT)${}\times{}$(IRREDUCIBLE POLYNOMIAL).
(I.e., experiments indicated that.) Moreover, they
knew that their (IRREDUCIBLE POLYNOMIAL), a polynomial in $m$,
$p_n(m)$ say, would have
the property $p_n(-m-2n+1)=p_n(m)$. Now, if we are asking ourselves
what a ``suitable" basis could be that has this property as well, and
which is built in the way of shifted factorials, then the most obvious
candidate is $(m+n-k)_{2k}=(m+n-k)(m+n-k+1)\cdots(m+n+k-1)$,
$k=0,1,\dots$. Indeed, it was very easy to guess the expansion
coefficients with respect to this basis. (See Theorems~1 and 2 in
\machSeite{CiKrAA}%
\cite{CiKrAA}. The polynomials that I was talking about are
represented by the
expression in big parentheses in 
\machSeite{CiKrAA}%
\cite[(1.2)]{CiKrAA}.)

\medskip
If the above ideas do not help, then I have nothing else to offer
than to {\em try some}, more or less arbitrary, {\em manipulations}. 
To illustrate what I could possibly mean, let us again consider
an example. In the course of working on
\machSeite{KratBD}%
\cite{KratBD}, I had to guess the result of a determinant evaluation
(which became Theorem~8 in 
\machSeite{KratBD}%
\cite{KratBD}; it is reproduced here as
Theorem~\ref{thm:TSSCPP1}). Again, the difficult part of guessing was
to guess the ``ugly" part of the result. As the dimension of the
determinant varied, this gave a certain sequence $p_n(x,y)$ of
polynomials in two variables, $x$ and $y$, of which I display
$p_4(x,y)$:
\MATH

In[1]:= VPol[4]

\leavevmode                  2      3    4                     2        3         2
\vskip-10pt%
Out[1]= 6 x + 11 x  + 6 x  + x  + 6 y - 10 x y - 6 x  y - 4 x  y + 11 y  - 
 
\leavevmode          2      2  2      3        3    4
\vskip-10pt%
>    6 x y  + 6 x  y  + 6 y  - 4 x y  + y

\endMATH
\noindent
(What I subsequently describe is the actual way in which the
expression for $p_n(x,y)$ in terms of the sum on the right-hand side
of \eqref{eq:TSSCPP1} was found.)
What caught my eyes was the part of the polynomial independent of
$y$, $x^4+6x^3+11x^2+6x$, which I recognized as
$(x)_4=x(x+1)(x+2)(x+3)$. For the fun of it, I subtracted that,
just to see what would happen:
\MATH
\goodbreakpoint%
In[2]:= Factor[\%-x(x+1)(x+2)(x+3)]
\goodbreakpoint%
\leavevmode                         2      3                     2        2        2
\vskip-10pt%
Out[2]= y (6 - 10 x - 6 x  - 4 x  + 11 y - 6 x y + 6 x  y + 6 y  - 4 x y  + 
 
\leavevmode        3
\vskip-10pt%
\MATHgroesser       y )

\endMATH
\noindent
Of course, a $y$ factors. Okay, let us cancel that:
\MATH
\goodbreakpoint%
In[3]:= \%/y
\goodbreakpoint%
\leavevmode                      2      3                     2        2        2    3
\vskip-10pt%
Out[3]= 6 - 10 x - 6 x  - 4 x  + 11 y - 6 x y + 6 x  y + 6 y  - 4 x y  + y

\endMATH
\noindent
One day I had the idea to continue in a ``systematic" manner: Let us
subtract/add an appropriate multiple of $(x)_3$\,! Perhaps,
``appropriate" in this context is to add $4(x)_3$, because that does
at least cancel the third powers of $x$:
\MATH
\goodbreakpoint%
In[4]:= Factor[\%+4x(x+1)(x+2)]
\goodbreakpoint%
\leavevmode                              2                  2
\vskip-10pt%
Out[4]= (1 + y) (6 - 2 x + 6 x  + 5 y - 4 x y + y )

\endMATH
\noindent
I assume that I do not have to comment the rest:
\MATH
\goodbreakpoint%
In[5]:= \%/(y+1)
\goodbreakpoint%
\leavevmode                     2
\vskip-10pt%
Out[5]= 6 - 2 x + 6 x  + 5 y - 4 x y + y
\goodbreakpoint%
In[6]:= Factor[\%-6x(x+1)]
\goodbreakpoint%
Out[6]= (2 + y) (3 - 4 x + y)
\goodbreakpoint%
In[7]:= \%/(y+2)
\goodbreakpoint%
Out[7]= 3 - 4 x + y
\goodbreakpoint%
In[8]:= Factor[\%+4x]
\goodbreakpoint%
Out[8]= 3 + y

\endMATH
\noindent
What this shows is that
$$p_4(x,y)=(x)_4-4(x)_3\,(y)_1+6(x)_2\,(y)_2-4(x)_1\,(y)_3+(y)_4.$$
No doubt that, at this point, you would have immediately guessed (as
I did)
that, in general, we ``must" have (compare \eqref{eq:TSSCPP1})
$$p_n(x,y)=\sum _{k=0} ^{n}(-1)^k\binom nk (x)_k\,(y)_{n-k}.$$

\section*{Appendix B: Turnbull's polarization of Bazin's theorem
implies most of the identities in Section~\ref{sec:general}}

\global\def\theequation{\mbox{B.\arabic{equation}}}
\setcounter{equation}{0}

In this appendix we show that all the determinant lemmas from
Section~\ref{sec:general}, with the exception of
Lemmas~\ref{lem:Krat6} and \ref{lem:Krat7}, follow from the
evaluation of a certain determinant of minors of a given matrix, 
an observation which I owe to Alain Lascoux. This evaluation, due to
Turnbull 
\machSeite{TurnAB}%
\cite[p.~505]{TurnAB}, is a polarized version of a theorem of
Bazin 
\machSeite{MuirAB}%
\cite[II, pp.~206--208]{MuirAB} (see also 
\machSeite{LeclAA}%
\cite[Sec.~3.1 and
3.4]{LeclAA}). 

For the statement of Turnbull's theorem we have to fix an $n$-rowed
matrix $A$, in which we label the columns,
slightly unconventionally, by
$a_2,\dots,a_m,b_{21},b_{31},b_{32},
b_{41},\dots, b_{n,n-1},\break
x_1,x_2,\dots,x_n$, 
for some $m\ge n$, i.e., $A$ is an $n\times (n+m-1+\binom
n2)$ matrix. Finally, let $[a,b,c,\dots]$ denote the minor formed by
concatenating columns $a,b,c,\dots$ of $A$, in that order.
\begin{Proposition} 
\label{prop:Turnbull}
{\em (Cf\@. 
\machSeite{TurnAB}%
\cite[p.~505]{TurnAB}, 
\machSeite{LeclAA}%
\cite[Sec.~3.4]{LeclAA})}.
With the notation as explained above, there holds
\begin{multline} \label{eq:Turnbull}
\det_{1\le i,j\le n}\big([b_{j,1},b_{j,2},\dots,b_{j,j-1},x_i,
a_{j+1},\dots,a_m]\big)\\
=[x_1,x_2,\dots,x_n,a_{n+1},\dots,a_m]
\prod _{j=2} ^{n}[b_{j,1},b_{j,2},\dots,b_{j,j-1},
a_{j},\dots,a_m].
\end{multline}
\quad \quad \qed
\end{Proposition}

Now, in order to derive Lemma~\ref{lem:Krat1} from
\eqref{eq:Turnbull}, we choose $m=n$ and for $A$ the matrix
\def\tbox#1#2{\hbox to#1pt{\hfil\scriptsize$#2$\hfil}}

{\scriptsize
\begin{gather*} 
\setcounter{MaxMatrixCols}{12}
\begin{matrix} 
\tbox{40.05382}{a_2}\hskip-2pt&\hskip-2pt
\tbox{12.33325}{\dots}\hskip-2pt&\hskip-2pt
\tbox{41.00523}{a_n}\hskip-2pt&\hskip-2pt
\tbox{40.11719}{b_{21}}\hskip-2pt&\hskip-2pt
\tbox{40.11719}{b_{31}}\hskip-2pt&\hskip-2pt
\tbox{40.11719}{b_{32}}\hskip-2pt&\hskip-2pt
\tbox{12.33325}{\dots}\hskip-2pt&\hskip-2pt
\tbox{41.0686}{b_{n,n-1}}\hskip-2pt&\hskip-2pt
\tbox{22.90746}{x_1}\hskip-2pt&\hskip-2pt
\tbox{22.90746}{x_2}\hskip-2pt&\hskip-2pt
\tbox{12.33325}{\dots}\hskip-2pt&\hskip-2pt
\tbox{22.90746}{x_n}
\end{matrix}
\hphantom{{},{}}
\\
\begin{pmatrix} 
1\hskip-2pt&\hskip-2pt
\dots\hskip-2pt&\hskip-2pt
1\hskip-2pt&\hskip-2pt
1\hskip-2pt&\hskip-2pt
1\hskip-2pt&\hskip-2pt
1\hskip-2pt&\hskip-2pt
\dots\hskip-2pt&\hskip-2pt
1\hskip-2pt&\hskip-2pt
1\hskip-2pt&\hskip-2pt
1\hskip-2pt&\hskip-2pt
\dots\hskip-2pt&\hskip-2pt
1\\
-A_2\hskip-2pt&\hskip-2pt
\dots\hskip-2pt&\hskip-2pt
-A_n\hskip-2pt&\hskip-2pt
-B_{2}\hskip-2pt&\hskip-2pt
-B_{2}\hskip-2pt&\hskip-2pt
-B_{3}\hskip-2pt&\hskip-2pt
\dots\hskip-2pt&\hskip-2pt
-B_{n}\hskip-2pt&\hskip-2pt
X_1\hskip-2pt&\hskip-2pt
X_2\hskip-2pt&\hskip-2pt
\dots\hskip-2pt&\hskip-2pt
X_n\\
(-A_2)^2\hskip-2pt&\hskip-2pt
\dots\hskip-2pt&\hskip-2pt
(-A_n)^2\hskip-2pt&\hskip-2pt
(-B_{2})^2\hskip-2pt&\hskip-2pt
(-B_{2})^2\hskip-2pt&\hskip-2pt
(-B_{3})^2\hskip-2pt&\hskip-2pt
\dots\hskip-2pt&\hskip-2pt
(-B_{n})^2\hskip-2pt&\hskip-2pt
X_1^2\hskip-2pt&\hskip-2pt
X_2^2\hskip-2pt&\hskip-2pt
\dots\hskip-2pt&\hskip-2pt
X_n^2\\
\hdotsfor{12}\\
(-A_2)^{n-1}\hskip-2pt&\hskip-2pt
\dots\hskip-2pt&\hskip-2pt
(-A_n)^{n-1}\hskip-2pt&\hskip-2pt
(-B_{2})^{n-1}\hskip-2pt&\hskip-2pt
(-B_{2})^{n-1}\hskip-2pt&\hskip-2pt
(-B_{3})^{n-1}\hskip-2pt&\hskip-2pt
\dots\hskip-2pt&\hskip-2pt
(-B_{n})^{n-1}\hskip-2pt&\hskip-2pt
X_1^{n-1}\hskip-2pt&\hskip-2pt
X_2^{n-1}\hskip-2pt&\hskip-2pt
\dots\hskip-2pt&\hskip-2pt
X_n^{n-1}
\end{pmatrix},
\end{gather*}
}

\noindent
with the unconventional labelling of the columns indicated on
top. I.e., column $b_{st}$ is filled with powers of $-B_{t+1}$, $1\le
t<s\le n$.
With this choice of $A$, all the minors in \eqref{eq:Turnbull}
are Vandermonde determinants. In particular, due to the Vandermonde
determinant evaluation \eqref{eq:Vandermonde}, we then have for the
$(i,j)$-entry of the determinant in \eqref{eq:Turnbull}
\begin{multline*} 
[b_{j,1},b_{j,2},\dots,b_{j,j-1},x_i,
a_{j+1},\dots,a_m]\\
=\prod _{2\le s<t\le j} ^{}(B_s-B_t)
\prod _{j+1\le s<t\le n} ^{}(A_s-A_t)
\prod _{s=2} ^{j} \prod _{t=j+1} ^{n}(A_t-B_s)\\
\times
\prod _{s=j+1} ^{n}(X_i+A_s)\prod _{s=2} ^{j}(X_i+B_s),
\end{multline*}
which is, up to factors that only depend on the column index $j$,
exactly the $(i,j)$-entry of the determinant in \eqref{eq:Krat1}.
Thus, Turnbull's identity \eqref{eq:Turnbull} gives the evaluation
\eqref{eq:Krat1} immediately, after some obvious simplification.

In order to derive Lemma~\ref{lem:Krat3} from
\eqref{eq:Turnbull}, we choose $m=n$ and for $A$ the matrix 

{\scriptsize
\begin{gather*} 
\begin{matrix} 
\tbox{69.44788}{a_2}\hskip-3pt&\hskip-3pt
\tbox{12.33325}{\dots}\hskip-3pt&\hskip-3pt
\tbox{71.35071}{a_n}\hskip-3pt&\hskip-3pt
\tbox{81.97723}{b_{21}}\hskip-3pt&\hskip-3pt
\tbox{81.97723}{b_{31}}
\end{matrix}
\hskip2cm
\\
\left(\begin{matrix} 
1\hskip-2pt&\hskip-2pt
\dots\hskip-2pt&\hskip-2pt
1\hskip-2pt&\hskip-2pt
1\hskip-2pt&\hskip-2pt
1
\\
-A_2-C/A_2\hskip-2pt&\hskip-2pt
\dots\hskip-2pt&\hskip-2pt
-A_n-C/A_n\hskip-2pt&\hskip-2pt
-B_{2,1}-C/B_{2,1}\hskip-2pt&\hskip-2pt
-B_{3,1}-C/B_{3,1}
\\
(-A_2-C/A_2)^2\hskip-2pt&\hskip-2pt
\dots\hskip-2pt&\hskip-2pt
(-A_n-C/A_n)^2\hskip-2pt&\hskip-2pt
(-B_{2,1}-C/B_{2,1})^2\hskip-2pt&\hskip-2pt
(-B_{3,1}-C/B_{3,1})^2
\\
\hdotsfor{5}\\
(-A_2-C/A_2)^{n-1}\hskip-2pt&\hskip-2pt
\dots\hskip-2pt&\hskip-2pt
(-A_n-C/A_n)^{n-1}\hskip-2pt&\hskip-2pt
(-B_{2,1}-C/B_{2,1})^{n-1}\hskip-2pt&\hskip-2pt
(-B_{3,1}-C/B_{3,1})^{n-1}
\end{matrix}\right.
\hskip2cm
\end{gather*}
\begin{gather*}
\hskip1.5cm
\begin{matrix} 
\tbox{81.97723}{b_{32}}\hskip-3pt&\hskip-3pt
\tbox{12.33325}{\dots}\hskip-3pt&\hskip-3pt
\tbox{120.67157}{b_{n,n-1}}\hskip-3pt&\hskip-3pt
\tbox{64.17642}{x_1}\hskip-3pt&\hskip-3pt
\tbox{12.33325}{\dots}\hskip-3pt&\hskip-3pt
\tbox{66.07925}{x_n}
\end{matrix}
\hphantom{{},{}}
\\
\hskip1.5cm
\left.\begin{matrix} 
1\hskip-2pt&\hskip-2pt
\dots\hskip-2pt&\hskip-2pt
1\hskip-2pt&\hskip-2pt
1\hskip-2pt&\hskip-2pt
\dots\hskip-2pt&\hskip-2pt
1
\\
-B_{3,2}-C/B_{3,2}\hskip-2pt&\hskip-2pt
\dots\hskip-2pt&\hskip-2pt
-B_{n,n-1}-C/B_{n,n-1}\hskip-2pt&\hskip-2pt
X_1+C/X_1\hskip-2pt&\hskip-2pt
\dots\hskip-2pt&\hskip-2pt
X_n+C/X_n
\\
(-B_{3,2}-C/B_{3,2})^2\hskip-2pt&\hskip-2pt
\dots\hskip-2pt&\hskip-2pt
(-B_{n,n-1}-C/B_{n,n-1})^2\hskip-2pt&\hskip-2pt
(X_1+C/X_1)^2\hskip-2pt&\hskip-2pt
\dots\hskip-2pt&\hskip-2pt
(X_n+C/X_n)^2
\\
\hdotsfor{6}
\\
(-B_{3,2}-C/B_{3,2})^{n-1}\hskip-2pt&\hskip-2pt
\dots\hskip-2pt&\hskip-2pt
(-B_{n,n-1}-C/B_{n,n-1})^{n-1}\hskip-2pt&\hskip-2pt
(X_1+C/X_1)^{n-1}\hskip-2pt&\hskip-2pt
\dots\hskip-2pt&\hskip-2pt
(X_n+C/X_n)^{n-1}
\end{matrix}\right).
\end{gather*}

}

\noindent
(In this display, 
the first line contains columns $a_2,\dots,b_{31}$ of $A$,
while the second line contains the remaining columns.)
Again, with this choice of $A$, all the minors in \eqref{eq:Turnbull}
are Vandermonde determinants.
Therefore, by noting that $(S+C/S)-(T+C/T)=(S-T)(C/S-T)/(-T)$, and 
by writing $p_{j-1}(X)$ for
\begin{equation} \label{eq:pj}
\prod _{s=1} ^{j-1}(X+B_{j,s})(C/X+B_{j,s}),
\end{equation}
we have for the $(i,j)$-entry of the determinant in
\eqref{eq:Turnbull}
\begin{multline*} 
[b_{j,1},b_{j,2},\dots,b_{j,j-1},x_i,
a_{j+1},\dots,a_m]
=\prod _{1\le s<t\le j-1}
^{}(B_{j,s}+C/B_{j,s}-B_{j,t}-C/B_{j,t})\\
\times
\prod _{j+1\le s<t\le n} ^{}(A_s+C/A_s-A_t-C/A_t)
\prod _{s=1} ^{j-1} \prod _{t=j+1} ^{n}(A_t+C/A_t-B_{j,s}-C/B_{j,s})\\
\times
\bigg(\prod _{s=j+1} ^{n}(X_i+A_s)(C/X_i+A_s)\,A_s^{-1}\bigg)
p_{j-1}(X_i)\prod _{s=1} ^{j-1}B_{j,s}^{-1}
\end{multline*}
for the $(i,j)$-entry of the determinant in \eqref{eq:Turnbull}.
This is, up to factors which depend only on the column index $j$, 
exactly the $(i,j)$-entry of the determinant in \eqref{eq:Krat3}.
The polynomials $p_{j-1}(X)$, $j=1,2,\dots,n$, 
can indeed be regarded as arbitrary
Laurent polynomials satisfying the conditions of Lemma~\ref{lem:Krat3},
because any Laurent polynomial $q_{j-1}(X)$ over the complex numbers 
of degree at most $j-1$
and with $q_{j-1}(X)=q_{j-1}(C/X)$ can be written in the form
\eqref{eq:pj}.
Thus, Turnbull's identity \eqref{eq:Turnbull} implies the evaluation
\eqref{eq:Krat3} as well.

Similar choices for $A$ are possible in order to derive
Lemmas~\ref{lem:Krat2}, \ref{lem:Krat3a} and \ref{lem:Krat5} 
(which are in fact just limiting cases of Lemma~\ref{lem:Krat3}) from
Proposition~\ref{prop:Turnbull}.

\section*{Appendix C: Jean-Yves Thibon's proof of Theorem~\ref{thm:Zagiermaj}}

\global\def\theequation{\mbox{C.\arabic{equation}}}
\setcounter{equation}{0}

Obviously, the determinant in \eqref{eq:Zagiermaj} is the
determinant of the linear operator\break 
$K_n(q):=\sum _{\si\in \frak S_n}
^{}q^{\maj \si}\si$ acting on the group algebra $\mathbb C[\frak S_n]$ of the
symmetric group. Thus, if we are able to determine all the
eigenvalues of this operator, together with their multiplicities, we
will be done. The determinant is then just the product of all the
eigenvalues (with multiplicities).

The operator $K_n(q)$ is also an element of
Solomon's descent algebra (because permutations with the same descent
set must necessarily have the same major index). The descent algebra
is canonically isomorphic to the algebra of noncommutative symmetric
functions (see 
\machSeite{GeKLAA}%
\cite[Sec.~5]{GeKLAA}). It is
shown in 
\machSeite{KrLTAA}%
\cite[Prop.~6.3]{KrLTAA} that, as a noncommutative symmetric function,
$K_n(q)$ is equal to $(q;q)_n\, S_n(A/(1-q))$, where $S_n(B)$ 
denotes the complete (noncommutative) symmetric function of degree $n$
of some alphabet $B$. 

The inverse element of $S_n(A/(1-q))$
happens to be $S_n((1-q)A)$, i.e., 
$S_n((1-q)A)*S_n(A/(1-q))=S_n(A)$,\footnote{By 
definition of the isomorphism between noncommutative symmetric
functions and elements in the descent algebra, 
$S_n(A)$ corresponds to the identity element in the
descent algebra of $\frak S_n$.}
with $*$ denoting the internal multiplication of noncommutative
symmetric functions (corresponding to the multiplication in the
descent algebra).
This is seen as follows. As in
\machSeite{KrLTAA}%
\cite[Sec.~2.1]{KrLTAA} let us write
$\si(B;t)=\sum _{n\ge0} ^{}S_n(B)t^n$
for the generating function for complete symmetric functions of some
alphabet $B$, and 
$\la(B;t)=\sum _{n\ge0} ^{}\La_n(B)t^n$
for the generating function for elementary symmetric functions, which
are related by $\la(B;-t)\si(B;t)=1$.
Then, by 
\machSeite{KrLTAA}%
\cite[Def.~4.7 and Prop.~4.15]{KrLTAA}, 
we have $\sigma((1-q)B;1)=\lambda(B;-q)\sigma(B;1)$.
Let $X$ be the ordered alphabet $\cdots < q^2 < q < 1$, so that
$XA=A/(1-q)$.
According to 
\machSeite{KrLTAA}%
\cite[Theorem~4.17]{KrLTAA}, it then follows that
\begin{multline*} 
\sigma((1-q)A;1)*\sigma(XA;1)=
\sigma((1-q)XA;1)=\lambda(XA;-q)\sigma(XA;1)\\
= \lambda(XA;-q)\sigma(XA;q)\sigma(A;1) =\sigma(A;1),
\end{multline*}
since by definition of $X$, $\sigma(XA;1)$ is equal to 
$\sigma(XA;q)\sigma(A;1)$ (see 
\machSeite{KrLTAA}%
\cite[Def.~6.1]{KrLTAA}).
Therefore, $S_n((1-q)A)*S_n(XA)=S_n(A)$, as required.

Hence, we infer that $K_n(q)$ is the inverse of
$S_n((1-q)A)/(q;q)_n$.

The eigenvalues of $S_n((1-q)A)$ are given in 
\machSeite{KrLTAA}%
\cite[Lemma~5.13]{KrLTAA}. Their multiplicities 
follow from a combination of Theorem~5.14 and
Theorem~3.24 in
\machSeite{KrLTAA}%
\cite{KrLTAA}, since the construction in Sec.~3.4 of
\machSeite{KrLTAA}%
\cite{KrLTAA} yields idempotents $e_\mu$
such that the commutative immage of $\alpha(e_\mu)$ is equal to
$p_\mu/z_\mu$. Explicitly, the eigenvalues of $S_n((1-q)A)$ are
$\prod _{i\ge1} ^{}(1-q^{\mu_i})$, where $\mu=(\mu_1,\mu_2,\dots)$
varies through all partitions of $n$, with corresponding multiplicities
$n!/z_\mu$, the number of
permutations of cycle type $\mu$, i.e., 
$z_\mu=1^{m_1}m_1!\,2^{m_2}m_2!\cdots$, where $m_i$ is the number of
occurences of $i$ in the partition $\mu$, $i=1,2,\dots$. 
Hence, the eigenvalues of
$K_n(q)$ are $(q;q)_n/\prod _{i\ge1} ^{}(1-q^{\mu_i})$, with the same
multiplicities. 

Knowing all the eigenvalues of $K_n(q)$ and their multiplicities
explicitly, it is now not extremely difficult to form the product of
all these and, after a short calculation, recover the right-hand side
of \eqref{eq:Zagiermaj}.

\section*{Acknowledgements}
I wish to thank an anonymous referee, 
Joris Van der Jeugt, Bernard Leclerc, Madan Lal Mehta,
Alf van der Poorten, Volker
Strehl, Jean-Yves Thibon, Alexander Varchenko, and especially Alain
Lascoux, for the many useful comments and discussions which helped to
improve the contents of this paper considerably.

\newpage

\immediate\closeout\Seiten

\input detsurv.sei

\message{!!!!!!!!!!!!!!!!!!!!!!!!!!!!!!!!!!!!!!!!!!!!!!!!}
\message{Achtung! Die Zitate von Muir in Schendel und Voigt}
\message{muessen immer angepasst werden!}
\message{!!!!!!!!!!!!!!!!!!!!!!!!!!!!!!!!!!!!!!!!!!!!!!!!}

\end{document}